\documentclass[12pt,a4paper,twoside]{book}

\usepackage{url}
\usepackage{hyperref}
\usepackage{graphicx}
\usepackage{epsfig}
\usepackage{fancyhdr,fancybox}
\usepackage{amsfonts}
\usepackage{amsmath}
\usepackage{amssymb}
\usepackage{verbatim}
\usepackage{amsfonts}
\usepackage{subfigure}
\usepackage[german,english]{babel}
\usepackage{latexsym}
\usepackage{multicol}
\usepackage{multirow}
\usepackage{CJK}

\newcommand{\R}{\mathbb{R}}
\newcommand{\Rd}{\mathbb R^d}
\newcommand{\Hs}{\mathcal H^s}
\newcommand{\N}{\mathbb{N}}

\newcommand{\eps}{\varepsilon}

\newcommand{\Fe}{F_\epsilon}
\newcommand{\Sd}{\mathcal{Q}^d}
\newcommand{\ti}{{\bf Fractal Curvature Measures \\and Image Analysis}
}

\newcommand{\tisimple}{{\bf Fractal Curvature Measures and Image Analysis}}

\newcounter{global}[chapter]

\newcounter{eqcounter}[chapter]

\newcommand{\supp}{{\rm supp}}
\newcommand{\bd}{{\rm bd}\;}

\newcommand{\be}{\begin{equation}}
\newcommand{\ee}{\end{equation}}
\newcommand{\bea}{\begin{eqnarray}}
\newcommand{\eea}{\end{eqnarray}}
\newcommand{\beastar}{\begin{eqnarray*}}
\newcommand{\eeastar}{\end{eqnarray*}}
\newcommand{\bee}{\begin{equation*}}
\newcommand{\eee}{\end{equation*}}
\newcommand{\bs}{\begin{split}}
\newcommand{\es}{\end{split}}
\newcommand{\br}{\begin{rm}}
\newcommand{\er}{\end{rm}}

\newlength{\MyNumberwidth}
\newlength{\MyHeaderWidth}

\newtheorem{theorem}{Theorem}[section]
\newtheorem{lemma}[theorem]{Lemma}
\newtheorem{lemma*}{Lemma}
\newtheorem{proposition}[theorem]{Proposition}
\newtheorem{corollary}[theorem]{Corollary}
\newtheorem{definition}[theorem]{Definition}
\newtheorem{definition*}[lemma*]{Definition}

\begin{document}


\thispagestyle{empty}

\begin{center}
{\Large \bf Diplomarbeit} \\
\vspace{1cm}
{zur Erlangung des akademischen Grades Diplom-Mathematiker \\
an der Fakult\"{a}t f\"{u}r Mathematik und Wirtschaftswissenschaften der \\
Universit\"{a}t Ulm\\} \vspace{1.5cm}
\begin{figure}[h]
\centering
\includegraphics[width=5cm]{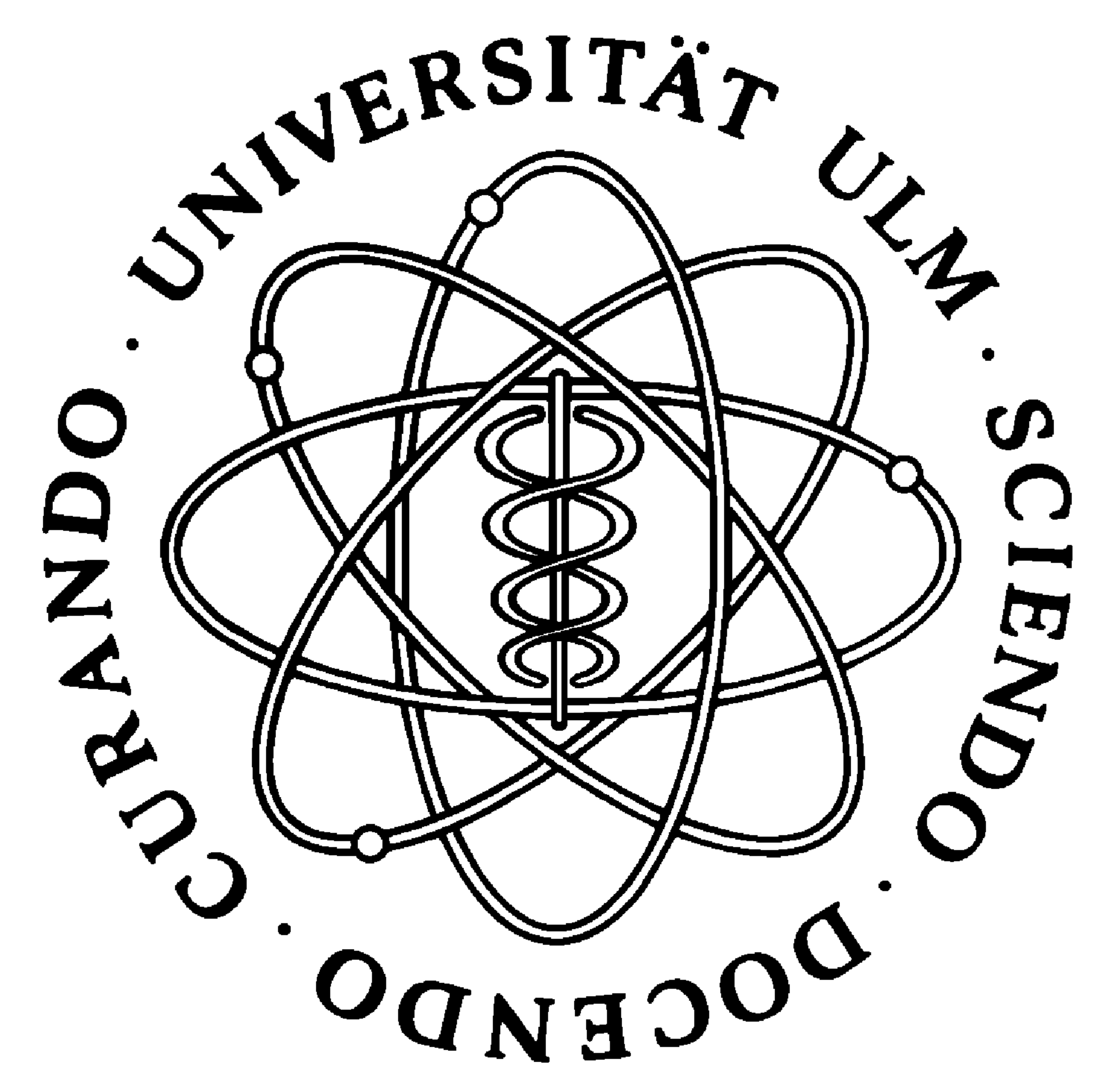} 
\end{figure}
\vspace{1cm} {\huge \ti \\} \vspace{1.5cm} {\large vorgelegt von}

{\large Peter Straka}\\

\vspace{1.5cm}

{\large {\bf Gutachter:} Prof.~Dr.~E.~Spodarev
, Prof.~Dr.~W.~Arendt}

\vspace{1.5cm}

Universit\"{a}t Ulm, Abteilung Stochastik

Oktober 2007

\end{center}
\clearpage


\thispagestyle{empty}


\clearpage

\thispagestyle{empty} 

\begin{minipage}[t]{15cm}
 \vspace*{7cm}
\begin{center}
{\large \textbf{Abstract}}
\end{center}

Since the recent dissertation by Steffen Winter, for certain self-similar sets $F$ the growth behaviour of the Minkowski functionals of the parallel sets $F_\eps := \{x\in \Rd : d(x,F)\leq \eps\}$ as $\eps \downarrow 0$ is known,  leading to the notion of fractal curvatures $C_k^f(F)$, $k\in \{0,\ldots,d\}$. The dependence of the growth behaviour on the fractal dimension $s = \dim F$ is exploited, and estimators for $s$ and $C_k^f(F)$ are derived. The performance of these estimators is tested on binary images of self-similar sets.

\end{minipage}


\thispagestyle{empty}

\hspace{1cm}
\begin{minipage}[t]{12cm}
{\LARGE {\itshape Danksagung}} \vspace*{1.5cm}

\selectlanguage{german}

Meinen tiefsten Dank m"ochte ich meinen Eltern aussprechen, die mich w"ahrend meiner langen Studienzeit immer bedingungslos unterst"utzt haben. Ohne sie w"are diese Arbeit nicht zustandegekommen. 

Prof. Dr. Evgeny Spodarev m"ochte ich f"ur dieses sch"one Thema danken, das es auch einem Neuling in der fraktalen Geometrie erm"oglicht hat innovative Arbeit zu leisten. Auch danke ich ihm f"ur seine stete Hilfsbereitschaft w"ahrend meiner gesamten Zeit an der Universit"at Ulm. 

Obwohl Dr. Steffen Winter keine offizielle Betreuerfunktion inne hatte, war er "uber Email und Telefon eine allgegenw"artige und prompte Hilfe bei meinen zahlreichen geometrischen Fragen. Ihm danke ich auch insbesondere f"ur die gewissenhafte Durchsicht dieser Arbeit. 

Danke an Daniel Meschenmoser und Johannes Mayer f"ur ihre Hilfe bei Programmierfragen und allgemeinen Fragen zum GeoStoch-Projekt. Besonders m"ochte ich Daniel danken f"ur die prompte Bereitstellung seiner Java-klasse ``Cell-Utils'' und f"ur das geduldige Ertragen meines Telefon-terrors.

Tobias Wunner m"ochte ich f"ur alle Zeitersparnis danken die mir seine Matlab-Kenntnisse boten, sowie f"ur sein Talent diese Zeit gleich wieder in ablenkendere Dinge zu investieren.

\selectlanguage{english}

Thanks also go to Dr. Simon Lloyd for his supporting me with food throughout the last phase of this thesis.

%
%
%
%

\vspace*{2cm}
\end{minipage}


\thispagestyle{empty}
\pagenumbering{roman} 
\lhead{Contents} \rhead{\thepage} \cfoot{} \tableofcontents
\clearpage 

\thispagestyle{empty} 
\addcontentsline{toc}{section}{List of Figures}
\lhead{List of Figur	es} \rhead{\thepage}
\cfoot{} \listoffigures 
\clearpage

\thispagestyle{empty}
\addcontentsline{toc}{section}{List of Tables}
\lhead{List of Tables} \rhead{\thepage}
\cfoot{} \listoftables 
\clearpage


\pagenumbering{arabic} 
\chapter*{Introduction}
\addcontentsline{toc}{chapter}{Introduction}
\pagestyle{myheadings}
\markboth{INTRODUCTION}{INTRODUCTION}

At the end of the 19th century, the common belief that every subset $A \subset \Rd$ of euclidean space could be assigned an integer number called ``the dimension'' of $A$ could not be supported anymore:
Examples like the Cantor-set (termed after the German mathematician Georg Cantor) showed that geometric measure theory needed to be generalized, which was then done essentially by Felix Hausdorff. 
The theory around these abstract sets, however, did not affect the natural sciences until some 50 years later, when strange attractors were discovered. Attractors are subsets of $\mathbb R^3$ with the property that there is a dynamical system in which trajectories of particles converge to this set, and they are termed ``strange'' if their dimension is non-integer (fractal). 

This development continued with a general reconsideration of how accurately length and area of concrete objects can be determined; a good example is the question ``How long is the coast of Britain?'' (B. Mandelbrot). A commonly expected answer would be a 4 digit number of miles, which corresponds to the length of a rectifiable curve that represents the coastline. Recall that measuring means adding up the lengths of the steps a divider has to take to run through the curve. The reciprocal value of the steplenghts is usually termed ``accuracy''. As the accuracy tends to $\infty$, the measured length will stay finite for rectifiable curves, and approach $\infty$ for non-rectifiable curves. Now the coast of Britain cannot ever be measured with arbitrarily high accuracy, but extrapolation to infinity \textit{suggests} that the measured lengths might rather tend to $\infty$ than to a finite value\footnote{
To visualize this, compare the distances i)~an airplane has to travel to view all of the coastline of Britain 
  ii)~a human has to walk to circle Britain with his feet touching the water and
  iii)~an ant has to travel if it is to cross every pebble of Britain's shore which touches the ocean (frozen in time).
}. An answer to whether or not this is right is impossible, of course, since the concrete coastline is something intrinsically different than an abstract mathematical curve, be it rectifiable or not. 

Note that this example can easily be taken into two dimensions, too: Consider a 1 cm$^2$ square piece of smooth glass. Under an electron microscope a very irregular shape is revealed, the measurement of the area of which on this scale is yielding a value far greater than 1 cm$^2$.

The observations above all have in common that measurements of areas and lengths grow as the scale shrinks. Fractal (not necessarily integer-valued) dimensions turn out to describe the rate of this growth; conversely, the knowledge of the growth-rate translates into knowledge of the fractal dimension of a set. To the best of the author's knowledge, so far the fractal character of subsets of $\mathbb R^3$ has only been described based on (the growth behaviour of) measurements of \textit{one} physical quantity only, i.e. either of volume, area or of length exclusively. However, for a big class of subsets $A$ of $\mathbb R^3$ one has the four important characteristics termed ``Minkowski-functionals'': These are denoted by $C_3(A)$, $C_2(A)$, $C_1(A)$ and $C_0(A)$ and have the physical interpretation \textit{volume}, \textit{surface area}, \textit{mean breadth} and \textit{Euler-characteristic} (i.e. the number of connected components minus the number of holes) respectively. The reason why, so far, fractal properties have been analysed on the basis of $C_d(A)$ only, must be the following: Because of the complexity of fractals in dimensions that are higher than one, nothing was known about the behaviour of the characteristics $C_k(A)$ for $k \leq d-1$ as the scale shrinks, until very recently S.Winter made the first step into this direction: For a subclass of the class of self-similar sets, which are probably the best-understood fractals, he thoroughly described the limiting behaviour of these characteristics. Now that this behaviour is known for a (yet very small) class of fractals, the theoretical foundation is laid for the analysis of fractals, \textit{also} based on $C_{d-1}(A),\ldots,C_0(A)$. Thus the question that this thesis is devoted to can now be formulated, here in the most abstract way possible:

\begin{quotation}
\textbf{If, on a shrinking scale, the simultaneous growth behaviour of \textit{all} Minkowski-functionals is known, how can this additional information be translated into further knowledge about geometrical properties of the fractal set?}
\end{quotation} 

Of course, one can only hope to find partial answers to this question, in a very restricted setting. We are only going to consider the class $\mathcal Q^2$ of sets which are subsets of $\mathbb R^2$, self-similar, satisfy the open set condition and have polyconvex parallel sets. By the $r$-parallel set of $A$ we mean the set $A_r$ of all points whose distance to $A$ is $\leq r$. The value $r$ reflects the accuracy of a measurement of $A$: If a set $A \in \mathcal{Q}^2$ existed in nature, (it would most probably be invisible and) measurements of the locations of its points would have an incertainty of $\pm r$ in each direction, so real-world measurements of the Minkowski-functionals $C_k(A)$ can be compared to those of $C_k(A_r)$. As the accuracy tends to infinity, i.e. as $r\downarrow 0$, Winter has shown that $\Sd$-sets would exhibit the growth behaviour $C_k(A_r) \sim C_k^f r^{s-k}$ (where $s = \dim A$), and that there is some analogy of the prefactors $C^f_k$ to the classical curvatures. The scope of this thesis will be to estimate the values $C^f_k$ and $s$ from given data $C_k(A_r)$. The construction of $A$, the dilation of $A$ by $r$ and the measurement of $C_k(A_r)$ were all performed by computer simulations, where the author strongly benefited from algorithms that were readily implemented in the GeoStoch (\cite{geostoch}) library of the University of Ulm.

This thesis is organized as follows: Part~I describes basic facts of fractal geometry in chapter~\ref{ch:fractal-geometry} and excerpts of the theory behind fractal curvature in chapter~\ref{ch:fractal_curvature}, which is mostly a summarization of theorems taken from Winter's dissertation \cite{diss_winter} that are necessary for chapter~\ref{ch:main}. In part~II, chapter~\ref{ch:discreteness}, a short introduction to the problem of squeezing fractals into binary images is given, followed by chapter~\ref{ch:review} which contains a review of methods that analyse fractal images. The author has found a great amount of literature on estimating fractal dimension, but nothing on estimates of the minkowski content, which would have been interesting with respect to chapter~\ref{ch:main}. Literature on the estimation of the dimension of graphs of functions is vast but has been ignored because it fits more in the scope of time series than image analysis. The main chapter, chapter~\ref{ch:main}, introduces a joint estimator of fractal dimension and of the three fractal curvature measures for sets in $\mathcal{Q}^2$. Finally, in chapter~\ref{ch:results}, the performance of the estimator is tested on several sample images.

\pagestyle{headings}

\part{Fractal Geometry}

\chapter{Basics of Fractal Geometry}
\label{ch:fractal-geometry}

In this chapter, notions from fractal geometry that are important for later use shall be summarized, and important theorems are cited. Proofs can be found in the books of Falconer~\cite{falconer90}, Mattila \cite{mattila95} and Rogers \cite{rogers99}.

\section{Hausdorff-measure and Fractal Dimension of Sets}

The maybe most elementary characteristic that is classifying a geometric object is its dimension. In classical geometry, the \textit{topological dimension} was most widely used: Wherever there exists a local isomorphism (a bijective mapping, both directions being continuous) from an open subset of $\mathbb R^k$ ($k \in \mathbb N_0$) to the object, the object is said to have dimension $k$. 

In 1918, \textit{Hausdorff} generalized this notion of dimension in such a way that a dimension could be assigned to \textit{every} Borel-set $F \in \mathbb R^d$, while the following demands were still met:
\begin{itemize}
	\item Monotonicity: $E \subset F \Rightarrow \dim(E) \leq \dim(F)$
	\item Open sets: If $E$ is a non-empty open subset of $\Rd$, then $\dim(E) = d$.
	\item Geometric invariance: If $f$ is a (non-singular) affine transformation, then 

		$\dim(f(E)) = \dim(E)$.
	\item Stability under countable unions: 
		$\dim\left(\bigcup_{i \in \mathbb N} E_i \right) = \sup_{i \in \mathbb N} \dim (E_i). $
\end{itemize}
A set whose Hausdorff-dimension is non-integer or different from its topological dimension is commonly termed a \textit{fractal}, though no commonly used mathematical definition of the word ``fractal'' exists.

We are now going to give a short introduction to Hausdorff-measures, which are the foundation of the Hausdorff-dimension, and we only consider the special case of the metric space $X=\mathbb R^d$ with the euclidean metric $d(\cdot,\cdot)$. First we introduce a class of \textit{dimension functions}:
\begin{equation*}
 \mathbb H:=\{h:[0,\infty]\rightarrow[0,\infty): h 
 \text{ is right-continuous and non-decreasing, } 
 h(t)=0\Leftrightarrow t=0 \}
\end{equation*}
Now for $\delta >0$, $A\subset X$, $h \in \mathbb H$ and 
$|A|=\sup\{d(x,y), x,y \in A\}$ being the diameter of $A$, let
\begin{equation}
 \mathcal H^h_\delta(A):=\inf \left\lbrace 
 \sum_{i=1}^\infty h(|C_i|): A\subset \bigcup_{i=1}^\infty {C_i}, 
 |C_i|<\delta
 \right\rbrace 
\label{eqn-Hhdelta}
\end{equation}
where the infimum is taken over all countable coverings $\{C_i\}$ of $A$ with arbitrary subsets $C_i \subset X$. $\mathcal H^h_\delta$ is not necessarily a Borel-measure, but 
\begin{equation*}
 \mathcal H^h(A):=\lim_{\delta\rightarrow 0} \mathcal H^h_\delta(A)
\end{equation*}
is. It is given the name \textit{Hausdorff-measure with dimension function h}. If 
$0 < \mathcal H^h(A) < \infty$ then $h$ is called the \textit{exact dimension function} of $A$. 

Consider now the subclass 
\[ \mathbb H_0 = \{h(t): h(t) = t^\alpha, \alpha \in (0,\infty)\}\]
of $\mathbb H$. If $h(t) = t^\alpha \in \mathbb H_0$, it is common to put $\mathcal H^\alpha := \mathcal H^h$ and to call it the $\alpha$-dimensional Hausdorff-measure. The reason for this is that for positive integers $d$, $\mathcal H^d$ and the $d$-dimensional Lebesgue measure $\lambda_d$ are multiples of each other. 

For $\alpha \leq \beta$ we have $\mathcal H^\beta_\delta(A) \leq \delta^{\beta - \alpha} \mathcal H^\alpha_\delta(A)$, which yields the implications
\begin{eqnarray*}
 \mathcal H^\alpha(A)< \infty &\Longrightarrow&  \mathcal H^{\alpha+\eps}(A)=0, \\
 \mathcal H^\alpha(A)>0 &\Longrightarrow&  \mathcal H^{\alpha-\eps}(A) = \infty
\end{eqnarray*}
whenever $0<\eps<\alpha$.
An immediate consequence is that given a subset $A \subset \mathbb R^d$, there is one distinct value $s$ such that 
\begin{equation*}
 \mathcal H^\alpha(A)= \begin{cases} 
\infty  & \text{ for }\alpha < s \\
0  & \text{ for }\alpha > s                             \end{cases}
\end{equation*}
We write $\dim_H(A)=s$ and say that \textit{A has Hausdorff-dimension} $s$. Note that the three cases $\mathcal H^s(A) = 0$, $0 < \mathcal H^s(A) < \infty$, $\mathcal H^s(A)=\infty$ are possible. The second case is of special interest, and we shall call $A \subset \mathbb R^d$ an $s$-set if $0 < \mathcal H^s(A) < \infty$, i.e. if $h(t)=t^s$ is an exact dimension function of $A$. In a later chapter we will see that self-similar sets satisfying the open set condition (OSC) are such sets. 

There are examples of sets $A$ for which $t^s$ is not an exact dimension function for any $s>0$, as almost all Brownian paths $\Xi_T$. In fact, $\dim_H (\Xi_T) = 2$ but $\mathcal H^2(\Xi_T) = 0$ almost surely; an exact dimension function would be e.g. $t^2\log\log(\frac{1}{t})$.

The Hausdorff dimension is very useful in a mathematical framework, but it
is not easy to calculate directly: The coverings from equation~\ref{eqn-Hhdelta} consist of possibly infinitely many sets of varying diameter.
Experimentalists prefer the following notion of dimension, as it is more easily calculated:
\begin{definition}
\label{def-1}
For $\delta>0$ and a bounded set $A$, let $N_\delta(A)$ be the smallest number of sets of diameter not greater than $\delta$ whose union covers $A$. The numbers 
\[
\underline\dim_B(A):= \liminf_{\delta\rightarrow 0+}\frac{\log N_\delta}{-\log \delta} \hspace{0.5cm} and \hspace{0.5cm}
\overline\dim_B(A):= \limsup_{\delta\rightarrow 0+}\frac{\log N_\delta}{-\log \delta}, 
\]
are called \textit{lower} and \textit{upper} Box-counting-dimensions (or lower and upper Minkowski-dimensions) of $A$. If they coincide, their common value $\dim_B(A)$ is simply called the Box-counting dimension (or Minkowski-dimension) of $A$. 
\end{definition}
Box-counting dimensions satisfy $\overline\dim_B(A) \geq \underline\dim_B(A) \geq \dim_H(A)$, and unfortunately they are not stable under countable unions. However, in many cases as e.g.~for self-similar sets satisfying the open-set-condition (OSC) (see Section~\ref{self-similar-sets}), the box-counting dimension exists and is equal to the Hausdorff-dimension.

Note that there is also the following representation of the above limits:
\begin{lemma}
\label{lem-1}
  \begin{eqnarray}
    \underline\dim_B(F) &=& \sup\left\lbrace t\in\mathbb R: N_\delta(F)\delta^t \rightarrow \infty ~(\delta \rightarrow 0) \right\rbrace \\
    \overline\dim_B(F) &=& \inf\left\lbrace t\in\mathbb R: N_\delta(F)\delta^t \rightarrow 0 ~(\delta \rightarrow 0) \right\rbrace 
  \end{eqnarray}
\end{lemma}
In the limiting case $s=t$, $\lim_{\delta \downarrow 0} N_\delta(F) \delta^s$ can be either $0$, a positive finite value or $\infty$.

Later (see Section~\ref{sec:sausage}) we will use the following alternative formula for the box-counting dimension, which goes back to Minkowski (and thus explains why there are two different names for the same dimension):

\begin{theorem}
\label{th:minkowski-dim}
	Let $A\subset \Rd$ be a Borel-set, write $\mathcal L^d$ for the Lebesgue-measure on $\Rd$ and denote 
by\footnote{Here $d(x,A) := \inf(d(x,a): a\in A)$ is the distance of the point $x$ to the set $A$.} 
\[A_r:=\{x\in\Rd: d(x,A)\leq r\}\] the $r$-parallel set of $A$. Then 
\begin{eqnarray}	
	\overline\dim_B(A)&=&d-\liminf_{r\downarrow 0}{\dfrac{\log \mathcal L^d(A_r)}{\log r}} \\
	\underline\dim_B(A)&=&d-\limsup_{r\downarrow 0}{\dfrac{\log \mathcal L^d(A_r)}{\log r}}
\end{eqnarray}
\end{theorem}

\textit{Proof:} In the definition of the Box-counting dimension, $N_r(A)$ can be substituted by $P_r(A)$, the greatest number of disjoint $r$-balls with centers in $A$ (see e.g. \cite[Ch.5]{mattila95}). Then the Lebesgue-volume of the parallel sets can be sandwiched via 
\begin{equation}
P_r(A)\alpha(d)r^d \leq \mathcal L^d(A(r)) \leq N_r(A)\alpha(d)(2r)^d,  
\end{equation}
($\alpha(d)$ denotes the volume of the $d$-dimensional unit ball). After taking logarithms, divide by $\log r$ and consider the superior and inferior limits. $\square$

\begin{definition}
\label{def:minkowski-content}
  Let $A\subset \Rd$ be such that $s:=\dim_B(A)$ exists, and let the parallel sets $A_r$ be defined as above. Then the limit
  \[ \mathcal M(A) := \lim_{r \downarrow 0} r^{s-d}\mathcal L^d(A_r) \]
is called \textit{Minkowski-content} if it exists. 
\end{definition}
A result that will be used later is that for self-similar sets $F$ satisfying the Open-Set-Condition we have 
\[ 0 < \mathcal M(F) < \infty. \]
Dimitris Gatzouras proved this in \cite{gatzouras} (without the assumption of $F$ having polyconvex parallel sets).


\section{Dimensions of Measures}

\label{sec:dim-of-measures}

Borel-measures can be viewed as a generalization of Borel-sets, in the sense that a set can only describe the shape of an object, whereas a measure can describe its shape and its mass distribution. Accordingly, measures admit a richer ``fractal theory'' than sets. Here, however, only the framework for understanding the Local Dimension Method in section~\ref{sec:localdim} shall be laid.

Let $\mu$ be a locally finite Borel-measure on $\Rd$ and let $B(x,r)$ denote the ball in $\Rd$ around $x$ of radius $r$. 

\begin{definition}
	The limit
\[
	\overline{D}^s\mu(x) := \limsup_{r \rightarrow 0+}{\dfrac{\mu(B(x,r))}{(2r)^s}}
\]
is called \textit{upper s-density of $\mu$ at $x$}. Similarly, the \textit{lower s-density of $\mu$ at $x$} is defined as
\[
	\underline{D}^s\mu(x) := \liminf_{r \rightarrow 0+}{\dfrac{\mu(B(x,r))}{(2r)^s}}.
\]
\end{definition}

The dimensions of a measure are now defined as follows:

\begin{definition}
Let $\mu$ be a Borel-measure on $\mathbb R^d$
\begin{eqnarray*}
	\underline{\dim}\mu(x) &=& \liminf_{r\rightarrow 0+} \dfrac{\log\mu(B(x,r))}{\log r} \\
	\overline{\dim}\mu(x) &=& \limsup_{r\rightarrow 0+} \dfrac{\log\mu(B(x,r))}{\log r}
\end{eqnarray*}	
are called the \textit{lower} and \textit{upper local dimension} of $\mu$ at $x \in \mathbb R^d$.
\end{definition}

These limits express the power law behaviour of $\mu(B(x,r))$ as $r \rightarrow 0$. If both lower and upper local dimension agree at $x$ then $\mu$ is called \textit{simple} at $x$, and one writes $\dim\mu(x)$ for the common value and calls it the \textit{local dimension of $\mu$ at $x$}. Note that for $x\in \mathbb R^d$ outside the support of $\mu$ one has $\dim \mu (x) = \infty$, and that $\dim \mu (x) = 0$ if $x$ is an atom of $\mu$, i.e. if $\mu(\{x\}) > 0$.

There is the following connection between the densities and dimensions of Borel-measures:

\begin{lemma}
	Let $\mu$ be a Borel-measure.
	\begin{eqnarray*}
		\overline{\dim}\mu(x) &=& \sup\{t\geq 0: \underline D^t\mu(x) = 0\} = \inf\{t\geq 0: \underline D^t\mu(x) = \infty \} \\
		\underline{\dim}\mu(x) &=& \sup\{t\geq 0: \overline D^t\mu(x) = 0\} = \inf\{t\geq 0: \overline D^t\mu(x) = \infty \}
	\end{eqnarray*}
\end{lemma}

In other words, as a function of $t$, the lower density $\underline D^t \mu(x)$ jumps from $0$ to $\infty$ at the value $\overline{\dim}\mu(x)$, whereas the upper density performs this jump at the lower dimension.

\section{Self-similar sets}

\label{self-similar-sets}

Self-similar sets are the best-understood fractal sets and serve as the starting point for the development of the theory of fractal curvature (see Chapter~\ref{ch:fractal_curvature}).

A map $f:X\rightarrow X$ is called a \textit{contraction} if 
\[
\text{Lip}f:=\sup_{x\neq y}\frac{d(f(x),f(y))}{d(x,y)}<1 ,
\]
and it is a \textit{similarity} if for some constant $r>0$ we have $d(f(x),f(y))=r~d(x,y)$ for all $x,y \in X$. An \textit{Iterated Function System} (IFS) \label{IFS} is a system $\mathcal S=\{S_i\}_{i\in \{1,\ldots,N\}}$ of maps that are both contractions and similarities. 

Instead of $(X,d)$ we are now going to consider the metric space $(\mathbb K^d, d_H)$:
\begin{definition}
Let    \[\mathbb K^d := \{K \subset \Rd: K\textnormal{ is compact}\}. \] and let $K_\delta$ be the $\delta$-parallel set of $K$ from definition~\ref{def-1}. Then
  \begin{eqnarray*}
        d_H : (\mathbb K^d,\mathbb K^d) &\longrightarrow& [0,\infty),\\
	(K,L)&\longmapsto&\inf\{\delta>0:K_\delta\supset L \textnormal{ and } L_\delta\supset K\}.
  \end{eqnarray*}
is called the \textnormal{Hausdorff-metric} on $\mathbb K^d$.
\end{definition}

Given the IFS $\mathcal S$, consider the map 
\[S:\mathbb K\longrightarrow \mathbb K, ~~~ K\longmapsto \bigcup_{i=1}^N S_i(K).\]

Since the $S_i$ are Lipschitz-continuous, they map compact sets to compact sets, and since the union is finite, $S$ maps compact sets to compact sets, so $S$ is well-defined. 

\begin{lemma}
\label{lem:contraction}
  $S$ is a contraction on $(\mathbb K, d_H)$ with $\textnormal{Lip }S = \max_{1\leq i \leq N}\textnormal{Lip}S_i$.
\end{lemma}
\textit{Proof:} If $(S_i(A))_\delta$ contains $S_i(B)$ for all $i$, then $(\bigcup_{i=1}^NS_i(A))_\delta$ contains $\bigcup_{i=1}^NS_i(B)$, and the same is true if $A$ and $B$ are swapped. Thus
\[d_H\left(\bigcup_{i=1}^NS_i(A),\bigcup_{i=1}^NS_i(B)\right) \leq \max_{1\leq i \leq N}d_H\left(S_i(A),S_i(B)\right)
  \leq \max_{1\leq i \leq N}\textnormal{Lip}S_i ~ d_H(A,B).\]

\begin{theorem}
  \label{th-attractor}
  For the Iterated Function System $\mathcal S$ let $S:\mathbb K^d \rightarrow \mathbb K^d$ be defined as above. Then there is a unique non-empty compact set $F \in \mathbb K^d$ called the attractor of $\mathcal S$ that satisfies 
  \[ S(F) = F. \]
  Furthermore, if $E \in \mathbb K^d$ and $E \supset S(E)$, then 
  \[ F = \bigcap_{k=0}^\infty S^k(E). \]
\end{theorem}
\textit{Proof:} This is a simple consequence of lemma \ref{lem:contraction} and Banachs Fixed-Point theorem on $\mathbb K^d$ endowed with the Hausdorff-metric. 

For later use we shall distinguish two types of attractors / self-similar sets:
\begin{definition}
  The Iterated Function System $\mathcal S = \{S_1,\ldots,S_N\}$ with similarity ratios $r_i$ and the corresponding self-similar attractor $F$ are called \textit{h-arithmetic} if 
 $h\in \R $ is the greatest number such that $\log r_i \in h\mathbb Z$ $\forall i\in \{1,\ldots,N\}$. If no such number $h$ exists, $\mathcal S$ and $F$ are called \textit{non-arithmetic.}
\end{definition}

The following is a regularity property of iterated function systems:
\begin{definition}[Open Set Condition]
\label{def:OSC}
  An iterated function system $\mathcal S$ is said to satisfy the open set condition (OSC) if there exists an open set $O$ such that
\[ S_i(O)\subset O \forall i \text{ and } S_i(O) \cap S_j(O) = \emptyset \text{ for }i\neq j. \]
\end{definition}
All sets from section~\ref{sec:images} satisfy the OSC. Now we can state the main theorem on the dimension of self-similar sets (see e.g. \cite{falconer90}, theorem 9.3). 
\begin{theorem}
  \label{th-9.3}
  Let $\mathcal S = \{S_1,\ldots,S_N\}$ be an iterated function system with corresponding similarity ratios $r_i$, and assume that the OSC holds. Furthermore, let $s$ be the solution of $\sum_{i=1}^N r_i^s=1.$ Then the attractor $F$ satisfies 
\begin{enumerate}
  \item $0 < \mathcal H^s(F) < \infty$
  \item $\dim F = \underline\dim_B F= \overline\dim_B F = s$
\end{enumerate}
\end{theorem}
Note that the OSC is guarantees that $\Hs(F)$ is bounded away from $0$ in the above theorem. There is a corresponding result for self-similar measures, also due to Hutchinson (\cite{hutchinson81}):
\begin{theorem}
\label{th-self-similar-measure}
Assume the same conditions as in theorem~\ref{th-9.3}. Then there is a unique Borel-measure $\mu$ with $\mu(\Rd)=1$, called the $\mathcal S$-invariant measure, such that 	\[\mu(\cdot) = \sum_{i=1}^N r_i^s \mu (S_i^{-1}(\cdot)), \] and it satisfies \[ \mu(\cdot) = \dfrac{\mathcal H^s(F \cap \cdot)}{\mathcal H^s(F)}.\]
\end{theorem}

\begin{corollary}
\label{cor:support}
  The measure $\mu$ is spreading the unit mass uniformly over all of $F$.
\end{corollary}
(Here we use the word ``uniformly'' because $\mu$ is the restriction of the rotationally and translationally invariant $s$-dimensional Hausdorff-measure.)

\textit{Proof:} It suffices to show that supp~$\mu = F$. $F$ is compact, so closed. Let $x \in F \cap U$ for some open set $U$. Then there exists a finite sequence $i_1, i_2, \ldots, i_k$ in $\{1,\ldots,N\}$ such that 
$U_{i_1,\ldots i_k} := S_{i_k}^{-1} \circ \ldots \circ S_{i_1}^{-1} (U) \supset F$, so 
$\mu (U) \geq \prod_{j=1}^k r_{i_j}^s \mu (U_{i_1,\ldots i_k}) > 0$. $\square$

The specialty of this measure is that it spreads the unit mass uniformly on its support, which is a self-similar set of dimension $s$. By sampling points according to this probability distribution we can generate attractors of any iterated function system, and in fact this method has been used for the generation of the images from section \ref{sec:images} on a computer. The uniformity of the distribution is important, as sampling via a non-uniform distribution can require many more points until every pixel of the attractor is finally covered. 

Before ending this section, the following result shall be noted, which justifies the Local Dimension method of section~\ref{sec:localdim}:

\begin{proposition}
  \label{prop:obvious}
Let $F$ be the attractor of the IFS $\mathcal S$ which satisfies the OSC. Then the self-similar measure $\mu$ from theorem~\ref{th-self-similar-measure} satisfies
\[\dim \mu (x) = s\]
for all $x \in F$.
\end{proposition}


\section{Generating Attractors of IFSs}
\label{sec:IFS}

For later use (section \ref{sec:images}) we now describe how to generate a binary image of the attractor $F$ of an iterated function system. Chapter \ref{ch:discreteness} decribes the problem of discretization more thoroughly, but for now it is sufficient to know that we represent a binary image of an attractor $F \subset \mathbb R^2$ by all the pixels that are intersecting it (these pixels shall be called \textit{the pixels of $F$}). 

Fix $x_0\in\mathbb R^2$, and consider the sequence of points 
\begin{equation}
\label{eq:sequence}
x_{k} := S_{I_k}(x_{k-1}), ~~ k = 1,2,\ldots 
\end{equation}
where $I_k$ is an iid\footnote{independent and identically distributed} sequence of random variables with values in the set $\{1,\ldots,N\}$ having the probability distribution
\[P(I_k=i) = r_i^s.\]
(Recall that $N$ is the number of similarities in the IFS.) In a few steps we will see that as $k \rightarrow \infty$, $x_k$ runs through (almost) all the pixels of $F$.

\begin{lemma}
\label{lem:approx}
Let $c := \max\{r_i, i \in \{1,\ldots, N\}\}$ .
\begin{enumerate} 
  \item If $d(x_0,F) \leq N$ for some $N > 0$ , then $d(x_k,F) \leq c^k N.$
  \item If $x_0\in F$, then $x_k \in F$ for all $k \in \mathbb N$.
\end{enumerate}
\end{lemma}

\textit{Proof:} Let $E = \{x \in \mathbb R^2: d(x,F) \leq N\}$. Then $x_k \in S^k(E)$, so $d(x_k,F) \leq d_H(S^k(E),F)$ which is not greater than $c^k d_H(E,F) = c^k N$ by lemma \ref{lem:contraction}. This shows 1. $S(F) = F$ shows 2. $\square$

This means that if $x_k$ is not already in $F$, at least it is approaching $F$ exponentially quickly. If the distance is measured in pixelwidths, then practically this means that if the diameter $\vert F \vert$  of $F$ equals $3000$, if $x_0$ is not further away from $F$ than $30000$ and if the similarities all have ratio not higher than $0.8$, then $x_{100}$ will be closer to $F$ than $10^{-5}$ (pixelwidths). This should be enough to assert that almost no pixels outside $F$ are marked black. Also, due to the contractive nature of the algorithm, the iterations of $x_k$ are computationally stable, as calculation errors decay exponentially with $k$. 

Now consider the Borel-measure 
\[\nu(B) := \lim_{n\rightarrow\infty}\frac{1}{n} \# \{k:k\leq n, x_k \in B\}\]
which denotes the relative frequency of occurences of $x_k$ in the set $B$.

\begin{theorem}
\label{th:dynamical-system}
Let $x_k$ be as in equation (\ref{eq:sequence}) and assume $x_0 \in F$. Then
\[\nu(B) = \frac {\Hs(B\cap F)} {\Hs(F)} = \mu(B)\]
for any Borel-set $B \subset \Rd$.
\end{theorem}

\textit{Idea of proof:} Let $O$ be a feasible open set for the Open Set Condition. Then $\bar O \supset F$, $\nu(\bar O) = 1$, and for any finite sequence $i_1, \ldots, i_k$ in $\{1,\ldots,N\}$ the set 
$\bar O_{i_1, \ldots, i_k} := S_{i_k} \circ \ldots \circ S_{i_1} (\bar O)$ will satisfy 
$\nu(\bar O_{i_1, \ldots, i_k}) = \prod_{j=1}^k r_j^s$ which is also the relative frequency of the occurences of the string $i_1, \ldots, i_k$ in the infinite sequence generated by the $I_k$. So $\nu$ is invariant in the sense of \cite[p.16]{hutchinson81}, and thus of the same form as $\mu$ in theorem \ref{th-self-similar-measure}. $\square$

Now it is clear that the sequence $x_{101}, x_{102}, \ldots$ starting at 101 and proceeding to infinity will enter all pixels of $F$: Lemma~\ref{lem:approx} proves that (almost) the same pixels will be run through by $x_k$, regardless of whether or not $x_0 \in F$. Corollary~\ref{cor:support} and theorem~\ref{th:dynamical-system} say that $F$ equals the support of $\nu$, and since every pixel of $F$ contains at least one point of $F$ each one of these has a positive probability of being marked black while iterating $x_k$.

The internet provides a great variety of software which is using this principle to generate self-similar and also self-affine sets. For the generation of the pictures in section \ref{sec:images}, the program ``Fractal Explorer'' \cite{fractalexplorer} has been used.


\chapter{Fractal Curvature}
\label{ch:fractal_curvature}

Curvature Measures are well understood for sets such as differentiable manifolds or, more generally, for sets of positive reach \cite{federer1959cm}. The first subsection summarizes the main results for curvature measures in classical geometry, but in a slightly less general setting that also serves our purposes: We are only going to consider sets that lie in the convex ring $\mathcal R$ of finite unions of compact convex sets. Then we describe how Fractal Curvatures can be defined for sets $F$ whose $\eps$-parallel sets $F_\eps$ are in $\mathcal R$ for arbitrarily small $\eps$, and finally formulas for the fractal curvatures are given if these sets are self-similar.

\section{Curvature measures on the convex ring}
\label{sec:convex-ring}

\paragraph{Steiner-formula.}

Most of the statements from this subsection are proved in \cite{schneiderweil}. We write 
$\mathbb K$ for the system of all compact and convex subsets of $\Rd$.
First of all, recall the definition of the $\epsilon$-parallel set $F_\epsilon$ for subsets $F\subset \mathbb{R}^d$ with the euclidean metric $d(\cdot,\cdot)$ (theorem \ref{th:minkowski-dim}), and assume that $K \in \mathbb K$. Then the volume of the parallel sets of $K$ is given by the \textit{Steiner-formula}:
\begin{equation}
\label{eq-steiner}
 C_d(K_\epsilon)=\sum_{j=0}^d \epsilon^{d-j}\kappa_{d-j}C_j(K)
\end{equation}
Here, $C_d$ denotes the $d$-dimensional Lebesgue-measure, and $\kappa_j$ is the volume of the $j$-dimensional unit ball. The numbers $C_j(K)$ that are determining the coefficients of the polynomial in $\eps$ on the right are called \textit{intrinsic volumes}. In the 3-dimensional case, these have the following interpretations: $C_3(K)$ corresponds to the volume of $K$, $C_2(K)$ is the surface area of $K$, $C_1(K)$ equals (up to a multiplicative constant) the mean breadth of $K$, 
 and $C_0(K)$ equals 1 (being the Euler-characteristic of the non-empty convex set $K$). 


\paragraph{Local Steiner-formula.}

For $x\in \Rd$ let $p_K(x)$ be the point in $K$ that is nearest to $x$. There is only one such point since $K \in \mathbb K$. Also note that $p_K$ is (Lipschitz-) continuous. Write $\mathcal B$ for the system of Borel-subsets of $\Rd$. Now, for any $A\in\mathcal B$ define the \textit{local parallel set of $K$ at $A$} as\footnote{$U_\eps(K,A)$ can be visualized as $K \cap A$ together with all line segments of length $\eps$ emanating orthogonally outwards starting from all $y \in A \cap \partial K$. At cusps, the orthogonal directions are comprised in the normal cone.} 
\begin{equation*}
 U_\eps(K,A)=p_K^{-1}(A)\cap K_\eps.
\end{equation*}
For fixed $K$, 
\[\rho_\eps(K,\cdot):=C_d(U_\eps(K,\cdot))\] 
is a Borel-measure concentrated on $K_\eps$. If both $K \in \mathbb K$ and $A \in \mathcal B$ are fixed, it can be shown (see \cite{schneiderweil}, Satz~2.3.3) that $\rho_\eps(K,A)$ is a polynomial in $\eps$, 
\begin{equation}
 \rho_\eps(K,A)=\sum_{j=0}^d \epsilon^{d-j}\kappa_{d-j}C_j(K,A),
 \label{eq-localSteiner}
\end{equation}
and that the coefficients $C_j(K,A)$ are finite Borel-measures for each fixed $K \in \mathbb K$. 

Equation (\ref{eq-localSteiner}) is known as the \textit{local Steiner-formula}. The term ``local'' arises from the observation that the ordinary Steiner-formula (\ref{eq-steiner}) arises as a special case of the local Steiner-formula (\ref{eq-localSteiner}), namely for $A = \Rd$. 
This also shows $C_j(K,\Rd) = C_j(K)$, i.e. the intrinsic volumes are the total masses of the curvature measures $C_j(K,\cdot)$, and will be called $j$-th total curvatures later.

\paragraph{Curvature measures.}

For each $j \in \{0,\ldots,d\}$, $C_j(K,\cdot)$ is called the \textit{$j$-th curvature measure of $K$.} We note some important properties:
\begin{enumerate}
\label{curvature-properties}
 \item $C_j(K,\mathbb{R}^d)=C_j(K)$, i.e. the intrinsic volumes are equal to the total curvatures.
 \item $C_d(K,\cdot)$ is supported by $K$, and for $j=0,\ldots,d-1$ $C_j(K,\cdot)$ is supported by the boundary $\partial K$ of $K$.
 \item Curvature measures are motion-covariant, i.e. given a rotation or translation $g$, one has $C_j(gK,gA)=C_j(K,A)$.
 \item $C_k(\cdot,\cdot)$ is homogeneous of degree k: For $\lambda > 0$, $C_k(\lambda K, \lambda B) = \lambda^k C_k (K, B)$.
 \item Curvature measures are additive: $C_j(K,\emptyset) = 0$, and for $K,M \in \mathbb K$ such that $K \cup M \in \mathbb K$ one has
\begin{equation}
 C_j(K\cup M,\cdot)=C_j(K,\cdot)+C_j(M,\cdot)-C_j(K\cap M,\cdot).
 \label{eq-additive}
\end{equation}
\end{enumerate}
Note that in the last property, we were assuming $K\cup M \in \mathbb K$ because curvature measures are initially defined only for compact convex sets. However, via the \textit{inclusion-exclusion principle}
\begin{equation}
 C_j(\bigcup_{i=1}^m K_i,\cdot)
=\sum_{I\subset \{1,\ldots,m\}}{(-1)^{\#I-1}C_j\left(\bigcap_{i \in I} K_i, \cdot \right)}
\end{equation}
one can formally assign a map from $\mathcal B$ to $\mathbb R$ for each finite union $\bigcup_{i=1}^m K_i$ of sets $K_i \in \mathbb K$, as the intersections 
$\bigcap_{i \in I} K_i$ on the right-hand side are members of $\mathbb K$. It turns out that for different representations of the same union, $\bigcup_{i=1}^m K_i=\bigcup_{i=1}^{m'} K'_i~$, the maps $C_j(\bigcup_{i=1}^m K_i,\cdot) = C_j(\bigcup_{i=1}^{m'} K'_i,\cdot)$ coincide (\cite[Satz~2.4.2]{schneiderweil}). This means that $C_j(F,\cdot)$ is well-defined for any 
$F\in \mathcal R$, where $\mathcal R$ is the ring of \textit{polyconvex sets}, i.e. finite unions of sets in $\mathbb K$.

The maps $C_j(F,\cdot)$ turn out to satisfy properties $2,3,4$ and 5. However, they are non-negative Borel-measures only for $j=d$ and $j=d-1$; for $j\in\{0,\ldots,d-2\}$ they are \textit{signed} Borel-measures in general, i.e. they can have a non-zero negative variatonal part. For a definition see p.\pageref{def:signed-measure} of the appendix.

\section{Rescaled Curvature Measures}
\label{sec:concept}

Fractal sets are in general neither polyconvex nor of positive reach, so the classical definition of curvature does not apply. However, a great amount of geometric information about a fractal set $F$ can be encoded in the behaviour of the curvatures $C_j(F_\eps,\cdot)$ of its parallel sets as $\eps \downarrow 0$, provided of course that $C_j(F_\eps)$ is defined. In this section and the next, a summary of relevant results taken out of Steffen Winter's dissertation \cite{diss_winter} will be given.

To make sure that $C_j(F_\eps,\cdot)$ is defined, we assume that for every compact $F$ there exists an $\eps_0 > 0$ such that $F_\eps \in \mathcal R$ for all $\eps \in (0,\eps_0]$, and we say that $F$ has polyconvex parallel sets. More generality could be achieved by assuming that $F_\eps$ is a finite union of sets of positive reach. However, this weaker assumption has not led to any deeper insights yet, so we stick to the setting $F_\eps \in \mathcal R$ to simplify the notation, knowing that the results can be easily extended to the positive-reach-setting.

Consider now the expressions $\eps^t C_j(F_\eps,A)$ where $t\in\mathbb R$ and their behaviour as $\eps\downarrow~0$. If $t$ is chosen too big, they are likely to tend to 0, and for a too small $t$, they will behave eratically and probably be unbounded.\footnote{Note the analogy to $\delta N_\delta^s$ and $r^{s-d}\mathcal L^d(A_r)$ for the Box-counts in lemma~\ref{lem-1} and for the Minkowski-content in definition~\ref{def:minkowski-content}.} The interesting values of $t$ are the values at which $\eps^t C_j(F_\eps,A)$ switches from being unbounded to being bounded. It seems wise to avoid cases in which $\eps^t C_j(F_\eps,A)$ stays bounded but where the positive and negative parts $\eps^t C_j^+(F_\eps,A)$ and $\eps^t C_j^-(F_\eps,A)$ are actually unbounded; thus the focus will be on $\eps^t C_j^{var}(F_\eps,A)$. To further simplify things, only the global scaling constant will be studied, i.e. only the case $A = \Rd$. Recall that we write $C_j(F,\Rd)=C_j(F)$ for the $j$-th total curvature.

\begin{definition}[scaling exponents]
\label{def:scaling-exponents}
The \textit{(upper) k-th curvature scaling exponent} $s_k$ and the \textit{lower k-th curvature scaling exponent} $\underline{s_k}$ of a subset $F\subset\mathbb R^d$ are respectively defined as
\begin{eqnarray*}
 s_k(F) & := &\inf\left\{t\in\mathbb R:\lim_{\eps\downarrow 0}\eps^tC_k^{var}(F_\eps)=0\right\} \\
 \underline{s_k}(F) &:= & \sup\left\{t\in\mathbb R:\lim_{\eps\downarrow0}\eps^tC_k^{var}(F_\eps)=\infty\right\}.
\end{eqnarray*} 
\end{definition}

Let $\delta > 0$. With this definition, one can be sure that for every $A \in \mathcal B$ and every symbol $\bullet \in \{+,-,var\}$
\[\sup\left\{\eps^{s_k+\delta}C_k^\bullet(F_\eps,A): ~\eps \in (0,1]\right\}<\infty,\]
and that there exists an $A \in \mathcal B$ and a symbol $\bullet \in \{+,-\}$ such that 
\[\sup\left\{\eps^{\underline s_k-\delta}C_k^\bullet(F_\eps,A): ~\eps \in (0,1]\right\}=\infty.\]

Note that an alternative way to write down definition~\ref{def:scaling-exponents} (also compare definition~\ref{def-1} and lemma~\ref{lem-1}) is
\begin{eqnarray*}
  s_k(F) &=& \limsup_{\eps\downarrow 0}\dfrac{\log C_k^{var}(\Fe)}{-\log\eps},\\
  \underline{s_k}(F) &=& \liminf_{\eps\downarrow 0}\dfrac{\log C_k^{var}(\Fe)}{-\log\eps}.\\
\end{eqnarray*}
In general, $\underline{s_k}(F)\leq s_k(F) $, but for many e.g. self-similar sets equality prevails. 

We can now state the central definition:

\begin{definition}[fractal curvature]
  Let $F \subset \Rd$ be compact, and assume that there exists an $\eps_0 > 0$ such that $F_\eps\in\mathcal R$ for all $\eps \in (0,\eps_0]$. If the limits
\begin{eqnarray}
C_k^f(F)&:=&\lim_{\eps\downarrow 0}{\eps^{s_k}C_k(F_\eps)} \textnormal{ and}\\
C_k^{f,var}(F)&:=&\lim_{\eps\downarrow 0}{\eps^{s_k}C_k^{var}(F_\eps)}
\end{eqnarray}
 exist, then they are called the \textit{k-th fractal (total) curvature} of $F$ and the \textit{k-th fractal (total) variational curvature} of $F$. 
\end{definition}

Unfortunately, this limit does not exist in many cases. It often happens that the total curvatures exhibit the growth behaviour $C_k(F_\eps) \sim \eps^{s_k}, ~ (\eps \downarrow 0)$, but the expressions $\eps^{s_k}C_k(F_\eps)$ keep oscillating as $\eps$ approaches 0. In these cases, an \textit{average} limit often exists:

\begin{definition}[average fractal curvature]
  Let $F$ be as above. If the limits
\begin{eqnarray}
\overline{C}_k^f(F):=\lim_{\delta\downarrow0}\frac{1}{\log\eps_0 - \log\delta}
  \int_\delta^{\eps_0}\eps^{s_k}C_k(F_\eps)\frac{d\eps}{\eps}\\
\overline{C}_k^{f,var}(F):=\lim_{\delta\downarrow0}\frac{1}{\log\eps_0 - \log\delta}
  \int_\delta^{\eps_0}\eps^{s_k}C_k^{var}(F_\eps)\frac{d\eps}{\eps}  
\end{eqnarray}
 exist, then they are called the \textit{k-th average fractal (total) curvature} of $F$ and $k$-th average fractal (total) variational curvature of $F$.
\end{definition}

$\overline{C}_k^f(F)$ can be considered as a generalization of $C_k^f(F)$, because if $C_k^f(F)$ exists then so does $\overline{C}_k^f(F)$, and in that case both their values coincide. The same statement applies to $\overline{C}_k^{f,var}(F)$ and $C_k^{f,var}(F)$. Note that for $k\in\{d-1,d\}$ the measures $C_k(F_\eps)$ and $C_k^{var}(F_\eps)$ coincide and so do their averages, so that in part~II where $d=2$ we only need to worry about $C_0^{var}(\cdot)$.

\paragraph{Consistency with classical curvature.}

Finally, it should be noted that the definition of fractal curvature $C_k^f(F)$ is an extension of the definition of the total curvature $C_j(F)$ from section~\ref{sec:convex-ring}: If $F\in\mathcal R$ and $C_k(F) \neq 0$, the scaling exponents $s_k$ all equal $0$, and in that case 
\[C_k^f(F)=\lim_{\eps\rightarrow 0}{\eps^{s_k}C_k(F_\eps)} = C_k(F)\] 
which follows from the continuity of the measure $C_k(\cdot)$ and the convergence $F_\eps \rightarrow F$ in the Haussdorff-metric $d_H(\cdot,\cdot)$. Motion invariance and homogeneity also carry over to the fractal versions:

\begin{proposition}
\label{prop:rescaling}
Let $F$ be a Borel-set such that the limit $\overline{C}_k^f(F)$ exists. Let $g$ be a motion in $\Rd$ and let $\lambda > 0$. Then the limits $\overline{C}_k^f(gF)$ and $\overline{C}_k^f(\lambda F)$ also exist, and
\[\overline{C}_k^f(gF)=\overline{C}_k^f(F) \textnormal{\hspace{1cm} and \hspace{1cm}} 
	\overline{C}_k^f(\lambda F)=\lambda^{s_k+k}\overline{C}_k^f(F).\]
\end{proposition}

\paragraph{Fractal curvature measures.} 

As seen above, the rescaled expressions $\eps^{s_k}C_k(F_\eps)$ lead to a suitable definition for the \textit{total} fractal curvature of a compact set $F$. Similarly, one can consider the expressions $\eps^{t}C_k(F_\eps,A)$ with suitable $t$ to arrive at a \textit{local} fractal curvature of $F$. The value $t$ is best chosen independently of $A$, as the limits 
$\lim_{\eps\downarrow0}\eps^{t}C_k(F_\eps,A)$ will then more likely be Borel-measures in $A$ if they exist. 

\begin{definition}[weak convergence of measures]
Let $\mu$ and $\mu_\eps$ be measures on the measurable space $(X,\Sigma)$ for every $\eps \in (0,\eps_0]$. The $\mu_\eps$ are said to converge weakly against $\mu$ as $\eps \downarrow 0$, if
\[\int f d\mu_\eps \rightarrow \int f d\mu ~~ (\eps \downarrow 0)\]
for every measurable function $f:X\rightarrow \R$. In this case we write
\[\lim_{\eps \downarrow 0} \mu_\eps = \mu.\]
\end{definition}

\begin{definition}[fractal curvature measure]
\label{def:fractal-curvature-measure}
If the limit
\[C^f_k(F,\cdot):=\lim_{\eps\downarrow 0}\eps^{s_k}C_k(F_\eps,\cdot)\]
exists, it is called $k$-th fractal curvature measure of $F$. Similarly, the limit
\[\overline C^f_k(F,\cdot):=\lim_{\delta\rightarrow0}\frac{1}{\log\eps_0 - \log\delta}
  \int_\delta^{\eps_0}\eps^{s_k}C_k(F_\eps,\cdot)\frac{d\eps}{\eps}\]
is called the average fractal curvature measure of $F$ if it exists.
\end{definition}

\section{Fractal curvature of self-similar sets}
\label{sec:self-similar-sets}

If fractal curvatures exist, and what values they have if they exist, turns out to be a difficult question. A natural starting point for these questions are (deterministically) self-similar sets, since they are probably best-understood among all fractal sets. In his dissertation \cite{diss_winter}, Steffen Winter showed the existence of all fractal curvature measures for a subclass of the self-similar sets and gave formulas for their calculation.

\paragraph{The class $\Sd$.}

Recall that the definition of fractal curvature requires that $F_\eps \in \mathcal R^d$ for all $\eps$ in some interval $(0,\eps_0]$. Now self-similar sets have the following fortunate property:

\begin{theorem}[\cite{llorentewinter07}, Proposition 4.6]
\label{th-par-sets}
Let $F$ be a self-similar set. Then 
\[\exists \eps > 0 : F_\eps \in \mathcal R^d \Longleftrightarrow \forall \eps > 0 : F_\eps \in \mathcal R^d. \]
\end{theorem}

In other words: Either all parallel sets of a self-similar set are polyconvex or none. This means that in order to find out whether or not there is an interval $(0,\eps_0]$ on which $F_\eps \in \mathcal R^d$ is satisfied, we just need to check $F_\eps \in \mathcal R^d$ for \textit{any one} $\eps> 0$. 

In order to simplify the notation we introduce the following

\begin{definition}[$\Sd$]
Denote the class of all compact self-similar subsets of $\Rd$ which satisfy the open set condition and have polyconvex parallel sets by $\Sd$. 
\end{definition}

\paragraph{Scaling exponents.}

For $\Sd$-sets, the scaling exponents always satisfy $s_k \leq s-k$. Even more can be said:

\begin{theorem}[\cite{diss_winter}, theorem 1.3.2]
\label{th:s_k-upper-bound}
  Let $F\in\Sd$. Then for all $k \in \{0,\ldots,d\}$ the expression $\eps^{s-k}C_k^{var}(F_\eps)$ is bounded.
\end{theorem}

In all known examples of $\Sd$-sets with non-integer dimension, in fact $s_k = s-k$, but an according theorem has not been proved yet. It seems as if for $s_k < s-k$ to happens it takes a ``degenerate'' fractal having non-empty interior.\footnote{An example for such a set is the $d$-dimensional unit cube, which can be viewed as the union of $2^d$ just-touching copies of size $\frac{1}{2}$.} 

Winter gave no theorem guaranteeing a lower bound for $s_k$ as universally as theorem~\ref{th:s_k-upper-bound}, but at least the following criterion:

\begin{theorem}[\cite{diss_winter}, Theorem 1.3.8]
\label{th:s_k-lower-bound}
  Let $F \in \Sd$ and $k\in\{0,\ldots,d\}$. Let $r_{min}$ be the smallest similarity ratio of the IFS of $F$. Suppose there exist constants $\eps_0,\beta>0$ and a Borel set $B \subseteq ((O^c)_{\eps_0})^c =: O_{-\eps_0}$ such that 
  \[C_k^{var}(F_\eps,B) \geq \beta\]
  for each $\eps \in (r_{min}\eps_0,\eps_0]$.
  Then for all $\eps<\eps_0$
  \[\eps^{s-k}C_k^{var}(F_\eps) \geq \beta\eps_0^{s-k}r_{min}^s > 0.\]
\end{theorem}

We apply this theorem to all sample sets $F$ from section \ref{sec:images} except the Cantor dust and the Koch curve (these are not members of $\Sd$), in order to determine their scaling exponents: Let ${\rm conv}(F)$ be the convex hull of $F$. Then for each $F$ the interior $O = {\rm conv}(F)^{\circ}$ of the convex hull is a feasible open set. Choose e.g. $\eps_0 = 0.01d$ where $d$ is the perimeter of $F$ and $B = O_{-\eps_0}$ and see that for $k\in\{0,1,2\}$ and $\eps \in (r_{min}\eps_0,\eps_0]$, $C_k^{var}(F_\eps,B) \geq \vert C_k(F_\eps,B) \vert \geq c$ for some value $c > 0$. Thus we note the following

\begin{corollary}
\label{cor:s_k}
For all $k\in\{0,1,2\}$ and all sets $F$ from section~\ref{sec:images} except the Cantor dust and the Koch curve,
  \[\underline{s_k}(F) = s_k(F) = s-k. \]
\end{corollary}

\paragraph{Formulas for (total) fractal curvatures.}

There are explicit formulas for the fractal curvatures of $\Sd$-sets. These involve the \textit{curvature scaling functions}:

\begin{definition}
  For a self-similar set $F$ associated with the IFS $\mathcal S = \{S_i,\ldots,S_N\}$, the \textit{$k$-th curvature scaling function} is defined by
\[R_k(\eps) := C_k(F_\eps) - \sum_{i=1}^N \textbf{\textup{1}}_{(0,r_i]}(\eps)C_k((S_iF)_\eps)\]
\end{definition}

These, together with the renewal theorem (see e.g. \cite[Theorem~7.2]{tecfracgeom}), are the tools which are making it possible to derive the following formulas: 

\begin{theorem}[\cite{diss_winter}, theorem 1.3.6]
\label{th:curvature-formula}
Let $F\in\Sd$ be the self-similar set corresponding to the IFS 
$\mathcal S=(f_1,\ldots,f_N)$ with similarity factors $(r_1,\ldots,r_N)$.
Then for $k=0,\ldots,d$ the following holds:
\begin{enumerate}
 \item The $k$-th average fractal total curvature $\overline{C_k^f}(F)$ 
	exists and is equal to 
	\begin{equation}
 	X_k=\frac{1}{\eta}\int_0^1\eps^{s-k-1}R_k(\eps)d\eps,
	\end{equation}
	where $\eta = -\sum_{i=1}^Nr_i^s\log r_i$ and $s=\dim_H(F)$.
  \item If $F$ is non-arithmetic, the $k$-th fractal total curvature 
	$C_k^f(F)$ exists and equals $X_k$.
\end{enumerate}
\end{theorem}

Note that for $k=d$, the above formula holds true even if $F$ does \textit{not} have polyconvex parallel sets. This has been shown by Gatzouras in \cite{gatzouras}. Recall that $C_2^f(F)$ is actually nothing else than $\mathcal M(F)$, the Minkowski-content of $F$ (see definition~\ref{def:minkowski-content}), and with this in mind we are going to call $\overline C_2^f(F)$ the \textit{average Minkowski-content} of $F$. Gatzouras also showed that for self-similar sets satisfying the open set condition always $X_d > 0$, a result which we are going to refer to by saying that ``$F$ has positive average Minkowski-content.''

\paragraph{Curvature measures.} Fractal curvature measures, which are the limits described in definition~\ref{def:fractal-curvature-measure}, have the following simple representation for $\Sd$-sets:

\begin{theorem}[\cite{diss_winter}, Theorem 1.5.1]
  Let $F\in \Sd$, and assume $s_k = s-k$. Then the average fractal curvature measures from definition~\ref{def:fractal-curvature-measure} exist and equal
\[\overline C^f_k(F,\cdot) = \overline C^f_k(F) \frac{\Hs(F\cap\cdot)}{\Hs(F)}.\]
If $F$ is non-arithmetic, the fractal curvature measures also exist and equal
\[C^f_k(F,\cdot) = C^f_k(F) \frac{\Hs(F\cap\cdot)}{\Hs(F)}.\]
\end{theorem}

Thus fractal curvature measures for $\Sd$-sets are just multiples of the self-similar measures we have already encountered in theorem~\ref{th-self-similar-measure}. The fact that they are all constant multiples of each other is the basis of the geometric invariant characteristic defined in section~\ref{sec:beyond}.

\part{Image Analysis}

\chapter{Limitations imposed by Discreteness}
\label{ch:discreteness}

The main reason for starting a new part here is the change of setting: Space is no longer $d$-dimensional and continuous, but 2-dimensional and discrete. The necessity to deal with discrete space, however, does not only stem from the fact that we are analysing images: If one thinks of real and concrete counterparts of the abstract fractals from part one, these have to be composed of finitely many atoms\footnote{This means that in nature, a 3D-Sierpi\'{n}ski gasket of diameter 1 meter can at most be iterated around 50 times before its smaller copies reach the size of atoms. Fractals of infinite recursion depth only exist in the theory of part one.}, just as images are composed of finitely many pixels. In this context, if nature is to be modeled by fractal geometry, the discretization problem arises sooner or later.

\paragraph{The Discretization Process.}

We assume a rectangular grid for the pixels, which we model by the set
\[\mathbb Z^2_\delta := \left\lbrace 
	(k\delta,l\delta) \in \mathbb R^2: k, l \in \mathbb Z
	\right\rbrace, \] 
where $\frac{1}{\delta}$ is called resolution.
Later it will turn out useful to endow $\mathbb Z^2_\delta$ with the euclidean metric $d(\cdot,\cdot)$ which is rescaled in such a way that the distance between horizontally and vertically neighbouring pixels is $1$. 

A binary image $I$ shall be represented by 
\[ I : \mathbb Z^2_\delta \rightarrow \{0,1\} \]
where the value $1$ corresponds to a black pixel and $0$ to a white pixel. We always assume that $I^{-1}(\{1\})$ is finite, i.e. that there are only finitely many black pixels and thus images are bounded. 

Finally, we represent a bounded subset $A\subset\mathbb R^2$ by the binary image $I_A$ which satisfies 
\[I_A\left((k\delta,l\delta)\right) = 1 \Longleftrightarrow 
	\left[k\delta,(k+1)\delta\right) \times \left[l\delta,(l+1)\delta\right) 
	\cap A \neq \emptyset,\]
i.e. exactly those pixels are black whose corresponding areas have non-empty intersection with $A$.

\paragraph{``Lacunarity''.}

Unsurprisingly, the fractal dimension (Hausdorff-dimension) of a fractal set is far from determining its structure completely. Figure~\ref{fig:stack} (taken from \cite{mandelbrot:mfl}) shows the first construction stages of one-dimensional Cantor-like dusts, with iterated function systems as follows:
\begin{figure}[h]
  \centering
  \includegraphics[width=10cm]{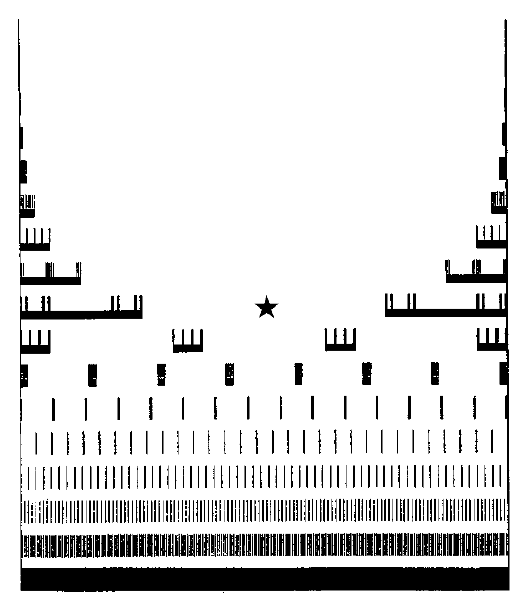}
  \caption{A collection of Cantor dusts of dimension $\frac{1}{2}$.}
  \label{fig:stack}
\end{figure}
In each line which is $k$ steps away from the starred line there are $2^{k+1}$ similarities, all of which have ratio $r = \frac{1}{4}^{k+1}$. In the upper half, the similarities place copies of the unit interval towards the ends of the interval, whereas in the lower half the similarities are spread uniformly. In every case, the similarity dimension is the value $s$ for which $2^{k+1} \times (\frac{1}{4}^{k+1})^s = 1$, i.e. $s =$ $1\over 2$. The size of gaps, determining the \textit{lacunarity} of a fractal, is high on the top and low at the bottom. 

Now note how the results of the discretization process differ from the topmost to the bottommost fractal: Going upwards from $\star$, at some point only the first and last pixel of the unit interval will be marked black, whereas going downwards far enough will result in \textit{every} pixel of the unit interval being marked black. Image analysers will not be able to differenciate between the Cantor-dust and 2 singletons resp. the whole unit interval. The obvious effect is that high lacunarity will result in an underestimation and low lacunarity in an overestimation of the fractal dimension. 

\paragraph{``Suitable'' binary images.}

The degree of fractality a binary image can convey is strongly limited by the length of chains of still visible iterates: If one can only see one smaller copy of the image in itself, almost all of the fractality has been lost in discretization, whereas not so much is lost if one can see e.g. {\tiny a smaller copy of} {\scriptsize a smaller copy of} {\footnotesize a smaller copy of} {\small a smaller copy}. Especially estimators based on the regression on the logarithmic scale\footnote{These are e.g. the box-counting algorithm, the sausage method and the method proposed in chapter~\ref{ch:main}.} depend on this sort of information. If the smaller copies are only supposed to be statistically similar to the bigger copy, additionally to a high recursion depth a large number of images is advantageous.

In this sense and with the above notion of lacunarity in mind, it seems that fractals are more suitable for image analysis the greater their similarity ratios are. It also seems that by excluding fractals of very low similarity ratios around $\frac{1}{10}$ one also avoids cases of extreme lacunarity as on the far ends of $\star$ in figure~\ref{fig:stack}. The sample images from section~\ref{sec:images} were chosen in this spirit.

\paragraph{Various definitions of lacunarity.}

The notion of lacunarity arose from the need to characterize fractals beyond their dimension. As this concept is far too complicated to be expressed by a mere number, there is still no uniformly accepted definition of lacunarity. 

In 1994 (\cite{mandelbrot:mfl}), B. Mandelbrot speaks of ``shell-lacunarity'' which is effectively the same as (the reciprocal value of) Minkowski-content. D. Gatzouras took this as a motivation to prove in 1999 (\cite{gatzouras}) that self-similar sets satisfying the OSC always have a well-defined average Minkowski-content (regardless of the parallel sets being polyconvex or not). With formulas from S. Winter's dissertation, the author calculated the average Minkowski-contents $C_1^f(F)$ for the sets $F$ from figure~\ref{fig:stack}. The somewhat surprising result is that at $\star$ $C_1^f(F)$ is minimal; proceding downwards $C_1^f(F)$ tends to $+\infty$; proceding upwards it increases first and from $k=3$ on it decreases and approaches the value $8.65...$. Thus the inverse average Minkowski-content does not describe the lacunarity effect from figure~\ref{fig:stack} ideally, as one would expect the starred dust to have lower lacunarity than the dust above it. However, it is consistent in that it is correctly monotonuous when there is no transition between two models (as in the step from $\star$ to the line above).

In the same work, B. Mandelbrot mentions ``gap-lacunarity'', which is (the reciprocal value of) $C_0^f(F)$, and he has already mentioned it to be a constant multiple of $C_1^f(F)$, without even the assumption of self-similarity of the dust on the real line. (The author does not know if this idea has been followed anywhere in the literature, in spaces of any dimension.)

In physics, lacunarity is preferably determined via the ``gliding box algorithm'' (see e.g. \cite{plotnick1996lag}): Given an object $M$ and a box $B(r)$ of size $r$, the mass $s$ of $M\cap B(r)$ is measured as $B$ is gliding through space (according to a never specified but seemingly uniform distribution). $Q(s,r)$ is the then obtained probability distribution of the masses $r$ according to the box size $s$. The first two moments are denoted by $Z_1(r)$ and $Z_2(r)$, respectively, and  
\[\Lambda(r) := \frac{Z_2(r)}{(Z_1(r))^2}\] is defined as the \textit{lacunarity} for the box-size $r$. Plots of $\Lambda(r)$ against $r$ on a logarithmic scale can be characteristic for certain random fractals. 

For sets as e.g. percolation clusters, the correlation of mass in different directions depicted by angular sectors has been examined by Mandelbrot et al, and it has been found that ``antipodal correlation'' can serve as a measure of lacunarity, too.

\chapter{A Review of Dimension Estimation Methods}
\label{ch:review}

In order to have a reference frame for the estimators $\hat s$ and $\Gamma_k$ of the method proposed in chapter~\ref{ch:main}, the author searched for methods of the estimation of dimension and of the (average) Minkowski-content in the literature. A large amount of literature has been found on dimension estimates, but nothing on the estimate of Minkowski-content. The reason for this probably is that, in natural sciences, the latter value is devoid of any meaning. However, in section~\ref{sec:beyond}, measurements of the (average) Minkowski-content will be used as a normalization factor for the 0th and 1st fractal curvatures, thus providing a new set of computable geometric invariants, at least for $\Sd$-sets.

Of the vast literature on dimension estimates the author ruled out the part which deals with dimension estimates of graphs of continuous functions, or more generally with estimates of the fractal index $\alpha$ of a stochastic process. It seems like drawing a binary image of a graph given a time series of data does not seem apropriate for estimates of the fractal index, and in recent literature estimators for $\alpha$ are all based on the data themselves and not on a binary image (see e.g.~\cite{sei}). 
Thus in this chapter only fractal dimension estimation methods will be reviewed. 

\section{The Box-Counting Method}

Recall definition~\ref{def-1} from chapter~\ref{ch:fractal-geometry}. If the box-counting dimension $\dim_B(F)=:s$ of a fractal $F$ exists, $N_\delta$ has the growth behaviour 
\begin{equation} 
\label{eq-box-growth}
N_\delta \sim C \delta^{-s} 
\end{equation}
where $C$ is some positive constant. Taking logarithms,
\begin{equation}
  \label{eq-box-loggrowth}
  \log N_\delta \sim \log C + s\log(1/\delta),
\end{equation}
so on logarithmic paper, the Box-counting-dimension appears as the slope of an asymptote. All box-counting algorithms exploit this behaviour of $N_\delta$, and a typical algorithm runs as follows: 
\begin{enumerate}
  \item Set the \textit{box-size in pixels} $\delta \in \mathbb N$ to a maximum value $\delta_{max}$.
  \item \label{step}Count the number $N_\delta$ of disjoint boxes of size $\delta$ that intersect the set $F$ in question.
  \item Record $N_\delta$ and decrease $\delta$. If $\delta \geq \delta_{min}$, go back to step 2.
  \item Finally, plot $\log N_\delta$ against $-\log\delta$, and fit a line to the data using the least squares method. Use the slope of the line as an estimator $\hat s^{(m)}$ for $s=\dim_B(F)$.
\end{enumerate}

Typically, $\delta_{max}$ is somewhere around one fourth to one third of the image diameter, whereas $\delta_{min}$ ranges from one to four pixels. In step 3, $\delta$ is decreased in such a way that the points in the logarithmic plot are being close to equidistant, i.e. 
$\delta_{new} = \lfloor \delta_{old}/a\rfloor$ for some constant $a > 1$, e.g. $a = 1.4$. 

Note that in step~\ref{step} the union of all disjoint boxes that intersect $F$ will certainly cover $F$. However, this union will not necessarily be made up of the minimal number of boxes possible, as would be required by the definition of $N_\delta$. Thus algorithms achieve more accurate (and usually higher) estimates if, each time step~\ref{step} is run through, the box-grid is shifted several times until the minimum covering number is found. Finding this minimal cover is incorporated in the box-counting algorithm of ``FracLac'' (\cite{fraclac}), a plug-in of the open-source image analysis software ``ImageJ'' \cite{imagej} which has been used for the box-counting estimates from chapter~\ref{ch:results}.

\section{Estimating Local Dimension}
\label{sec:localdim}

As pointed out in section \ref{sec:dim-of-measures}, Borel-sets can be generalized by measures which are supported by these sets, and in this setting the notion of local dimension arises quite naturally. In dynamical systems, an attractor can be assigned the measure of the relative frequency
\[\mu(B) := \lim_{n\rightarrow\infty}\frac{1}{n} \# \{k:k\leq n, x_k \in B\}\]
of a particle with coordinates $x$ being inside a Borel-set $B$ at discrete-time measurements $k\in\mathbb N$. If the attractor is ``strange'', typically the local dimension of this measure will vary throughout the set and will not be an integer.

In \cite{cutler1989eds} Cutler and Dawson describe statistical methods of estimating the local dimension of a measure on the basis of points $x_k$ sampled according to this measure. An important application of their method that they had in mind probably was to describe the geometrical properties of an attractor.

\paragraph{The estimator.} Cutler and Dawson consider the nearest neighbour statistics
\[\rho_n(x) := \min_{1 \leq i \leq n}{\Vert X_i - x \Vert}\]
where $x \in \supp \mu$ is a point of the support of the measure $\mu$ and the $X_i$ are independent random variables sampled according to the law $\mu$. They showed that if the local dimension is ``simple at $x$'', i.e. if 
\[\overline \dim \mu (x) = \underline \dim \mu (x) = \alpha(x),\]
then\footnote{With the conventions $\log 0 = -\infty$, $\frac{1}{0}=+\infty$, $\frac{1}{\infty}=0$ and admitting $\infty$ as a valid limit, this statement is also valid for $x$ outside of the support of $\mu$. In that case, $\alpha(x)=+\infty$; if $x$ is an atom of $\mu$, then $\alpha(x)=0$.} 
\[l_n(x) := \dfrac{\log (b\rho_n(x))}{-a \log n} 
\stackrel{a.s.}{\longrightarrow} \dfrac{1}{\alpha (x)}.\]
The positive parameters $a$ and $b$ are not affecting the limit of $l_n(x)$, however they have a biasing influence on $l_n(x)$ as an estimator for $\alpha(x)$. Cutler and Dawson suggest the values $a = 1$ and $b= 2$ as for a uniform distribution on an interval this choice is neutralizing the bias term with the lowest order in $n$. 

\paragraph{Binary images.}

If the binary image is representing an $s$-set $F$ (i.e. $0 < \mathcal H^s (F) < \infty$), then the simplest $s$-dimensional measure that can be assigned to $F$ is maybe 
\[\mu(\cdot) := \dfrac{\mathcal H^s (\cdot \cap F)}{\mathcal H^s (F)},\]
which spreads the unit mass uniformly over $F$. 
In order to be able to utilize the above method, we have to find a discrete analogon $\hat \mu$ of $\mu$ for the binary image; to this end, we spread the unit mass uniformly over the black pixels. Then we recall from section~\ref{sec:IFS} that we assume that the black pixels are exactly those which contain at least one point of the support of $\mu$. But since this point could be located anywhere in the pixel, we spread the mass of each black pixel evenly across its square. The author believes this to be the best way of discretizing the measure $\mu$ from above.

\paragraph{The algorithm.} 

Part of the work on this thesis was the creation of the java-class \texttt{LocalDimension} in GeoStoch (\cite{geostoch}) that is implementing the above method for binary images. It calculates a histogram of estimated local dimensions of the measure $\hat \mu$ at randomly picked points. 

The algorithm implemented in the above mentioned java-class runs as follows:
\begin{enumerate}
 \item Sample $m$ points $y_1,..., y_m$ (the ``test-points'') according to $\hat \mu$ in the following way: Record the coordinates of the black pixels in one big pixel-array \texttt B. Then pick one pixel according to a uniform distribution on (1,length(\texttt B)) and add a uniformly distributed number between 0 and 1 to both coordinates. 
 \item Create an array \texttt L of length $m$ containing the nearest-neighbour statistics $l(y_i)$ for each $i = 1,\ldots,m$, and initialize all entries \texttt{L[i]} to have the value $l_0(y_i) = d$ which is the diameter of the image.
 \item As $j$ ranges from $1$ to $n$, sample the point $x_j$ according to $\hat \mu$ as in step 1 and update the array \texttt L via \texttt{L[i] = }$\min\{\texttt {L[i]}, d(x_j,y_i)\}$ for each $1\leq i \leq m$.
 \item Draw a histogram for the values $\hat \alpha(y_i) := l_n(y_i)^{-1}$, $1 \leq i \leq m$, where 
$l_n(y_i) = \texttt{L[i]}$.
 \item Choose either the modal value of the histogram or the arithmetic mean of the 
$\hat \alpha(y_i)$ to be an estimate for the global Hausdorff-dimension of $F$.
\end{enumerate}
\begin{figure}[!ht]
	\centering
\includegraphics[width=\textwidth]{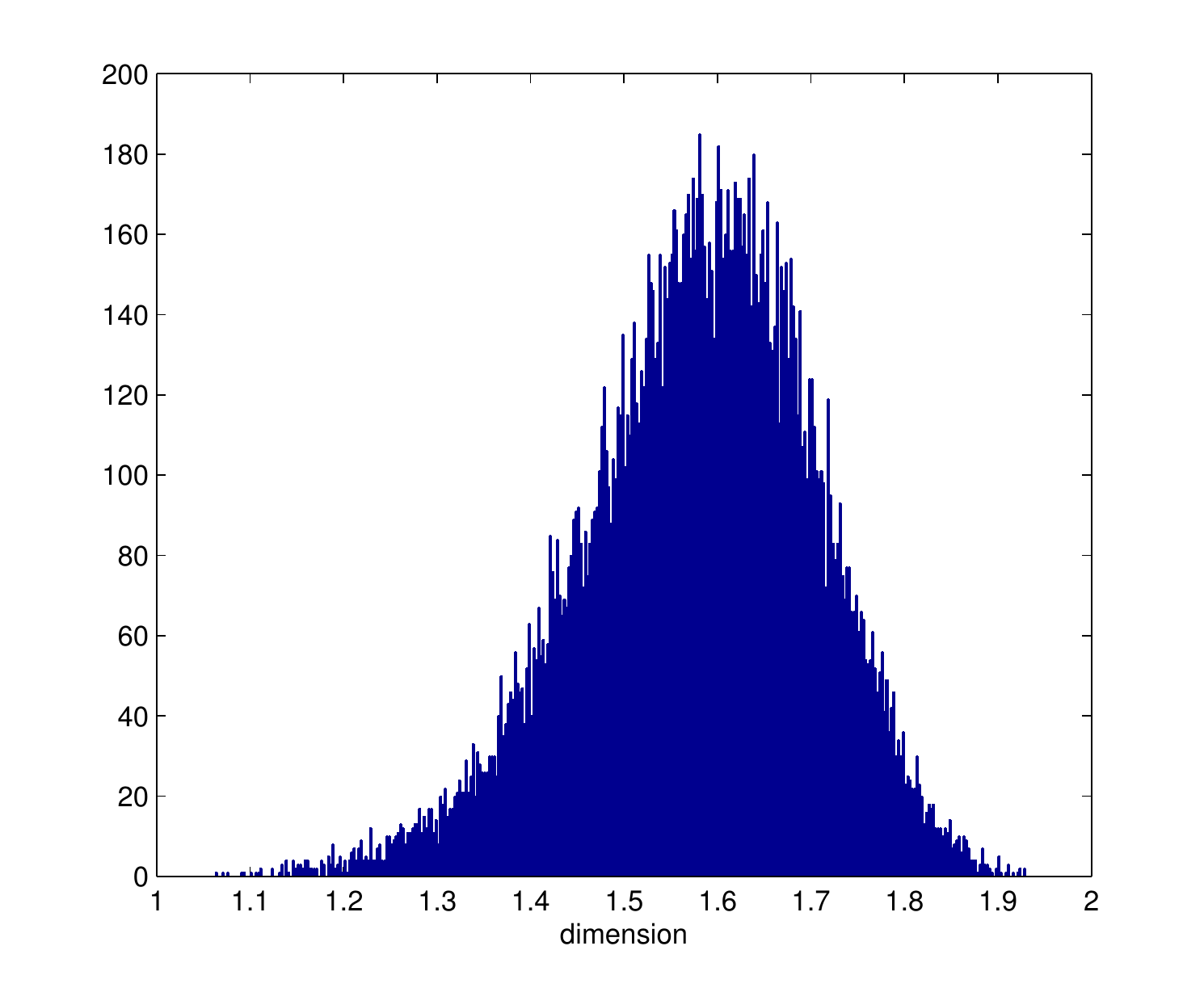}
\caption[A local dimension histogram]{A typical histogram of local dimensions of the self-similar triangle $\triangle$. The measured local dimensions at 20000 test points has been represented by columns of width 0.0025. The image has 383313 black pixels, and 300000 sample points were used. The average of the dimension estimates lies at 1.578, close to the real value 1.588...}
        \label{fig:hist}
\end{figure}
Figure~\ref{fig:hist} shows a typical histogram returned by the above algorithm for the self-similar triangle $\triangle$ from section \ref{sec:images}. 

\paragraph{Limitations.} Note that the histogram is far from its desired form, which would be one single column at the fractal dimension $s\approx 1.588$. However, since the local dimension of the measure $\hat \mu$ actually equals 2 everywhere on its support, we are glad to achieve this result at all. The number $m$ of test points can be increased almost arbitrarily, meaning that the histogram can actually approximate the distribution of the $l_n(y)^{-1}$ as $y$ varies over the support of $\hat \mu$. However, the number $n$ of sample points cannot be increased arbitrarily. If it increases beyond the number of black pixels, the pixel-scale local dimension character of $\hat \mu$ begins to shine through, and dimension estimates are approaching 2. The best value for $n$ seems to be somewhere around 80\% of the number of black pixels. 

For self-similar pictures which have a constant local dimension, the arithmetic mean of the values $\hat \alpha (y_i)$ yielded good estimates for the fractal dimension, except for the Koch-curve whose dimension was overestimated by 0.09, and the Tripet, where the value was 0.10 too low.

Finally, note that measures, as opposed to sets, can be represented more accurately if grey-values are admitted for the pixels, meaning that a pixel can be set darker the more mass the measure has at this pixel, whereas a set either intersects a pixel or not. If in the construction procedure of the self-similar images which is using the self-similar measure from section~\ref{sec:IFS} the relative frequency of entries of points in a pixel is measured and the grey value is calculated accordingly, it might be that a \textit{refined} discretized measure $\hat \mu$ will yield better dimension estimates. 

\section{The Sausage Method}
\label{sec:sausage}

Recall the Minkowski-representation of the Box-counting-dimension, theorem~\ref{th:minkowski-dim}: If $\dim_B(F)$ exists, the shrinking rate of the Lebesgue-volume of the parallel sets $\mathcal L^d(F_r)$ for $r \rightarrow 0$ is directly connected to the Box-counting-dimension $\dim_B(F)$ via
\[ \mathcal L^d(F_r) \sim \mathcal M(F) r^{d-\dim_B(F)} \]
where $\mathcal M(F)$ is the Minkowski-content of $F$. As before (see equations \ref{eq-box-growth} and \ref{eq-box-loggrowth}), on a logarithmic scale the growth behaviour is linear,
\[ \log\mathcal L^d(F_r) \sim \log\mathcal M(F) + (d-\dim_B(F)) \log r,\]
which is the basis of the following 

\paragraph{Algorithm.} 

\begin{figure}[h]
  \centering
  \includegraphics[width=\textwidth]{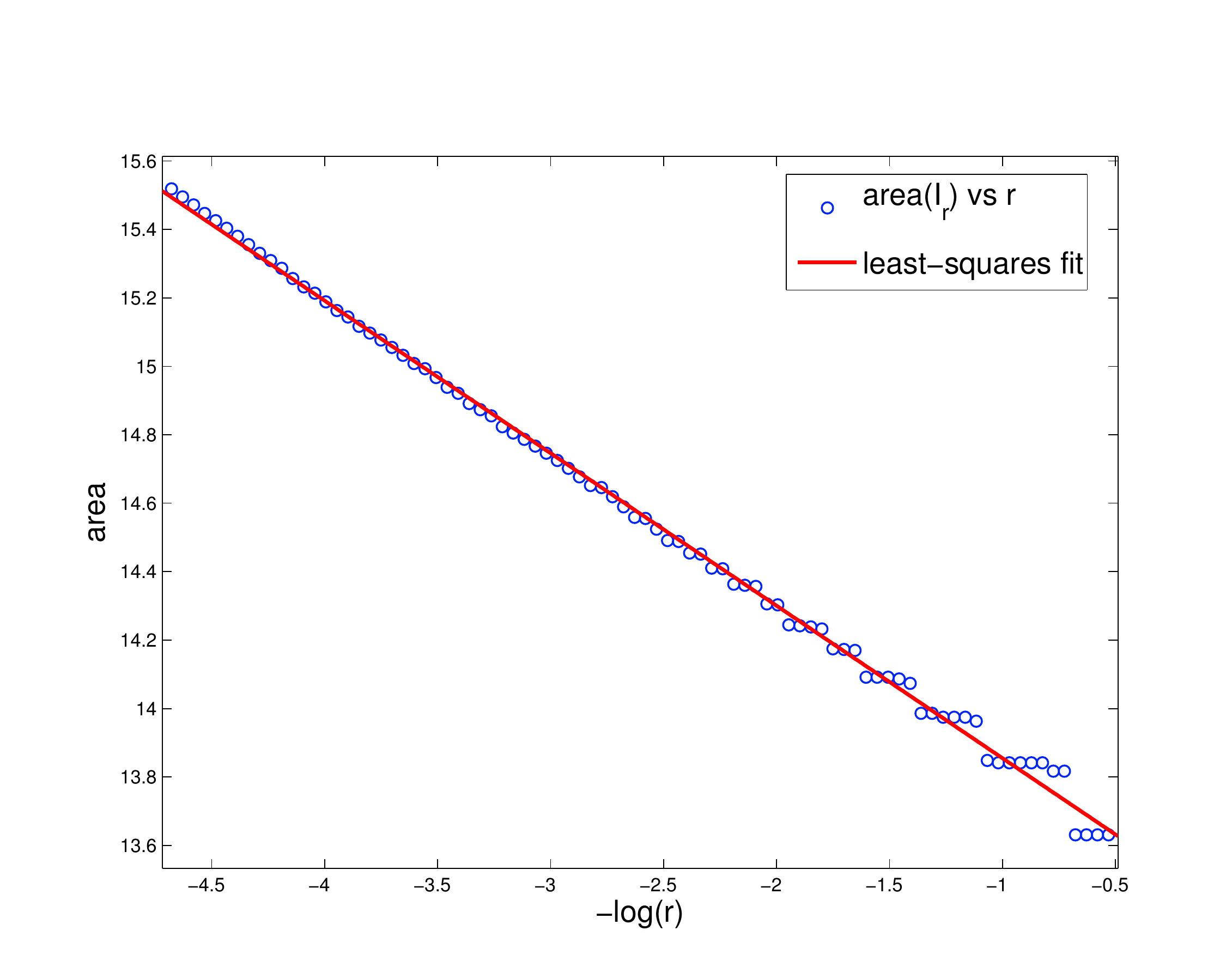}
\caption[A plot for the sausage method]{A plot of the areas of the dilated binary image against the dilation radii. The underlying fractal set is the Sierpi\'{n}ski Tree from section~\ref{sec:canonical}.}
  \label{fig:sausage-plot}
\end{figure}

Now proceed according to the following steps:

\begin{enumerate}
  \item Calculate the distance-transform $D_I$ of the binary image $I$. This is an image that records the distances to the nearest black pixels:
  \begin{eqnarray*}
  D_I : \mathbf P &\rightarrow	& \R \\
                p &\mapsto	& d\left( p,I^{-1}\left( \{1\}\right) \right)
  \end{eqnarray*}
\item For radii $r$ that are equidistantly distributed on the logarithmic scale (the maximal and minimal radius can be chosen the same way as $\delta_{max}$ and $\delta_{min}$ in the Box-counting method) calculate the areas of the $r$-dilated binary images 
\[I_r := D_I^{-1}\left( [0,r] \right).\]
\item Plot these areas against $r$ on logarithmic paper and fit a least squares line to the data points. The slope serves as an estimate for $\dim_B(F) - 2$ (see figure~\ref{fig:sausage-plot}). 
\end{enumerate}

This procedure has been implemented by the author in the class \texttt{Curvature2D} of the GeoStoch-library as the special case 
\[(\texttt{useEuler, useBdlength, useArea}) = (\texttt{false, false, true}),\] 
for a description see section~\ref{sec:algorithms}.
\begin{figure}[!h]
  \centering
  \includegraphics[width=\textwidth]{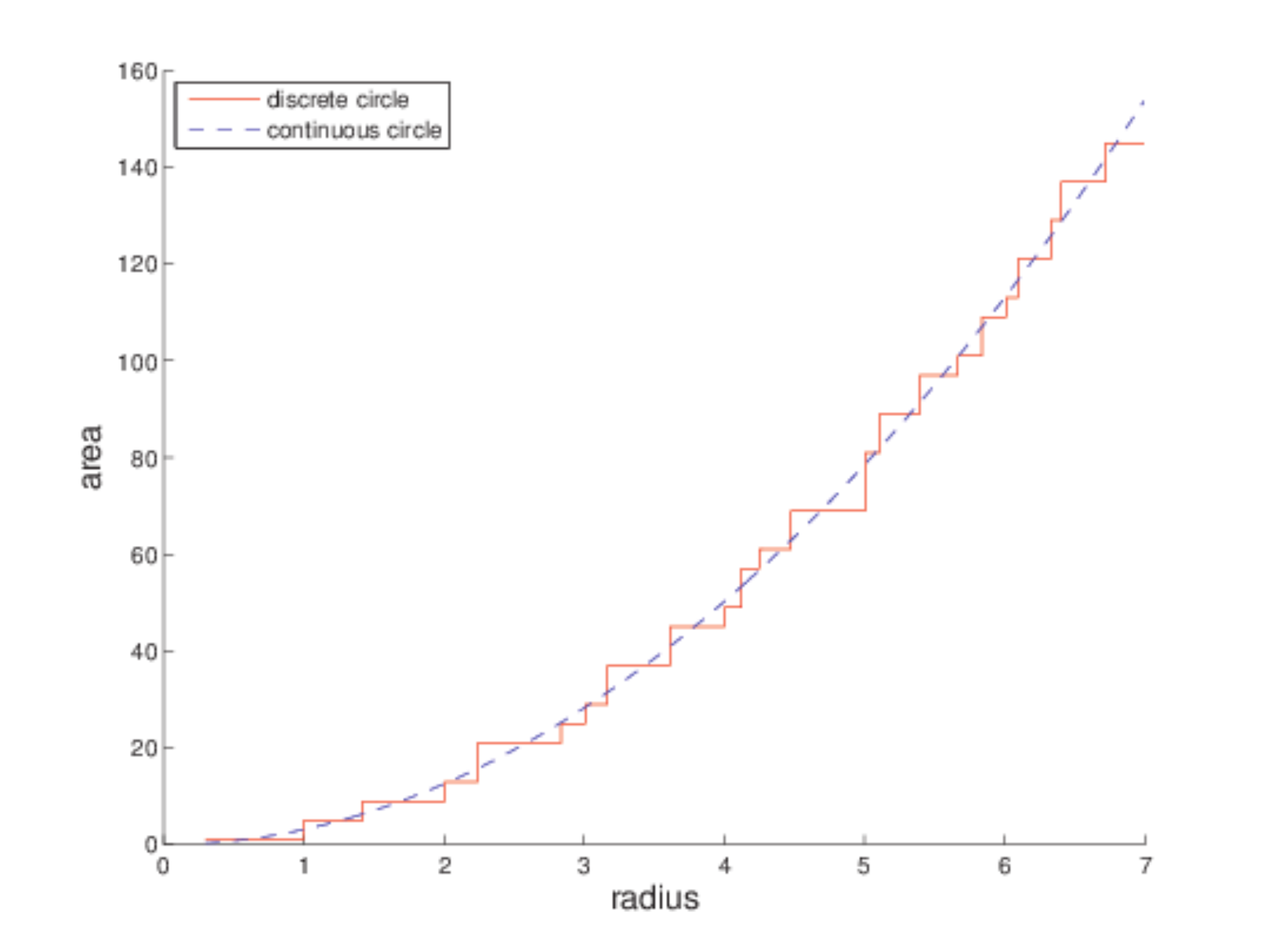}
  \includegraphics[width=\textwidth]{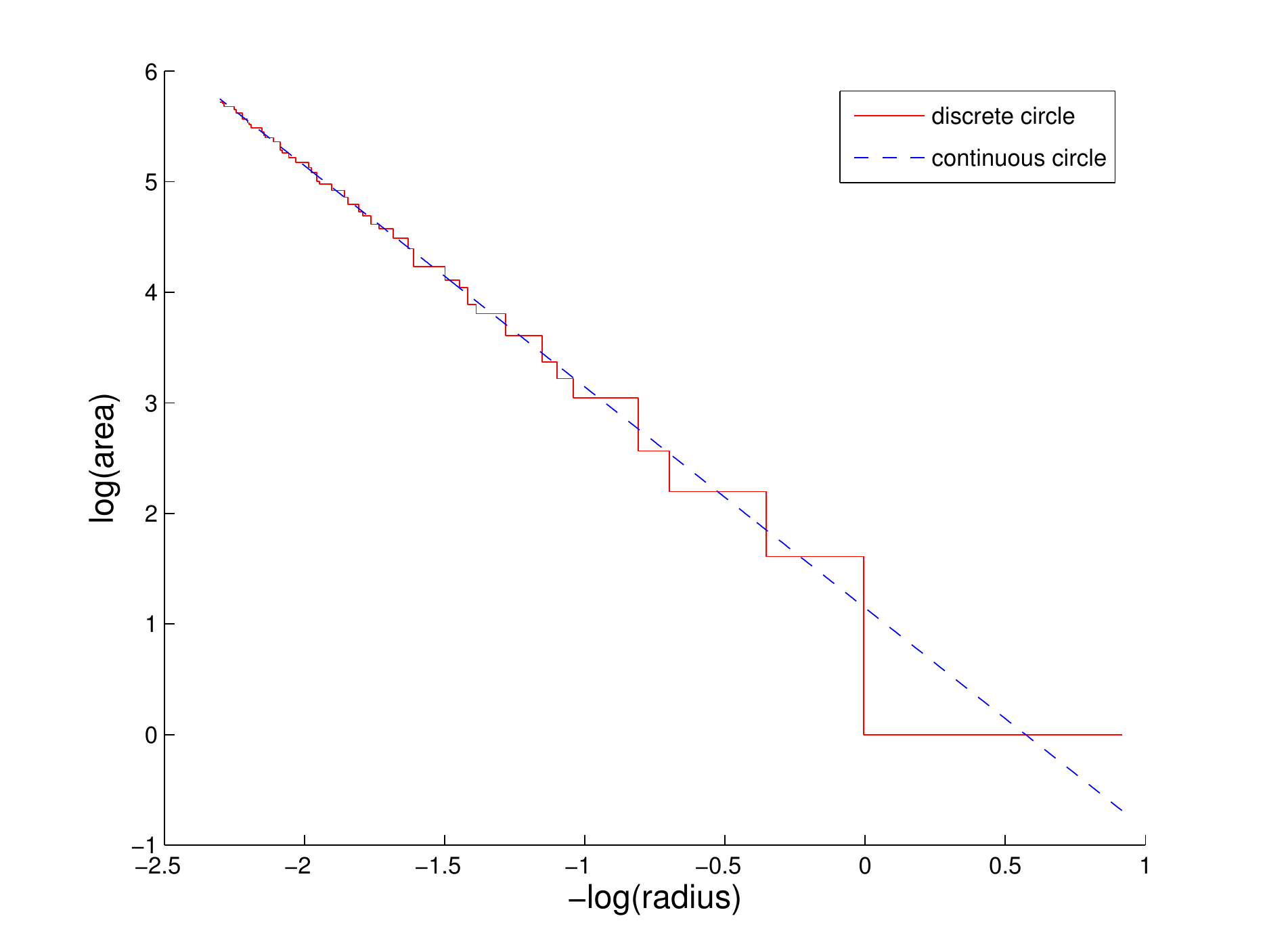}
  \caption[Area of discretized circle]{The areas of the discrete and continuous disk, as a function of the radius.}
  \label{fig:discrete-area}
\end{figure}

It should be noted how in figure~\ref{fig:sausage-plot} stairs appear for small radii. In fact, since the dilated images $I_r$ are all discrete, their area is piecewise constant as $r$ decreases. On a logarithmic scale this effect is more prominent for small radii than for big radii, as figure~\ref{fig:discrete-area} illustrates. 
This means that, whenever sparsely distributed dilation radii are used (e.g. to reduce the run-time of the algorithm), care should be taken as not to use those radii which are close to a jump of the area. A good strategy is to use the \textit{optimal-area radii}, i.e. those radii for which the discretized disks $B_{discrete}(o,r)$ centered at the pixel in the origin $o$ have an area which is exactly equal to $\pi r^2$. With the discretization mentioned above, i.e. 
\[ B_{discrete}(o,r) := \left\lbrace p \in \mathbf P: \vert p - o \vert \leq r \right\rbrace, \]
where $p$ and $o$ are the coordinates of the center of two pixels, these radii can be found as the horizontal intersections of the two graphs in figure~\ref{fig:discrete-area}, whereas the jumps occur at the vertical intersections. 

The first optimal-area radii (rounded to 4 digits) are
 \[0.5642,~~ 1.262,~~ 1.696,~~ 2.585,~~ 3.432,~~ 3.785,~~ 4.406,~~ 4.687,~~ 5.322, \ldots\]
the exact values of which are 
\[ \sqrt{\frac{1}{\pi}}, \sqrt{\frac{5}{\pi}}, \sqrt{\frac{9}{\pi}}, \sqrt{\frac{37}{\pi}}, \sqrt{\frac{45}{\pi}}, \sqrt{\frac{61}{\pi}}, \sqrt{\frac{69}{\pi}}, \sqrt{\frac{89}{\pi}},\ldots\]
corresponding to the areas 
\[1,5,9,21,37,45,61,69,89,\ldots\]
More optimal area radii of size up to 400 are listed in the class \texttt{OptimalRadii} of the GeoStoch library (\cite{geostoch}).
These numbers have been found numerically with Matlab as the zeros of the function 
\[ \mathcal D(r) =  \mathcal L^d (B_{discrete}(o,r)) - \mathcal L^d ( B(o,r)). \]
Also see the paragaph on the sausage method by Stoyan and Stoyan \cite{stoyan1994frs}, who suggest admitting corner points of pixels as possible centers of discrete circles.

\section{Further Methods}

\paragraph{Variations of Box-counting.}
For images of deterministic fractals, variations of the box-counting method have been proposed: Sandau and Kurz suggest the "Extended Counting Method" (\cite{sandau1997mfd}), which is depicting the highest dimension measurable by box-counting if the position of the coarsest grid is varied. Mart\'inez-L\'opez, Cabrerizo-V\'ilchez and Hidalgo-\'Alvarez (\cite{martinezlopez2001ime}) suggest replacing the number $N_\delta(F)$ of boxes of size $\delta$  by the number \[N'_{\delta,s}(F) := \sum_{j=1}^\infty  \left\vert \frac{A_j^{(\delta)}\cap F}{\delta} \right\vert^s,\] 
where $A_j^{(\delta)}$ are the (just touching) boxes of size $\delta$ used for covering $F$, $\vert \cdot \vert$ is the diameter of a set and $s$ is the (unknown) Hausdorff-dimension. An initial value $s_0$ is chosen, and iteratively $s_{i+1}$ is computed as the least squares fit of a line to the plot of $N'_{\delta,s_{i}}(F)$ versus $-\log\delta$ on a logarithmic scale. It was shown empirically that $s_i$ converges very quickly. 

\paragraph{Stochastic processes of fractal index.}

A large class of random fractals for which dimension estimates are important is generated by stochastic processes. As an example, rough surfaces are often modeled by graphs of stationary gaussian processes: These are families of real-valued random variables satisfying
\[X_t\sim \mathcal N(0,\sigma^2)~ \forall t\in \Rd. \]
Adler (\cite{adler1981grf}) showed that if the covariance function
\[\gamma(t) := cov(X_0,X_t) \]
satisfies
\[\gamma(t) = \gamma(0) - c\vert t \vert ^\alpha + o(\vert t \vert ^\alpha)
  ~~~(t \rightarrow 0)\]
for some $0 < \alpha < 2$ and $c>0$, then the graph 
\[G_X := \left\lbrace (t,X_t) \in \mathbb R^{d+1}: t \in \Rd \right\rbrace \]
of $X_t$ has Hausdorff-dimension $d - \frac{\alpha}{2}$ almost surely. Hall and Roy (\cite{hall1994rbf}) extended this result to processes of the form $g(X_t)$ for certain smooth $g$ and $X_t$ is as above.

Other important examples are the paths of fractional Brownian motion (fBm) with Hurst-parameter $H \in (0,1)$: This is a (non-stationary) centered gaussian process $B_t$, i.e. 
\[B_t \sim \mathcal N(0,\sigma_t^2) \]
with covariance function
\[\gamma(s,t) = cov(B_t,B_s) = E(B_tB_s) = \frac{1}{2} \left(s^{2H} + t^{2H} - \vert t-s \vert ^{2H}  \right). \] 
It is not hard to see that the two processes $\{a^{-H}B_{at},t\geq0\}$ and $\{B_t, t\geq 0\}$ have the same distribution, and thus $B_t$ is said to be statistically self-similar. The paths of $B_t$ are fractal sets of Hausdorff-dimension $1 + H$. 

Thus the knowledge of the dimension of these random fractals is equivalent to the "fractal index" $\alpha$ resp. $H$ of the according processes, and in the last 20 years a big effort has been put into constructing estimators of the fractal index. A good summary can be found in the recent PhD-thesis of Sei Tomonari (\cite{sei}).

\chapter{Image Analysis with Fractal Curvature}
\label{ch:main}

In the following we are going to propose a dimension estimation method that exploits the connection between the scaling exponents $s_k$ from definition~\ref{def:scaling-exponents} and the dimension $s:=\dim(F)$ of a fractal set $F$. Similarly to the sausage method, we are going to consider the parallel sets $F_\eps$, but additionally to the area $C_2(\Fe)$ we are going to measure the other (total) curvatures $C_1(\Fe)$ and $C_0(\Fe)$ as well (i.e. half the boundary length and the euler number), and perform a linear regression simultaneously on all three data sets. It appears that the quality of the measurements of the $C_k(\Fe)$ is decreasing with $k$, but in several cases, especially for non-arithmetic $F$, the additional data increase the quality and robustness of the dimension estimate. We assume 
$F \in \Sd$ throughout this chapter, and that the scaling exponents satisfy $s_k = s-k$. As noted in section~\ref{sec:self-similar-sets}, this seems to be the case whenever $F \in \Sd$ and $\dim F \notin \mathbb Z$, and it has been verified for the sample sets of section \ref{sec:images} in corollary~\ref{cor:s_k}.

\section{The Regression Model}


\paragraph{The logarithmic scale.}

Since iterated function systems shrink a set $E$ with the ratio $r_i$ every time a similarity $S_i$ is applied, the size of the sets $S_k\circ\ldots\circ S_1(E)$ decays exponentially with $k$. Therefore logarithmic rescaling is useful to describe growth phenomena for these fractals. Now let us express the growth behaviour of the curvatures of the parallel sets on a logarithmic scale. 
First, recall that two functions $f$ and $g:D\rightarrow \R$ are said to be asymptotic to each other for $x\rightarrow a \in \bar D$ if\footnote{By $o(1)$ we just mean a sequence converging to 0.} 
\[\frac{f(x)}{g(x)} = 1 + o(1), ~ x\rightarrow a,\]
provided of course $g(x)\neq0$ in some open set containing $a$.

Now the renewal theorem tells us that  
\begin{equation}
\label{eq:parallel-sets-growth}
\eps^{s-k}C_k^{var}(\Fe) \sim p_k(\eps) ~~~(\eps \downarrow 0),
\end{equation}
where $p_k: (0,1] \rightarrow \R$ is either constant (if $F$ is non-arithmetic) or periodic of multiplicative period $e^{-h}$ (if $F$ is $h$-arithmetic), i.e. $p_k(e^{-h}\eps) = p_k(\eps)$. Note that the expression on the left of (\ref{eq:parallel-sets-growth}) is positive for all $\eps > 0$ and all $k$, so ``$\sim$'' is well-defined.

For brevity we introduce the notation 
\begin{equation}
\label{eq:average}
\bar f := \lim_{\delta\downarrow 0} \frac{1}{-\log\delta}\int_\delta^1f(\eps)\frac{d\eps}{\eps}
=\lim_{T\rightarrow\infty}\frac{1}{T}\int_0^Tf(e^{-t})dt  
\end{equation} 
if the limits exist for a function $f:(0,1]\rightarrow \R$, and call it the \textit{average of $f$}. Note that functions that are asymptotic to each other have the same averages, so with $g_k(\eps):=\eps^{s-k}C_k^{var}(\Fe)$ we have
\[\bar p_k = \bar g_k = \overline C_k^{f,var}(F).\]
For a concise notation we decompose $p_k$ into two parts:
\begin{equation}
\label{eq:normalized-periodic}
p_k(\eps) = \overline C_k^{f,var}(F) p_k^0(\eps),  
\end{equation}

so that $\bar p_k^0 = 1$. Now (\ref{eq:parallel-sets-growth}) can be rewritten as 
\begin{equation}
  \dfrac{C_k^{var}(\Fe)}{\eps^k} \sim \overline C_k^{f,var}(F) p_k^0(\eps) \eps^{-s},
  ~~~ (\eps \downarrow 0), 
\end{equation}
and after taking logarithms and setting $\eps := e^{-x}$ we finally arrive at 
\begin{equation} 
\label{eq:log-master}
  \log \left( \dfrac{C_k^{var}(F_{e^{-x}})}{e^{-xk}} \right) 
		\sim
 \log \overline C_k^{f,var}(F) + \log p_k^0(e^{-x}) +sx, ~~~ (x \rightarrow \infty).
\end{equation}
Recall that $C_d^{var}(F_\eps)$ is just the area of $F_\eps$, and $C_{d-1}^{var}(F_\eps)$ is just half the length of the boundary of $F_\eps$. In general we have no way of measuring $C_0^{var}(\cdot)$ with an image analyser, because it is not additive on $\mathbb K$ and thus the inclusion-exclusion principle does not work. Sometimes, however, $C_0^{var}(F_\eps)$ can be precisely calculated as it is just equal to 
$2N(F_\eps)-C_0(F_\eps)$, where $N$ is the number of connected components (see p.\pageref{app:C0var} of the appendix). In case $C_0^{var}(F_\eps)$ cannot be calculated, it will be excluded in the following regression.

\paragraph{Simultaneous regression.}

Let $\mathcal X = \left\lbrace x_1,\ldots, x_m \right\rbrace $ be a set of $x$-values. Having equation (\ref{eq:log-master}) in mind, introduce the following variables:
\begin{eqnarray}
\label{eq:def-variables}
  y_{kj}  &:=& \log \left( \dfrac{C_k^{var}(F_{e^{-x_j}})}{e^{-x_jk}} \right) \\
  D_k &:=& \log \overline C_k^{f,var}(F) \\
  q_{kj} &:=& \log p_k^0(x_j).
\end{eqnarray} 

\begin{figure}[!h]
  \centering
  \includegraphics[width=\textwidth]{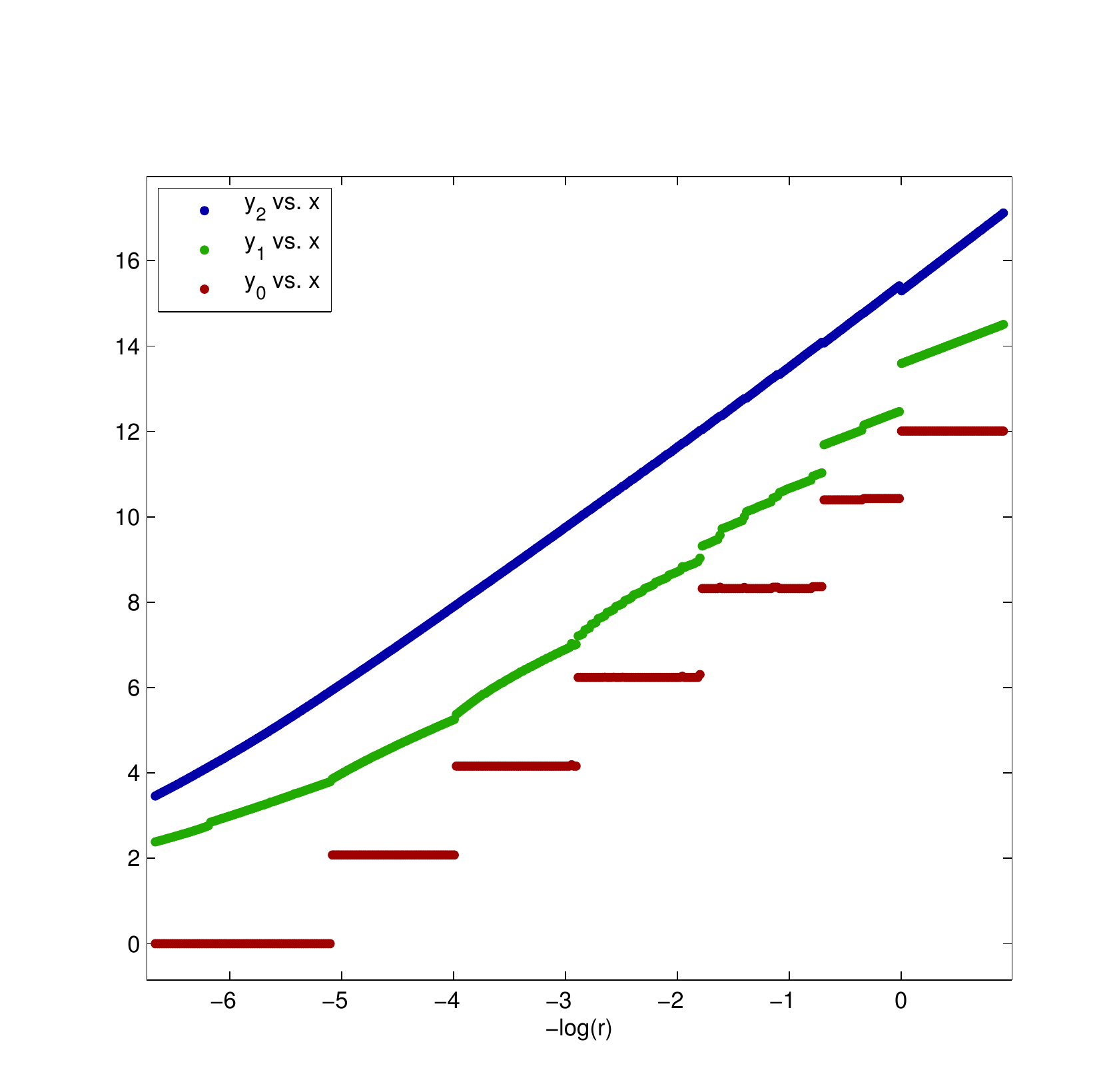}
  \caption[Plot of $y_k$ against $x$]{A plot of $y_k$ against $x$ for the 2D Sierpi\'{n}ski Carpet. Notice that for too big dilation radii (i.e. $r>\exp(5)\Leftrightarrow x<-5$) $C_0^{var}(F_r)$ turns to $1$, and the logarithms become $0$. For too small dilation radii (i.e. $r<1 \Leftrightarrow x>0$) the measured curvature $C_k(F_r)$ is piecewise constant as the maximum pixel-resolution is reached, yielding the slopes 0,1,2 for $y_0, y_1, y_2$ respectively.}
  \label{fig:linplot}
\end{figure}

Equation (\ref{eq:log-master}) suggests that, given that the approximation $"\sim"$ is good enough, for fixed $k$ the $y_{kj}$ all lie on the graph of a function which is the sum of an affine function $D_k + sx$ and a periodic function $q_k(x)$ (of additive period $h$ if $F$ is $h$-arithmetic). Figure~\ref{fig:linplot} shows a plot for the Sierpi\'{n}ski Carpet of size $3000\times3000$ pixels, where the corresponding $y_k$-values have been calculated with an image analyser.
This suggests the following regression model:
\begin{equation}
\label{eq:regression-model}
    \begin{pmatrix}
  Y_{0j} \\ Y_{1j} \\ \vdots \\ Y_{dj}
  \end{pmatrix}
 =
 \begin{pmatrix}
   D_0 \\ D_1 \\ \vdots \\ D_d
 \end{pmatrix}
 +
 \begin{pmatrix}
   q_{0j} \\ q_{1j} \\ \vdots \\ q_{dj}
 \end{pmatrix}
 + 
 \begin{pmatrix}
   s \\ s \\ \vdots \\ s
 \end{pmatrix}
  x_j
 +
 \begin{pmatrix}
   \delta_{0j} \\ \delta_{1j} \\ \vdots \\ \delta_{dj} 
  \end{pmatrix}
\end{equation}
where the $\delta_{kj}$ are random variables which are modelling discretization and measurement errors. Since we do not know any better, we assume (rather optimistically) that for each $k$ and $j$ they are normally distributed 
$\delta_{kj}\sim\mathcal N(0,\sigma^2)$ for some $\sigma^2$, and that they are mutually independent. 

The problem is now to estimate the unknown values $s$ and\footnote{The variables $D_k$ turn out to be averages of the logarithms of the $k$-th total variational curvatures of the parallel sets. 
By putting $\widehat{C_k^f(F)}:=\exp(D_k)$ one could define another average of the fractal curvatures; the author did not follow this idea as he had no desire to reformulate all the results of Winter's dissertation.} 
$D_0, \ldots, D_d$, given some $x$-values $x_j$, $j\in \{1,\ldots, m\}$ and the corresponding $d+1$-dimensional point cloud 
\[\mathcal Y := \left\lbrace
	\left( 
	Y_{0j},\ldots,Y_{dj}\right)^T: Y_{kj} = Y_k(x_j),~ j=1,...,m \right\rbrace. \]
\begin{figure}[!h]
  \centering
  \includegraphics[width=\textwidth]{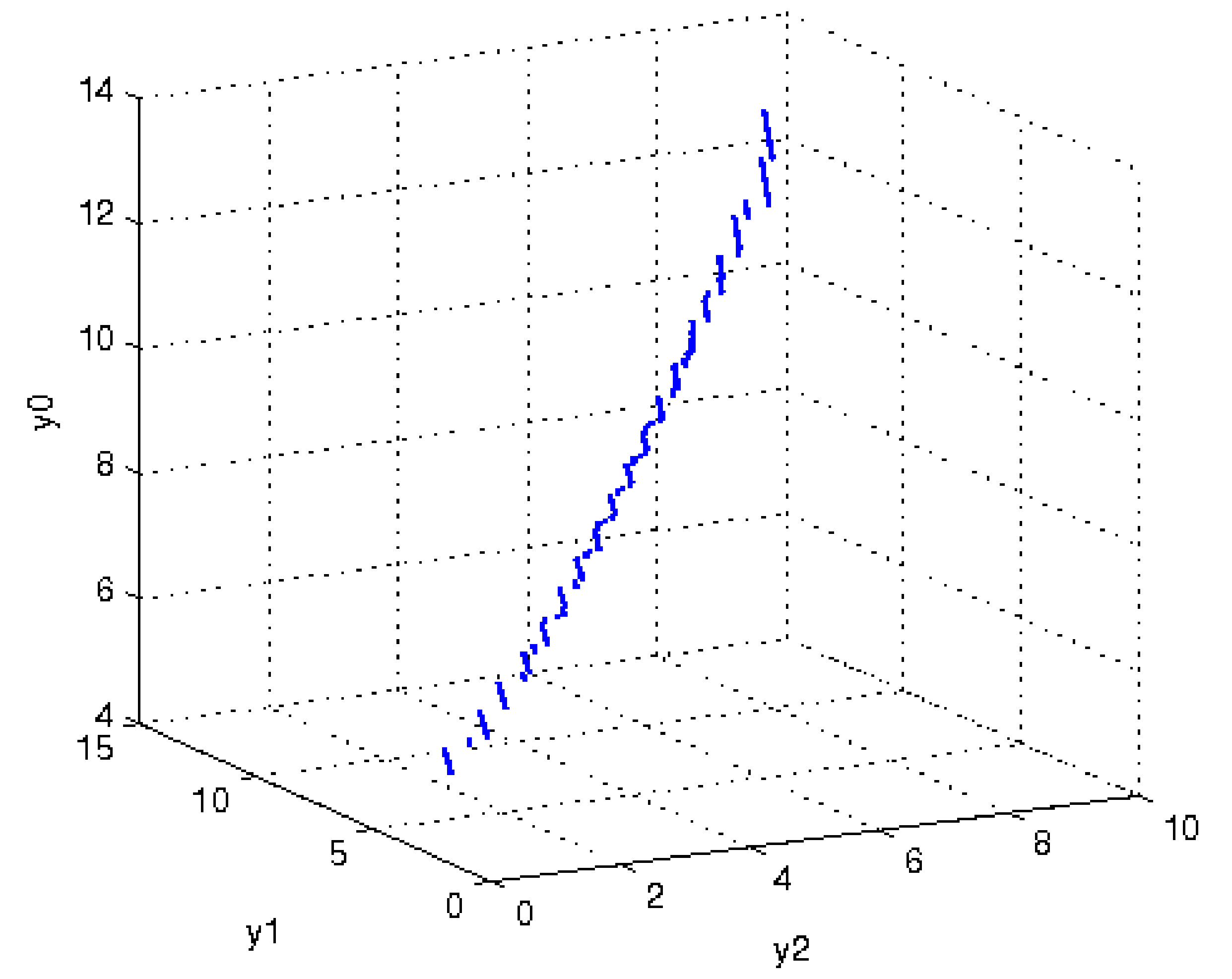}
\caption[Plot of point-cloud]{The point cloud $\mathcal Y$}
  \label{fig:3d-plot}
\end{figure}
Note that the the points of $\mathcal Y$ concentrate around a one-dimensional line segment in $\mathbb R^{d+1}$ whose length is related to the fractal dimension $s$ and whose base point is given by $(D_0,\ldots,D_d)$ (assuming $x_0 = 0$, see figure~\ref{fig:3d-plot}). One way of obtaining estimators $\hat s^{(m)}$ and $\hat D^{(m)}_k$ is to ignore the periodic components of $y_{kj}$ and to fit a line to $\mathcal Y$ so that the sum of squared distances is minimal. 

\begin{proposition}
\label{prop:simplified-regression}
  For a self-similar set $F \in \Sd$, consider the simplified regression model 
  \begin{equation}
  \label{eq:simplified-regression}
       \begin{pmatrix}
  Y_{0j} \\ Y_{1j} \\ \vdots \\ Y_{dj}
  \end{pmatrix}
 =
 \begin{pmatrix}
   D_0 \\ D_1 \\ \vdots \\ D_d
 \end{pmatrix}
 + 
 \begin{pmatrix}
   s \\ s \\ \vdots \\ s
 \end{pmatrix}
  x_j
 +
 \begin{pmatrix}
   \delta_{0j} \\ \delta_{1j} \\ \vdots \\ \delta_{dj} 
  \end{pmatrix}
  \end{equation}
  with variables as defined in (\ref{eq:regression-model}). 
Then the least squares estimate   $(\hat s^{(m)}, \hat D^{(m)}_0,\ldots, \hat D^{(m)}_d)$ for $(s,D_0,\ldots,D_d)$ is given by
  \begin{eqnarray*}
 \hat s^{(\hat m)} &=& \dfrac{\sum_{i=0}^d \left( \sum_{j=1}^m Y_{ij} 
	(x_j - \bar x)\right)} {m(d+1)(\overline{x^2}-\bar{x}^2)} \\
 \hat D_0^{(\hat m)} &=& \overline{Y_0} - \bar{x} \hat s \\
     &\vdots& \\
 \hat D_d^{(\hat m)} &=& \overline{Y_d} - \bar{x} \hat s
\end{eqnarray*}
where $\bar{x} := \frac{1}{m}\sum_{j=1}^m x_j$, $~~\overline{x^2} := \frac{1}{m}\sum_{j=1}^m x_j^2$ and $\overline{Y_i} := \frac{1}{m}\sum_{j=1}^m Y_{ij}$.
  \end{proposition}
\textit{Calculation:}  
 Differentiation of the sum of squared errors 
\[
  S_m^2(\hat s^{(m)},\hat D^{(m)}_0,\ldots,\hat D^{(m)}_d) := \frac{1}{m(d+1)} \sum_{i=0}^d \sum_{j=1}^m (Y_{ij}-(\hat D^{(m)}_i+\hat s^{(m)} x_j))^2 , 
\]
yields the critical point $(\hat s^{(m)}, \hat D^{(m)}_0,\ldots, \hat D^{(m)}_d)$ which has to satisfy
\begin{equation*}
\begin{array}{rrrrrrrcc}
  (d+1) \overline{x^2} \hat s^{(m)} 	&+& \bar{x} \hat D^{(m)}_0 	&+& \ldots 	&+& \bar{x} \hat D^{(m)}_d &=&  \frac{1}{m}\sum_{i=0}^d \sum_{j=1}^m  Y_{ij}x_j \\
   \bar{x} \hat s^{(m)} &+&  \hat D^{(m)}_0 && && &=& \overline{Y_0} \\
\vdots &&  && & \ddots &  &\vdots& \vdots  \\
   \bar{x} \hat s^{(m)} &+& && &&  \hat D^{(m)}_d &=& \overline{Y_d}, \\
\end{array}  
\end{equation*}
a linear system of equations, the solution of which yields the above results. 
$\hfill \square$

If the periodic components $p_k$ of (\ref{eq:regression-model}) are ignored, one can nevertheless arrive at a reasonable estimate $\hat s$ for $s$:

\begin{proposition}
Assume that the regression model (\ref{eq:regression-model}) holds true and let $\mathcal X = \{x_1, x_2, \ldots\}$ be such that $x_j \rightarrow \infty$, and write $\hat s^{(m)}$ for the least-square estimator for $s$ corresponding to the first $m$ data points. Then, with probability 1,
  \[ \hat s^{(m)} \rightarrow s, ~~~ (m\rightarrow \infty). \]
\end{proposition}
\textit{Proof:} This is a consequence of 
\[\frac{Y_{kj}}{x_j} \rightarrow s\] 
almost surely for all $k=0,1,2$ and $j\rightarrow \infty$. $\hfill \square$

However, for a choice of $\mathcal X$ as above, the least-squares estimates $\hat D^{(m)}_k$ that correspond to the first $m$ data points converge to $D_k$ \textit{only for non-arithmetic sets $F$} in general:

\begin{proposition}
  Let $F \in \Sd$ be a non-arithmetic self similar set, and assume the regression model (\ref{eq:regression-model}) holds true. Let $\hat s^{(m)}$ and $\hat D^{(m)}_k$ be the least squares estimators for the simplified regression model (\ref{eq:simplified-regression}) as defined in Prop.\ref{prop:simplified-regression}, associated to the first $m$ data points. Then 
   \[\hat D^{(m)}_k \rightarrow D_k ~~ (m\rightarrow \infty) ~~ \text{ for all } k \in \{0,\ldots,d\}\text{ with probability }1.\]
\end{proposition}
\textit{Proof:} Since $F$ is non-arithmetic, the normalized periodic functions $p_k^0(\eps)$ from (\ref{eq:normalized-periodic}) are constant and equal 1, and thus the simplified regression model (\ref{eq:simplified-regression}) is equivalent to the regression model (\ref{eq:regression-model}). Using standard results on linear regression one sees that the estimators converge almost surely. 
$\hfill \square$

\begin{proposition}
Let $F \in \Sd$ be an $h$-arithmetic self similar set, and let $\hat s^{(m)}$ and $\hat D^{(m)}_k$ be as above. Furthermore, assume that the set of $x$-values $\mathcal X = \{x_1, x_2, \ldots \}$ is uniformly distributed$\mod h$ and that it is converging to $+\infty$. Then 
\[\textnormal{ess sup } p_k^0 > \textnormal{ess inf } p_k^0\] with respect to the Lebesgue-measure on $\R$ implies 
  \[\limsup_{m\rightarrow\infty} \hat D^{(m)}_k < D_k \text{ for all } k\in\{0,\ldots,d\} \text{ with probability 1,} \]
  i.e. $\hat D^{(m)}_k$ is asymptotically biased for each $k$. 
\end{proposition}
A typical set of $x$-values with the above property is $\mathcal X = \{aj: a>0, j \in \mathbb N, \frac{a}{h} \notin \mathbb Q\}$. 

\textit{Idea of proof:} If $F$ is $h$-arithmetic, the estimators $\hat s^{(m)}$ and $\hat D^{(m)}_k$, $k\in\{0,\ldots,d\}$ are minimizers of the sum of squared deviations
\[S_m^2(\hat s^{(m)},\hat D^{(m)}_0,\ldots,\hat D^{(m)}_d)
                  = 
\inf_{(t,E)\in\R\times\R^{d+1}}\frac{1}{m(d+1)} \sum_{i=0}^d \sum_{j=1}^m 
    (Y_{ij}- (E_i + t x_j))^2.\]
Now, as $m\rightarrow \infty$, we can change the order of the infimum and the limit, and  
\begin{eqnarray*}
\lefteqn{\lim_{m\rightarrow\infty}S_m^2(\hat s^{(m)},\hat D^{(m)}_0,\ldots,\hat                 D^{(m)}_d)} \\ 
&=& \inf_{(t,E)\in\R\times\R^{d+1}}
    \frac{1}{d+1} \sum_{i=0}^d
     \lim_{m\rightarrow\infty}\frac{1}{m} \sum_{j=1}^m 
    (Y_{ij}- (E_i + t x_j))^2 \\
&=& \inf_{E\in\R^{d+1}} \frac{1}{d+1} \sum_{i=0}^d
     \lim_{m\rightarrow\infty}\frac{1}{m} \sum_{j=1}^m 
    (Y_{ij}- (E_i + s x_j))^2 \\
&=& \inf_{E\in\R^{d+1}} \frac{1}{d+1} \sum_{i=0}^d
     \lim_{m\rightarrow\infty}\frac{1}{m} \sum_{j=1}^m 
    (D_i-E_i + q_{ij} + \delta_{ij})^2 \\
&=& \inf_{E\in\R^{d+1}} \frac{1}{d+1} \sum_{i=0}^d 
    \left(  (D_i-E_i)^2 + \sigma^2 + \overline{q_i ^2} + 2 (D_i-E_i) \bar q_i \right)\\
&=& \inf_{E\in\R^{d+1}} \frac{1}{d+1} \sum_{i=0}^d 
    \left( \left( (D_i-E_i)+\bar q_i \right)^2 +(\overline{q_i^2}-\bar q_i^2)
          + \sigma^2 \right),
\end{eqnarray*}
a unique infimum which is attained if $E_i = D_i + \bar q_i$ for all $i\in\{0,\ldots,d\}$. Note that in the second equality without loss of generality we assumed $t=s$, as for $t \neq s$ the limit is $\infty$. In the fourth equality we have put $\bar q_i := \int_0^h \log p^0_i(e^{-x})dx$ and $\overline{q_i^2} := \int_0^h \vert\log p_i^0(e^{-x})\vert^2dx$ and used the ergodic theorem. Thus $\hat D_k^{(m)}$ converges to $D_i + \bar q_i$. With the Jensen-inequality, one sees that 
\[\bar q_i = \int_0^h \log p^0_i(e^{-x})dx < \log\int_0^h p^0_i(e^{-x})dx = 0\]
where the inequality is strict because the logarithm is a strictly concave function and because we have assumed $p^0_k(e^{-x})$ not to be constant. 
$\hfill \square$

\paragraph{Estimating Fractal Curvature.} As we have seen, for arithmetic sets $F\in\Sd$ the estimator $\hat D^{(m)}_k$ will systematically underestimate $\log \overline C_k^f(F)$, even if we assume that $x$ can tend to $\infty$, i.e. that we can represent $\eps$-parallel sets for arbitrarily small $\eps$. This is essentially due to oscillations of the rescaled curvatures that average out negatively under the concave function $\log$. But even for non-arithmetic sets, plots of rescaled curvatures may also show oscillations, since these might only vanish for infinitely large $x$. Therefore, a reasonable estimator for the fractal curvature $\overline C_k^f(F)$ must smoothe out these oscillations, and this can be achieved via averaging as in (\ref{eq:average}): Recall that $\overline C_k^f(F)$ equals the average of the rescaled curvature function, i.e.
\begin{equation}
\label{eq:5.10}
\overline C_k^f(F) 
  = \lim_{T\rightarrow \infty}\frac{1}{T}\int_0^T \exp(-x(s-k))C_k(F_{e^{-x}}) dx.
\end{equation}
In practice, however, we can only estimate this integral, as  $g_k(\eps)$ can only be determined for finitely many values of $\eps$, and because we can only use an estimate $\hat s$ for $s$ in the above formula. The integral will be estimated by the area of a histogram whose $j$-th column is centered at $x_j$ and has height equal to $y_{kj}$ for $j\in\{1,\ldots,m\}$. To this end we define the endpoints of the histogram to be
\begin{eqnarray*}
t_0 &:=& x_1 - a, \\
t_j &:=& \frac{x_j + x_{j+1}}{2}, ~~~1 \leq j \leq m-1, \\
t_m &:=& x_m +a
\end{eqnarray*}
where $a$ is half a typical stepwidth between the $x_j$. Thus the limit in (\ref{eq:5.10}) can be estimated by \footnote{Equation (\ref{eq:gamma}) seems to be computationally unstable; instead, the author used the form
\[\frac{1}{t_m - t_0} \sum_{j=1}^m  {\rm sgn}(C_k(F_{e^{-x_j}})) \exp \left(-sx_j+\log\frac{(C_k(F_{e^{-x_j}}))}{e^{-kx_j}}\right) (t_{j}-t_{j-1})\]}
\begin{equation}
\label{eq:gamma}
\Gamma_k^{(m)} := \frac{1}{t_m - t_0} \sum_{j=1}^m \exp (-x_j (s-k))C_k(F_{e^{-x_j}}) (t_{j}-t_{j-1}).
\end{equation}

As seen before, $\vert C^f_k(F) \vert$ is underestimated by $e^{\hat D_k}$ , under the assumption of course that the regression model above is valid. However, note that
\[e^{\hat D_k} \leq \Gamma_k\]
whenever both estimators have the same underlying data set. The size of their relative difference is inversely related to the quality of the fit of a line to the according data set: Usually the fit of a line to $y_{2\cdot}$ is better than to $y_{0\cdot}$, meaning that for $k=2$ the relative difference between the above two estimators lies somewhere between $1\%$ and $2\%$, whereas for $k=0$ it typically ranges from $5\%$ to $12\%$. 

Applied to the sample images from section \ref{sec:images}, $\Gamma_k^{(m)}$ is slightly more accurate than $e^{\hat D_k}$.




\section{Algorithms}
\label{sec:algorithms}
From now on we assume $d=2$ until the end of this thesis, as we describe the implementation of algorithms for the calculation of $\hat s^{(m)}$ and $\hat D^{(m)}_k$ for $k\in\{0,\ldots,d\}$ which we have yet only done for 2D-images. A generalization to higher dimensions is straightforward, it should be noted however that for $k \leq d-2$ the total variational curvatures $C_k^{var}(\Fe)$ may be hard to obtain.

An algorithm for the estimation of $s=\dim F$ and $\overline C_k^f(F)$ for $k\in\{0,1,2\}$ has been implemented by the author and included in the GeoStoch \cite{geostoch} library. It runs as follows:
\begin{enumerate}
\item Set the parameters:
\begin{enumerate}
  \item Set the boolean variables \texttt{useEuler, useBdlength} and \texttt{useArea} to \texttt{true} if the corresponding data $y_{k\cdot}$ for $k=0,1,2$ shall be taken into account in the regression. By default, all three variables are set to the values \texttt{false, true, true} respectively\footnote{Note that a choice of \texttt{false}, \texttt{false}, \texttt{true} will return a dimension estimate by the sausage method from above and an estimate of the average Minkowski-content.}.
  \item If \texttt{quickEvaluate} is set to \texttt{true}, the array of dilation radii $\mathbf R$ is chosen as explained in the paragraph below. If not, the radii of $\mathbf R$ will be chosen uniformly on the logarithmic scale, according to the parameters \texttt{r\_min}, \texttt{step} and \texttt{r\_max}. 
\end{enumerate}

\item Calculate the distance-transform $D_I$ of the binary image $I$ which represents the fractal set $F$. This is an image that records the distances of each pixel to the nearest black pixel:
  \begin{eqnarray*}
  D_I : \mathbf P &\rightarrow	& \R \\
                p &\mapsto	& d\left( p,I^{-1}\left( \{1\}\right) \right)
  \end{eqnarray*}
\item \label{item:break} As $r_j$ runs through $\mathbf R$ in an increasing order, represent the parallel set $F_{r_j}$ by 
  \[\hat F_{r_j} := D_I^{-1}\left( [0,r] \right),\]
and for each $k \in \{0,1,2\}$ estimate $C_k^{var}(F_{r_j})$. For $k\in\{1,2\}$ this is done by calculating the $k$-th Minkowski-functional $C_k(\hat F_{r_j})$;\footnote{and dividing by 2 if $k=1$.} the 0-th variational curvature $C_0^{var}(\hat F_{r_j})$ is estimated by $N + Q = 2N - C_0(\hat F_{r_j})$ where $N$ is the number of connected components of $F_{r_j}$ and $Q$ is the number of holes of $F_{r_j}$. (This goes wrong quite frequently, see p.\pageref{app:C0var} of the appendix.) 

For the calculations of the Minkowski-functionals the algorithm proposed by Klenk, Spodarev and Schmidt (see~\cite{minkdiscr}) and Guderlei (\cite{guderleietal}) is employed, which for each binary image calculates all of the above three functionals simultaneously. Now the values $y_{kj}$ are calculated as in equation~(\ref{eq:def-variables}), and $x_j$ is set to $-\log r_j$. If the boolean variable \texttt{brk} is set to \texttt{true}, then the algorithm will jump to the next step as soon as $N+Q\leq 2$. If not, a warning will be displayed, saying that the dilation radius might already have grown too big for $F_{r_j}$ to reveal any fractal structure.
\item The estimator $\hat s$ from the simplified regression model  (Proposition~\ref{prop:simplified-regression}) is calculated, based on each data set $y_{k\cdot}$ for which the corresponding boolean variable (\texttt{useArea, useBdlength, useEuler}) is set. 
\item The estimator $\Gamma_k^{(m)}$ of the $k$-th total fractal curvature of $F$ is calculated via averaging as in (\ref{eq:gamma}). 
\end{enumerate}
Note that the data set $\{y_{0j}: j\in\{1,\ldots,m\}\}$ should be excluded by setting \texttt{useEuler = false} if $C_0^{var}(F_{r_j})$ cannot be calculated properly, as in that case the $y_{0j}$ do not show growth behaviour suitable for regression. The break condition 
\[\text{``}N + Q \leq 2 \text{''} \] in step~\ref{item:break} helps to avoid the case $y_{0j} = \log 0$.   

\paragraph{The choice of dilation radii.}

The proper choice of dilation radii is essential to the accuracy of the estimates $\hat s$ and $\Gamma_k$, $k \in \{0,1,2\}$. Some thought has to be given to the range and to the distribution of dilation radii:
\begin{itemize}
\item If data for $r \in (0,1)$ are included in the regression, the author has observed that the slope of $y_2$ will be slightly underestimated, while the slopes of $y_1$ and $y_0$ will be overestimated, meaning that the overall accuracy decreases. Thus only radii $\geq 1$ should be used. The default value for $r_{min}$ is 1.261.
\item If the radii are increased too far, $F_r$ will completely lose its fractal character. A reasonable upper bound for the radii can be e.g. 
\[\max\{0.06 * \min\{w, h\},20\}\]
where $w$ is the width and $h$ is the height of the image. This is also the default value.
\item Since the parallel sets $F_r$ are represented by binary images $\hat F_r$, the curvatures $C_k(\hat F_r)$ will only assume discrete values, and $C_k(\hat F_r)$ will be piecewise constant\footnote{Note that on a logarithmic scale, this effect is more prominent for small radii than for big radii, see figure~\ref{fig:discrete-area}. Also observe this effect in figure~\ref{fig:linplot}, where the $y_k$ are piecewise linear as a function of $x=-\log r$, the slope being equal to $k$.}. 
Thus the main amount of information, encoded in $y_{k}$ on an interval where $C_k(F_{e^{-x}})$ is constant, is extracted by the first $x$ that falls into this interval; Further values $x$ that fall into this interval will only have an averaging effect. It seems reasonable to pick radii somewhere out of the middle of these intervals, and not close to their ends which are the points of discontinuity of the $C_k(\hat F_r)$. These points are given e.g. by the optimal area radii from section~\ref{sec:sausage}. 
\end{itemize}

With this in mind, the author suggests two different choices for the set of dilation radii $\mathbf R$, depending on whether or not the run-time of the algorithm is of importance:
\begin{description}
\item [Lower run-time:] The hand-picked radii \texttt{quickRadii} from the class \texttt{OptimalRadii} are used. These are radii yielding the optimal discretized area, starting from $r_1 = r_{min} = 0.5641$ and increasing in multiplicative steps of around $1.5$.
\item [Higher precision:] Choose 
\[\mathbf R = \left\lbrace r_j : r_j = r_{min} \texttt{step}^\texttt{j}, j\in\mathbb N, r_j \leq r_{max} \right\rbrace \]
where the multiplicative stepwidth $\texttt{step}>1$ can be chosen almost arbitrarily close to 1. Default values are $r_{min} = 1.2616$, $r_{max} = \max\{0.06 * \min\{w, h\},20\}$ and $\texttt{step}=1.05$.
\end{description}
 
\paragraph{Problems.}

\begin{figure}[!h]
  \centering
  \includegraphics[height=11cm]{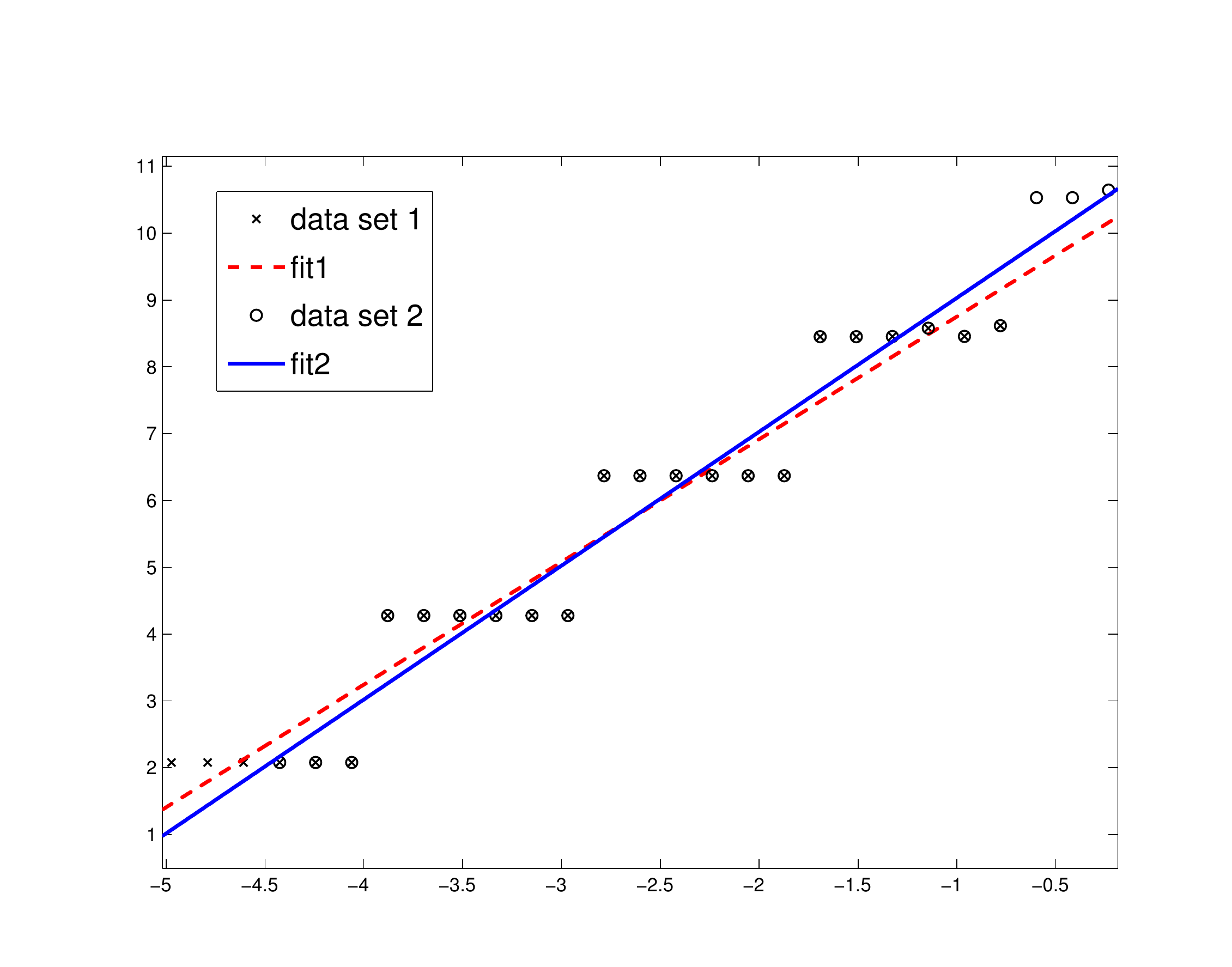}
  \caption[A periodic plot for different data sets]{A plot of $y_0$ for the Sierpi\'{n}ski Carpet SC and for two different data sets. Both data sets can be thought of as being comprised of 4 full periods. However, extrapolating to infinity one sees that the slope of fit 1 is preferable to the slope of fit 2. Indeed, the slopes are 1.87 (fit 1) and 2.04 (fit 2), and the dimension $s$ of the Sierpi\'{n}ski Carpet is 1.89. With the breaking condition ``$N+Q \leq 2$'' one achieves that the lowest stair is not aborted, i.e. too high slopes as in fit 2 above will be avoided.}
  \label{fig:badplot}
\end{figure}
For arithmetic self-similar sets the following problem arises: As the range of radii, for which the values $y_{kj}$ are calculated, is bounded, there will most probably be an incomplete period in the data set. 
But even if there is no incomplete period, the slope can severely deviate from its correct value, depending on where the beginning of a period is put 
(see figure~\ref{fig:badplot}).
The consolating thought here is that examples of arithmetic self-similar sets in nature arise even less often than examples of non-arithmetic self-similar sets.

\chapter{Results}
\label{ch:results}

\section{The Sample Images}
\label{sec:images}

\begin{figure}[!h]
  \centering
  \includegraphics[width=\textwidth]{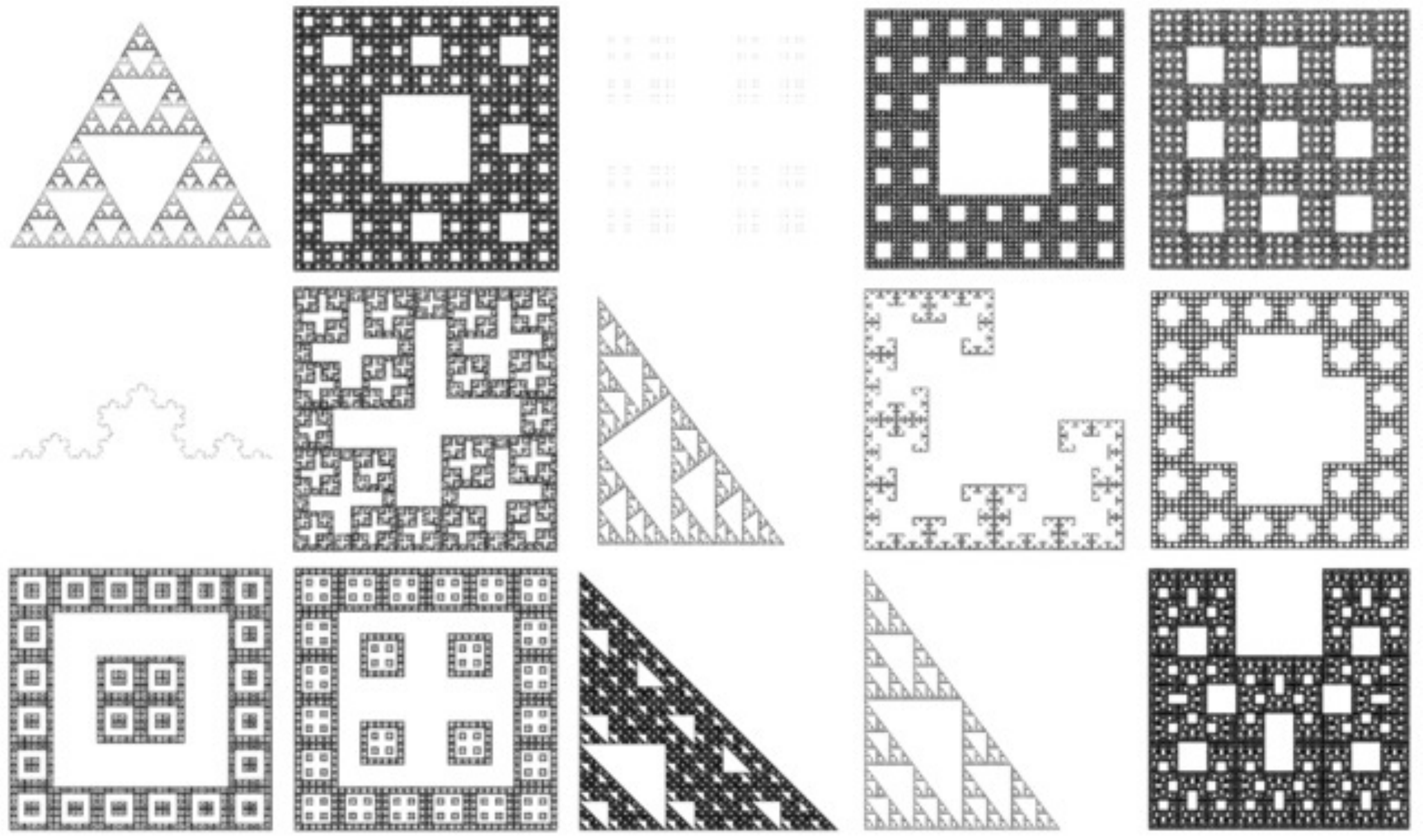}
  \caption{The Sample Images}
  \label{fig:montage}
\end{figure}2

As most canonical examples of self-similar sets as e.g. the Sierpi\'nski gasket, Sierpi\'nski carpet, Cantor dust and Sierpi\'nski tree are arithmetic in the sense of section~\ref{sec:self-similar-sets}, some effort has been put into constructing self-similar sets that are non-arithmetic (and have polyconvex parallel sets). These sets are the non-arithmetic self-similar triangle and square below. 
It turns out that, as $\eps \downarrow 0$, the growth of especially  $C_0^{k,var}(F_\eps$ is more uniform for non-arithmetic sets than for arithmetic sets, which allows a better fit in the regression analysis of chapter~\ref{ch:main}. 

The sets from this chapter all have polyconvex parallel sets except the Cantor dust and the Koch curve. To see that, imagine a dilation of the set in question by a large $\eps$; it will be equal to the $\eps$-parallel set of the set which is the union of the filled-in smaller copies, which is of course polyconvex. Applying theorem~\ref{th-par-sets} one now sees that all parallel sets are polyconvex.

\subsection{Canonical Examples}
\label{sec:canonical}
These are the Sierpi\'nski Gasket, Sierpi\'nski Carpet, Sierpi\'nski Tree, Koch curve and the Cantor Dust.

\paragraph{Sierpi\'nski Gasket SG.}

\begin{figure}[!h]
\begin{center}
  \includegraphics[height=6cm]{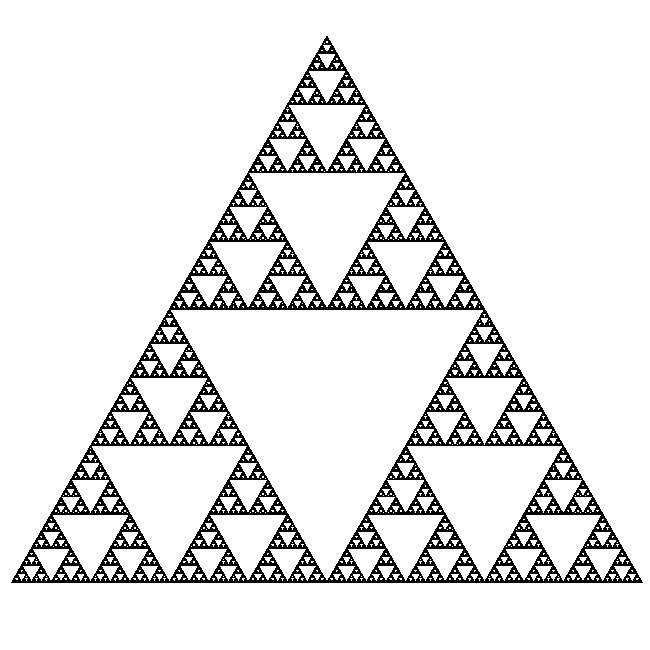}
\end{center}
\caption{The Sierpi\'nski Gasket}
\end{figure}

Here the Iterated Function System consists of three similarities which map (equally oriented) half-sized copies of the gasket towards its vertices. All three similarities have the ratio $\frac{1}{2}$, and thus the gasket is $\log 2$-arithmetic, and its dimension $s$ is the solution of $3 \left(\frac{1}{2}\right)^s = 1 \Leftrightarrow s=\frac{\log 3}{\log2} \approx 1.585... $. It has polyconvex parallel sets. The 0-th variational measure $C_0^{var}(SG_\eps)$ of the parallel sets can be calculated accurately by an image analyser, being equal to $1-C_0(SG_\eps)$

The author recalculated the fractal curvatures of SG in his dissertation (\cite{diss_winter}). If the base has length $1$, the correct values are:
\begin{eqnarray*}
  C_0^f(SG) &\approx& -0.042345... \\
  C_1^f(SG) &\approx& 0.37615... \\
  C_2^f(SG) &\approx& 1.81
\end{eqnarray*}

\paragraph{Sierpi\'nski Carpet SC.}

\begin{figure}[h!]
\begin{center}
  \includegraphics[height=6cm]{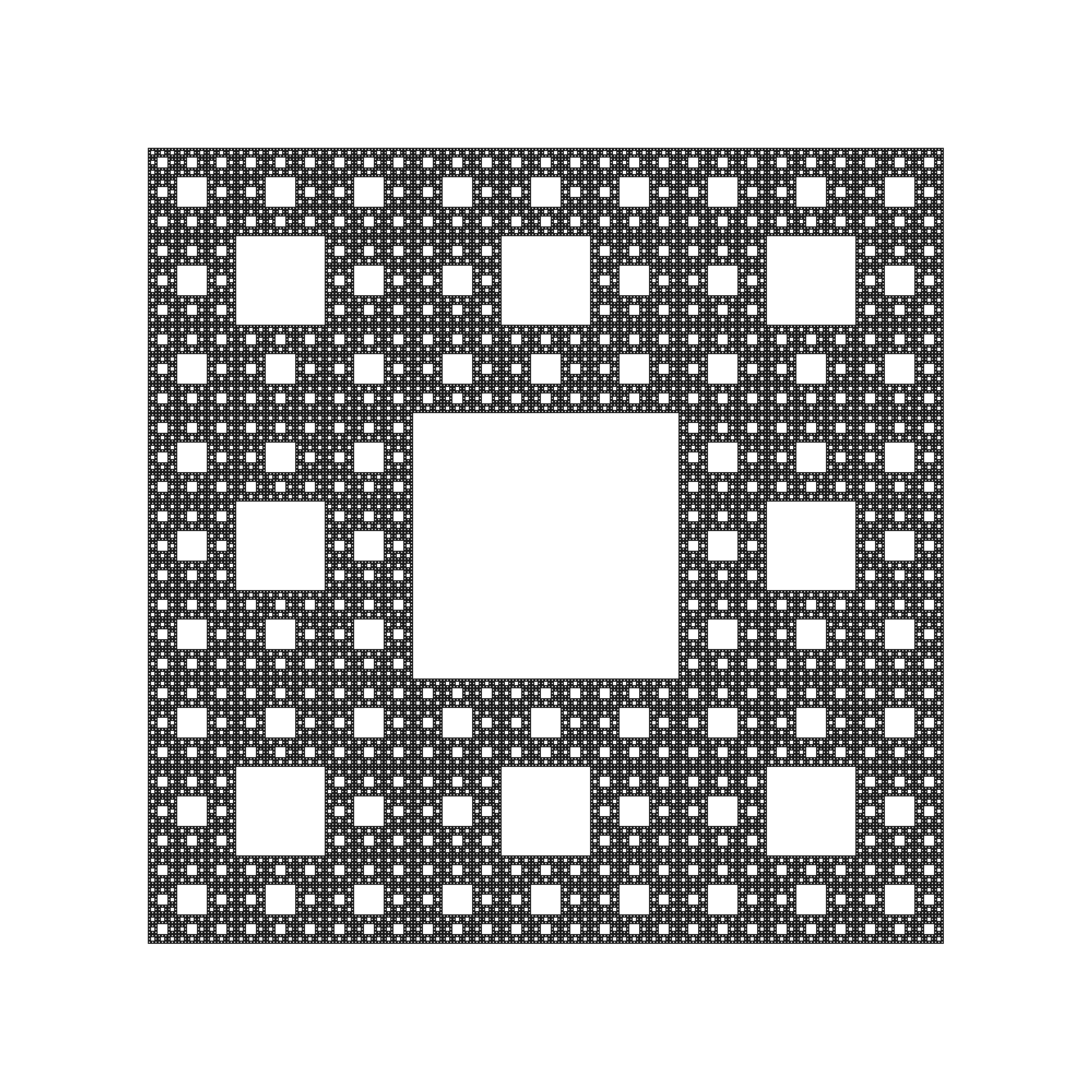}
\end{center}
\caption{The Sierpi\'nski Carpet}
\end{figure}

The Sierpi\'nski Carpet SC 
is the attractor of eight similarities of ratio $\frac{1}{3}$ and thus has dimension $s = \frac{\log 8}{\log 3} \approx 1.893...$.  It is log 3-arithmetic and its parallel sets are polyconvex. 

On a scale where the base line has length 1, the average total fractal curvatures are 
\begin{eqnarray*}
  C_0^f(SC) &\approx& -0.0162 \\
  C_1^f(SC) &\approx&  0.0725 \\
  C_2^f(SC) &\approx&  1.352
\end{eqnarray*}
as computed by Winter, and $C_1^f(SC)$ was rechecked by the author.

\paragraph{Sierpi\'nski Tree ST.}

\begin{figure}[h!]
\begin{center}
  \includegraphics[height=6cm]{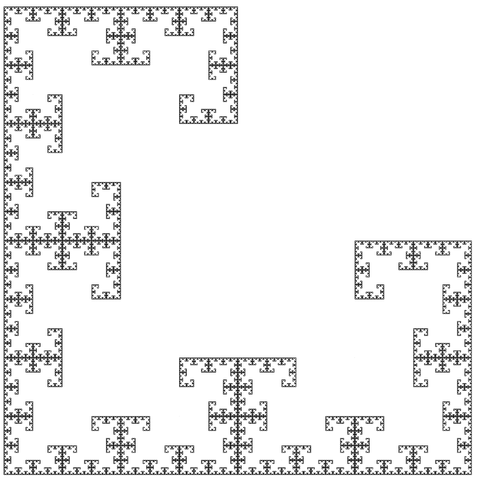}
\end{center}
\caption{The Sierpinski Tree}
\end{figure}

ST is the attractor of three similarities, all of ratio $\frac{1}{2}$: 
\begin{enumerate}
  \item rotation by $-\frac{\pi}{2}$ and placement in upper left corner
  \item placement in lower left corner
  \item rotation by $+\frac{\pi}{2}$ and placement in lower right corner.
\end{enumerate}

Similarly to SG, ST has dimension $s = \frac{\log3}{\log2} \approx 1.585... $ and is $\log 2$-arithmetic. It has polyconvex parallel sets. But unlike the Sierpi\'nski gasket, here we have $X_0(ST) = 0$ and in fact the $0$-th fractal curvature $C_0^f(ST)=0$. Furthermore, $C_0^{var}(ST_\eps)$ cannot be measured correctly, and thus in the regression process of chapter~\ref{ch:main} \texttt{useEuler} has to be set to \texttt{false}.

\paragraph{Cantor Dust CD.}

\begin{figure}[h!]
  \centering
  \includegraphics[width=6cm]{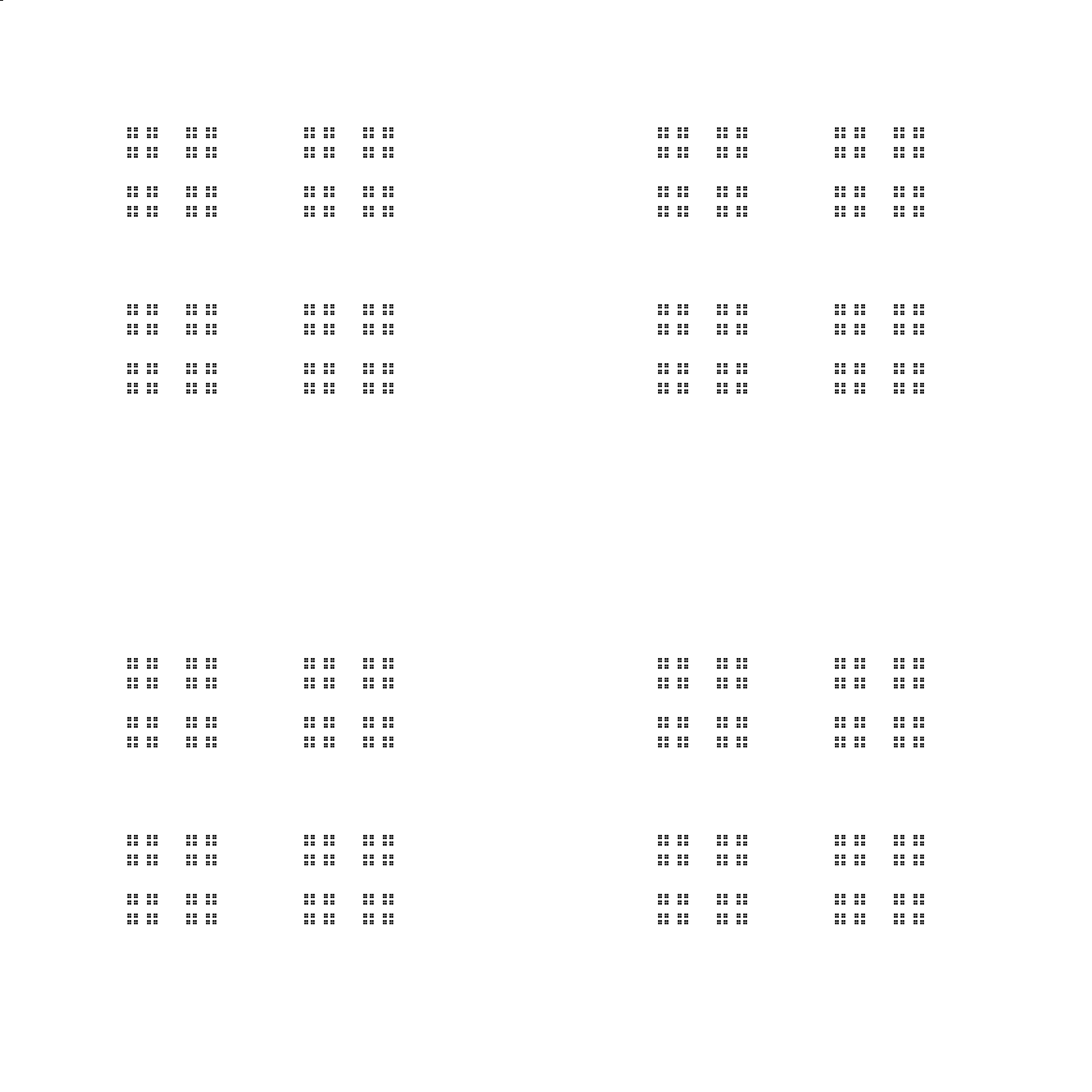}
  \caption{The middle third Cantor dust.}
  \label{fig:CD}
\end{figure}

Recall the Cantor Dust in one dimension: In each step of the construction, the middle third of the unit intervall is erased. In the two-dimensional case, 
in each step a cross of thickness one third is erased. 

The corresponding Iterated Function System consists of 4 similarities of ratio $\frac{1}{3}$, mapping CD in the upper left, upper right, lower left and lower right corner. As all similarities have ratio $\frac{1}{3}$, CD is $\log 3$-arithmetic, and $\dim(\textnormal{CD}) = \frac{\log 4}{\log 3} \approx 1.2619$. This is an example of a self-similar set with no polyconvex parallel sets, so strictly speaking the theory of chapter~\ref{ch:fractal_curvature} does not apply. However, there does seem to be the correct scaling behaviour $s_2 = s-2$, $s_1 = s-1$, and dimension estimates based on the regression method of chapter~\ref{ch:main} with 
\texttt{(useEuler, useBdlength, useArea) = (false, true, true)} are quite exact. 

\paragraph{Koch curve KC.}

\begin{figure}[!h]
  \centering
  \includegraphics[width=9cm]{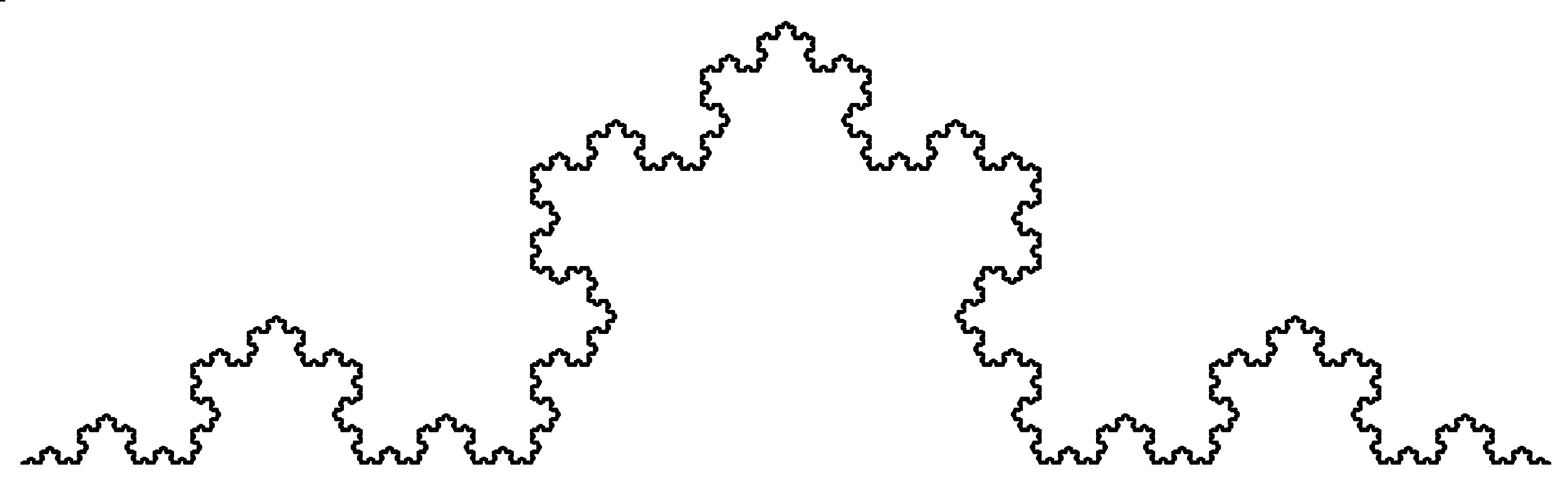}
  \caption{The Koch-curve}
  \label{fig:kochcurve}
\end{figure}

The Koch-curve KC shares the dimension-, arithmeticity- and polyconvexity- properties of the Cantor dust CD, and also seems to have the same scaling behaviour for $s_2$ and $s_1$ with equally good dimension estimates using the above parameters. It is noteworthy that the local dimension estimates are too high for this set.

\subsection{Further Examples}

\paragraph{Modified Sierpi\'nski Carpet MC.}
\begin{figure}[ht]
\begin{center}
  \includegraphics[height=6cm]{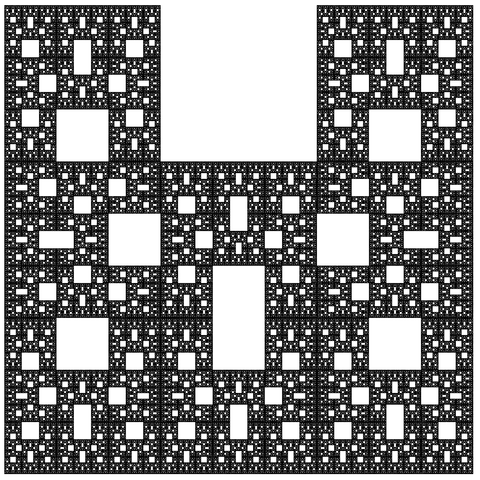}
\end{center}
\caption{The modified Sierpi\'nski Carpet}
\end{figure}

The iterated function system of the modified carpet 
consists of 8 similarities of ratio $1 \over 3$, just as for the original carpet, so we have the same dimension, $s = \frac{\log 8}{\log 3} \approx 1.893...$, and $\log 3$ arithmeticity. There are some differences in the positioning of the smaller copies (there is a copy placed in the middle now) and some of the smaller copies are rotated by multiples of $\frac{\pi}{4}$. The parallel sets are polyconvex. Fractal curvatures are taken from \cite{diss_winter}.

\begin{eqnarray*}
  C^f_0({\rm MC}) & \approx & -0.014...\\
  C^f_1({\rm MC}) & \approx &  0.0720...\\
  C^f_2({\rm MC}) & \approx &  1.344...
\end{eqnarray*}

\paragraph{The Tripet.}

\begin{figure}[!h]
\centering
  \subfigure[The Tripet (half triangle, half carpet)]{
 	\includegraphics[height=6cm]{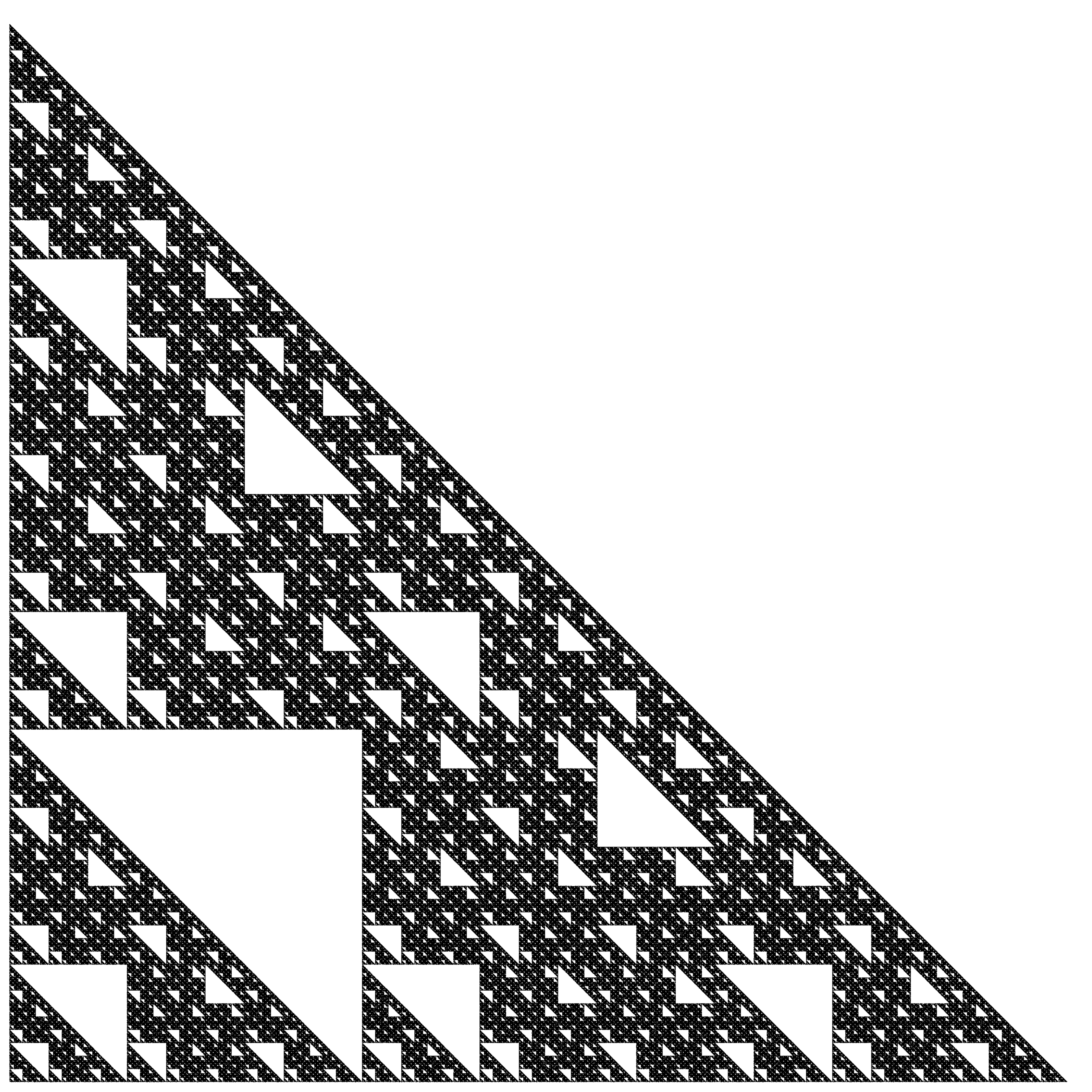}
  	\label{fig:tricarpet}}
  \hspace{1cm}
  \subfigure[The similarities generating the Tripet. Thanks go to Steffen Winter for 
	this idea and picture.]{
	\includegraphics[height=6cm]{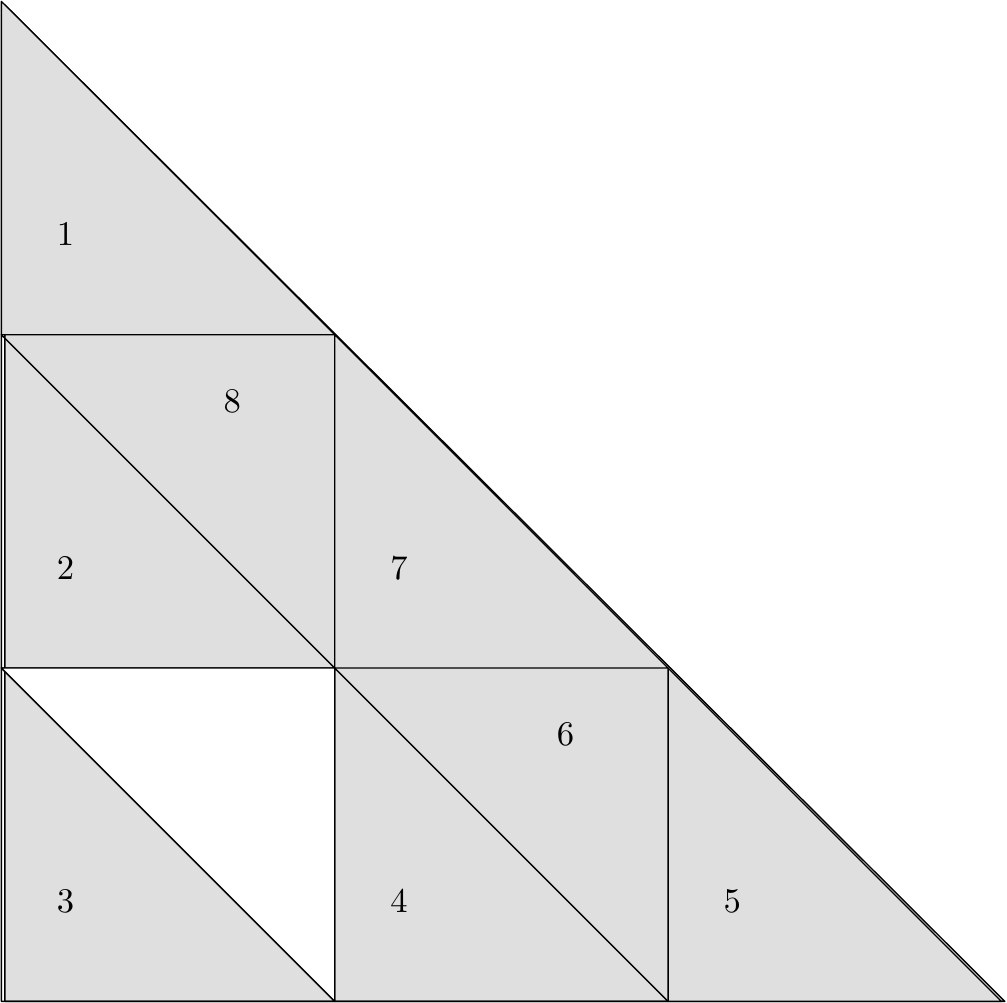}
  	\label{fig:tricarpet-similarities}
	}
\caption{The Tripet.}
\end{figure}

Like the Sierpi\'nski Carpet and the modified Sierpi\'nski Carpet, the tripet is the attractor resulting from 8 similarities of ratio $\frac{1}{3}$ and thus has dimension $s = \frac{\log 8}{\log 3} \approx 1.893...$. It is also $\log 3$-arithmetic. The specialty about this set is that standard methods seem to underestimate the dimension.

\paragraph{Non-arithmetic self-similar triangle $\triangle$.}

\begin{figure}[h!]
  \centering
  \includegraphics[height=8cm]{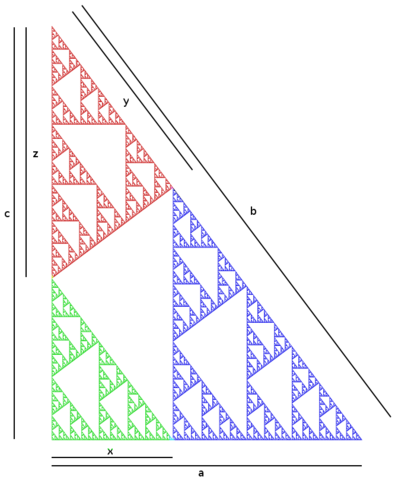}
  \caption{The non-arithmetic self-similar triangle $\triangle$}
  \label{fig:triangle}
\end{figure}

Consider the self-similar triangle $\triangle$ from the above figure 
with a right angle at the bottom left corner. Two smaller, positively oriented copies (green and blue) of the triangle are placed in the bottom corners, while the similarity that maps the triangle to the upper vertice (red) is a scaled reflection and thus negatively oriented. A few calculations show that, in order for the vertices of the smaller copies of the triangle to meet as they do in figure~\ref{fig:triangle}, necessarily
\[
 x=\frac{ac^2}{b^2+c^2}, \hspace{1cm} y=\frac{bc^2}{b^2+c^2}, \hspace{1cm} z=\frac{b^2c}{b^2+c^2},
\]
and thus the above triangle is unique up to similarity transformations. To fix the size of $\triangle$ we choose $a=0.6$, $b=1$ and $c=0.8$, then the similarity ratios $r_2$ for the reflection and $r_1$, $r_3$ for the remaining two similarities become
\[
 r_1=\frac{25}{41}, \hspace{1cm} r_2 = \frac{20}{41}, \hspace{1cm} r_3 = \frac{16}{41}.
\]
The similarity dimension $s$ which is the solution to
\[r_1^s + r_2^s + r_3^s = 1\]
was computed numerically with Matlab, yielding the result 
\[s \approx 1.5882...\]
The set $\triangle$ is ``properly non-arithmetic'', in the sense that 
\[i \neq j \Rightarrow \frac{\log r_i}{\log r_j} \notin \mathbb Q.\]
To see that e.g. $\frac{\log r_1}{\log r_2} \notin \mathbb Q$, suppose that 
\[\frac{\log r_1}{\log r_2}=\frac{p}{q} \text { for some }p \in \mathbb Z, q \in \mathbb \N,\]
then (using the uniqueness of prime number decompositions)
\[r_1^q = r_2^p \Rightarrow  41^{q-p} \times 5^{-q} \times 2^{4p-2q} = 1 \Rightarrow  p=q=0,\]
a contradiction to $q \in \mathbb N$.

\begin{figure}[!h]
  \centering
  \includegraphics[width=9cm]{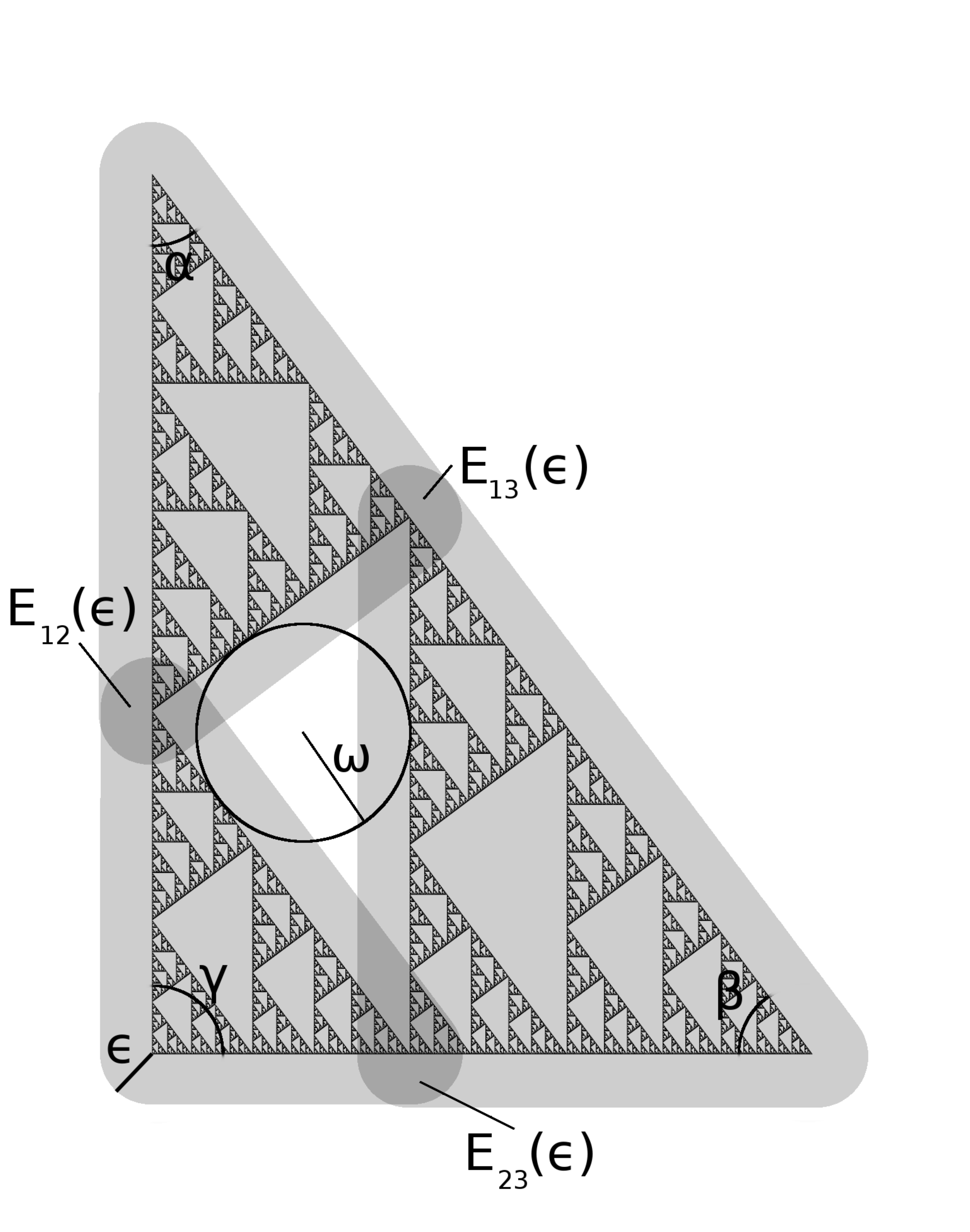}
  \caption{The parallel sets of the (smaller copies of) $\triangle$.}
  \label{fig:triangle-dil}
\end{figure}

As $\triangle$ is a (non-arithmetic) self-similar set with polyconvex parallel sets, its fractal curvatures $C^f_k(\triangle)$ are well-defined. The tedious task of computing them via theorem~\ref{th:curvature-formula} is described in the appendix on page~\pageref{app:triangle}; here only the results shall be given. 

\begin{eqnarray*}
  C^f_0(\triangle) & \approx & -0.023459108...\\
  C^f_1(\triangle) & \approx &  0.239312913...\\
  C^f_2(\triangle) & \approx &  1.162171558....,
\end{eqnarray*}
on a scale where the hypotenuse has length 1.

\paragraph{The sheared Sierpinski Gasket SSG}

\begin{figure}[!h]
  \centering
  \includegraphics[height=8cm]{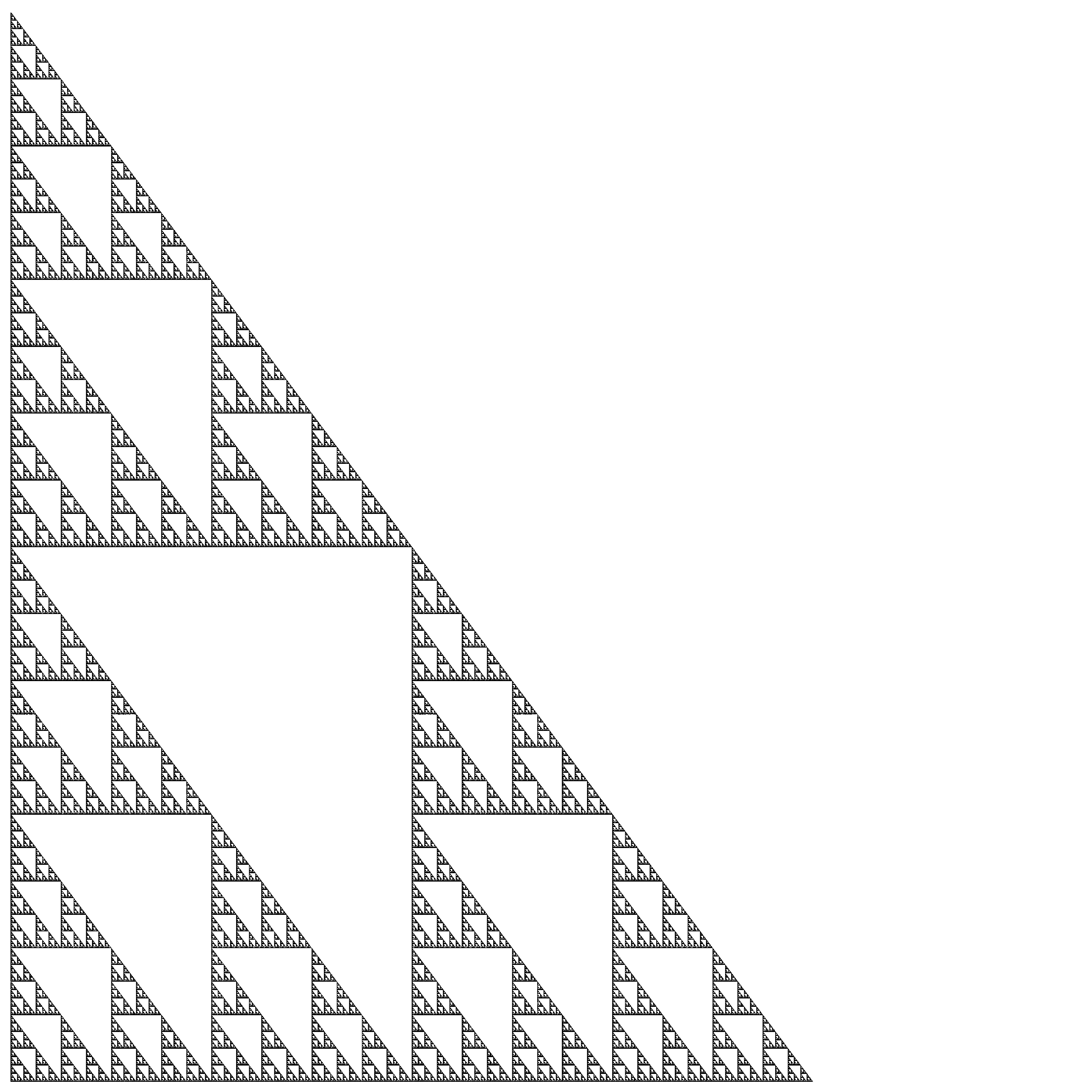}
  \caption{The sheared Sierpi\'nski Gasket}
  \label{fig:gaske_schief}
\end{figure}

The Sierpi\'nski Gasket can be ``sheared'' in such a way that it has the same convex hull as the self-similar triangle from above. The dimension stays the same: $\dim(\textnormal{SSG}) = \dim(\textnormal{SG}) = \frac{\log 3}{\log 2} \approx 1.558...$. For a comparison of measured total curvatures see pages~\pageref{tab:curvature-estimates} and \pageref{tab:specific-curvatures}.

\paragraph{The sets M1, M2 and M3.}
\begin{figure}[h!]
\label{fig:M123}
\begin{center}
  \includegraphics[width=\textwidth,bb=14 14 1719 583]{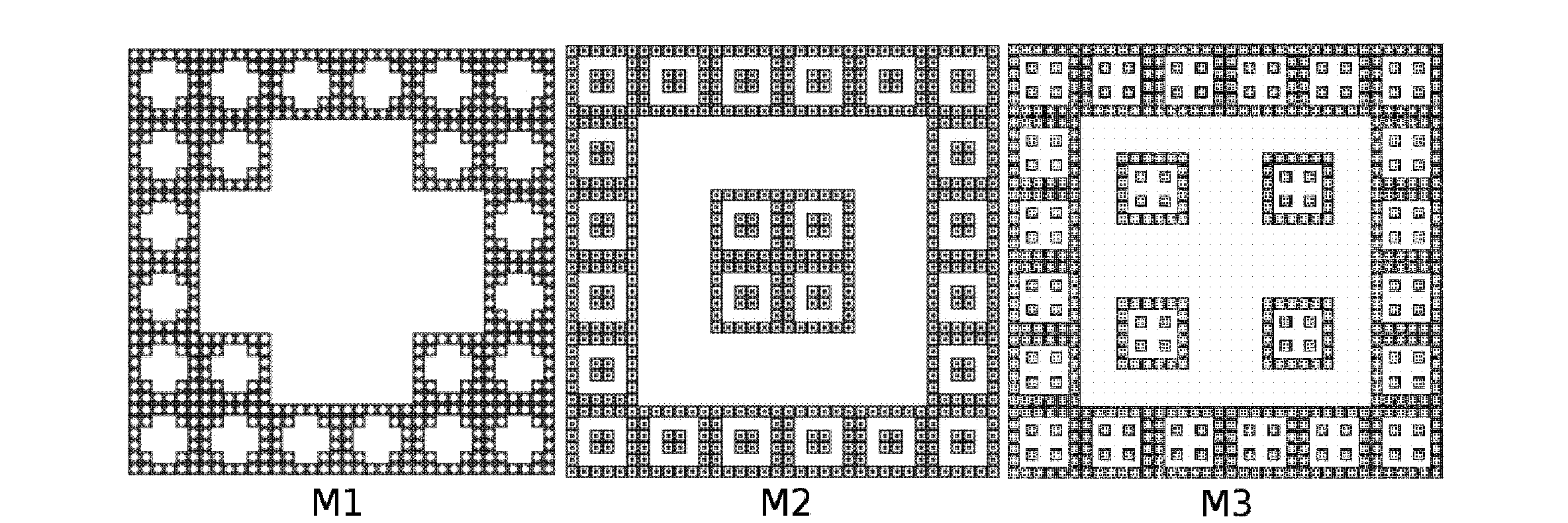}
\end{center}
\caption{The sets M1, M2 and M3}
\end{figure}

The three polyconvex self-similar sets from figure~\ref{fig:M123} are the attractors of $24$ similarities of ratio $\frac{1}{6}$ and thus all have the dimension 
\[s = \frac{\log24}{\log6} \approx 1.7737...\]
Since the similarity ratios are relatively low, there are not so many iterations possible before the maximum resolution is reached, and naturally this means that discretization errors will be higher compared to the other sets (except the window and gate below which similarity ratios $\frac{1}{7}$). The ``gaps-sizes'' are decreasing from M1 to M3, which is reflected in increasing estimates of dimension and of curvatures. 

\paragraph{The Window and Gate.}

\begin{figure}[h]
\centering
  \subfigure[The Window]{
  	\includegraphics[width=5cm]{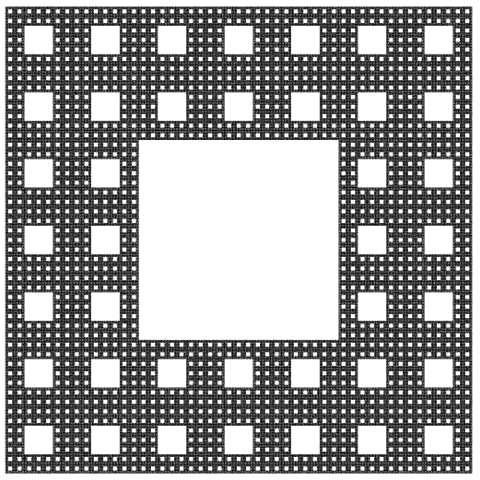}
  }
  \subfigure[The Gate]{
	\includegraphics[width=5cm]{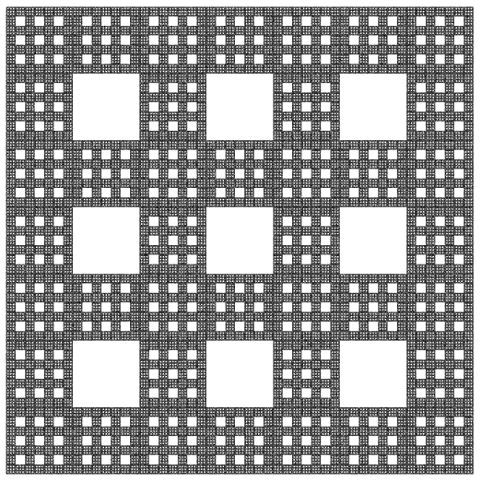}
  }
  \caption{Two versions of the Sierpi\'nski Carpet}
\end{figure}

The dimension of these two sets is $\frac{\log 40}{\log 7} \approx 1.8957...$ as there are 40 similarities of ratio $1 \over 7$. The author calculated their curvatures in calculations which were even longer than those for the set $\triangle$ and which are not given here for the sake of brevity. The results are:
\begin{center}
\begin{tabular}{rcrrcr}
  $C^f_0(W)$ & $\approx$ & $-0.0146171712902$ & $C^f_0(G)$& $\approx$ &$-0.0163916537451$	\\
  $C^f_1(W)$ & $\approx$ & $ 0.0652764265706$ & $C^f_1(G)$& $\approx$ &$ 0.0732007965716$	\\
  $C^f_2(W)$ & $\approx$ & $ 1.251813666054$ & $C^f_2(G)$& $\approx$ &$1.403780236274$  
\end{tabular}
\end{center}
On the pixelscale, i.e. where the base has length 2920, the fractal curvatures rescale with 
$2920^{\frac{\log 40}{\log 7}}$, yielding
\begin{center}
\begin{tabular}{rcrrcr}
  $C^f_0(W)$ & $\approx$ & $-54228$   &~~~ $C^f_0(G)$& $\approx$ &$-60811  $	\\
  $C^f_1(W)$ & $\approx$ & $ 242166$  &~~~ $C^f_1(G)$& $\approx$ &$ 271565 $	\\
  $C^f_2(W)$ & $\approx$ & $ 4644054$ &~~~ $C^f_2(G)$& $\approx$ &$ 5207828$  
\end{tabular}
\end{center}
It is noteworthy that the specific fractal curvatures (i.e. $0$th and $1$st curvatures if the $2$nd curvatures are normalized to $1$, see section~\ref{sec:beyond}) are identical up to 14 significant decimal figures:
\begin{eqnarray*}
  \Xi_0(W) = \Xi_0(G) & \approx & -0.011676794 \\
  \Xi_1(W) = \Xi_1(G) & \approx &  0.052145481.
\end{eqnarray*}

\paragraph{Non-arithmetic self-similar square $\square$.}

\begin{figure}[h!]
 \centering
  \includegraphics[width=6cm]{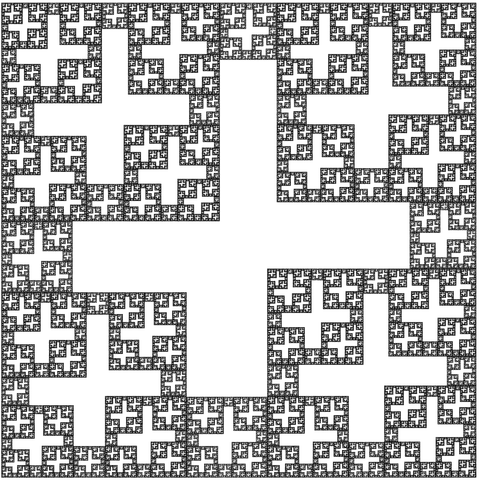}
  \includegraphics[width=6cm]{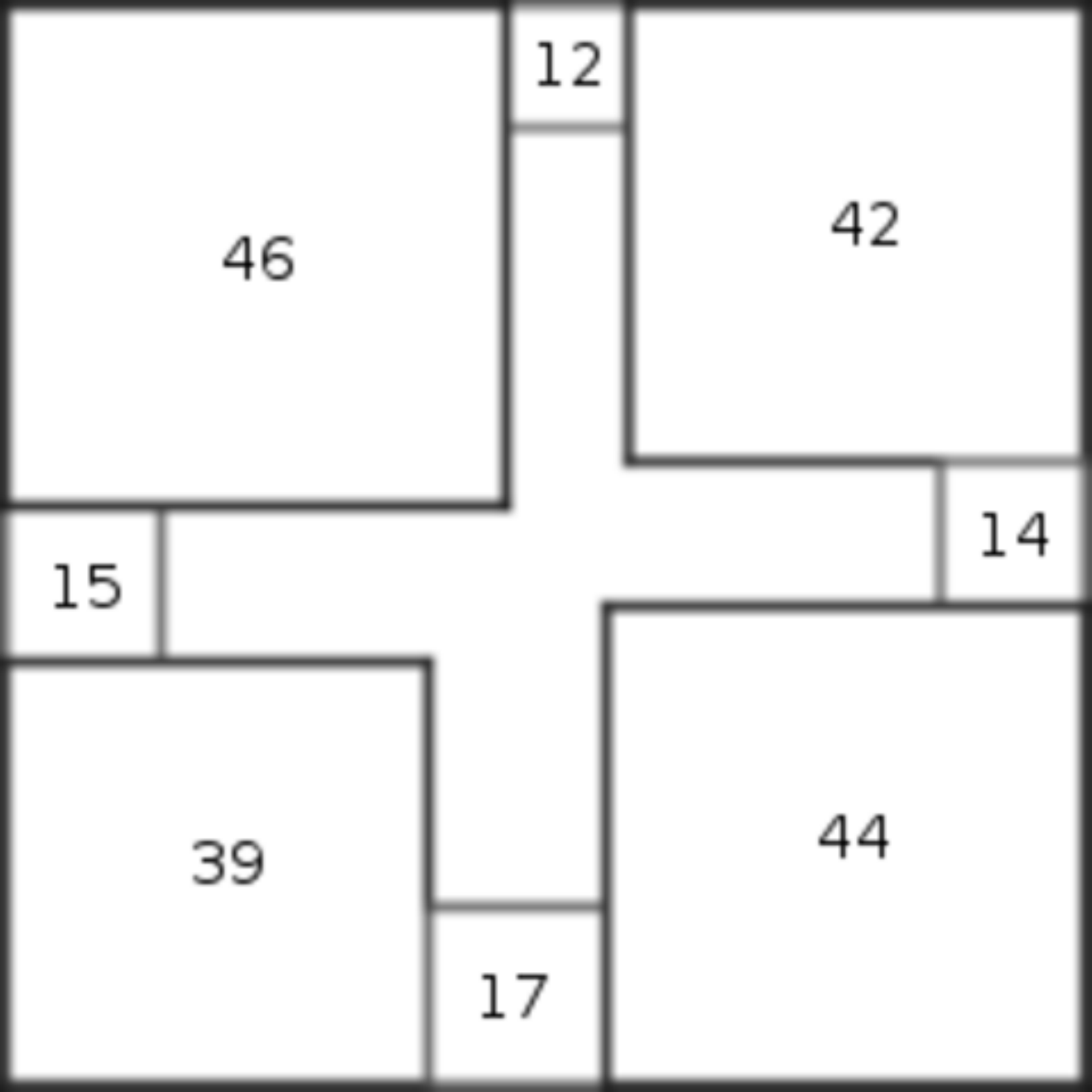}
 \caption[Self-similar Square]{Left: a non-arithmetic self-similar square. Right: its subdivision by similarities. The big square has base-length 100, the numbers show the base-lengths of the small squares.}
 \label{fig:non-arithmetic-square}
\end{figure}

Figure \ref{fig:non-arithmetic-square} shows a self-similar square that consists of eight small copies of itself like the Siepi\'nski Carpet, but the sizes of the copies vary according to the numbers on the right. Computing the similarity dimension numerically yields the value $s=1.7937...$. Note that $\square$ is non-arithmetic, and that $C_0^{var}(\square)$ cannot be measured properly. 

\section{Estimates of Dimension}

The estimator $\hat s$ of the fractal dimension which is based on the multiple regression algorithm from chapter~\ref{ch:main} will be compared to the already known estimators, a short description of which can be found in chapter~\ref{ch:review}. The sample sets from the previous section will be binary images in the versions ``large'' ($3000 \times 3000$ pixels) and ``small'' ($250 \times 250$ pixels). For the small versions, cutouts of the large images have been chosen, and care has been taken that the cutouts do not disappear completely in a hole of the big fractal. Tables \ref{tab:big} and \ref{tab:small} show comparative results. 

For the box-counting estimates we chose the program
``FracLac'' (\cite{fraclac}), a plug-in of the open-source image analysis software ``ImageJ'' \cite{imagej}. The local dimension estimates were obtained from an algorithm that the author has implemented himself (see section~\ref{sec:localdim}). It has been sorted into the ``GeoStoch''-library (\cite{geostoch}) as the class \texttt{LocalDimension}. 
The default parameters for this method were used in the below measurements: 1050 test points in which the local dimension is estimated, and $0.8\times$ the number of black pixels has been used as the number of sample points for the nearest neighbour characteristics. Unlike in \cite{stoyan1994frs}, we do not choose the modal value of the histogram of dimension estimates as the estimate for the global dimension, but the arithmetic mean of the 1050 estimates.

In tables \ref{tab:big} and \ref{tab:small}, by ``area'', ``b'dary'' and ``euler'' we mean that in the method of chapter~\ref{ch:main} the regression has been based exclusively on one of the according data sets $y_2$, $y_1$ and $y_0$ (see chapter~\ref{ch:main}). Note that the case ``area'' corresponds to the sausage method from section~\ref{sec:sausage}. By joint2 we mean the regression estimate based on both $y_2$ and $y_1$, and by joint3 the estimate based on all three data sets. 

If the set $F$ does not belong to $\Sd$, there is no theoretical foundation for $y_1$ and $y_0$ to be used in a regression estimate for the dimension (yet). So, strictly speaking, for the Koch-curve and for the Cantor dust, not only euler and joint3 should be marked ``N/A'' but also joint2 and b'dary. However, $y_1$ seems to admit good fits for \textit{all} self-similar sets, and thus the estimates have been kept in the table. Also note that the Sierpi\'nski tree is a member of $\Sd$, yet still $y_0$ cannot be used for a regression here as it does not represent the total variational measure $C_0^{var}(F_{x_j})$ very well.

\begin{table}[!h]
\label{tab:big}
\centering
\begin{tabular}{c||*{7}{p{1.1cm}|}c}
	\multicolumn{8}{c}{images of size $3000\times3000$} \\ \hline
& theor.& box 	& local	& area	& b'dary& euler	& joint2& joint3 \\ \hline\hline 
\includegraphics[height=18pt]{CD.png} 
& 1.262 & 1.254	& 1.281	& 1.211	&(1.274)&(1.394)&(1.242)&(1.291)\\ \hline
\includegraphics[height=18pt]{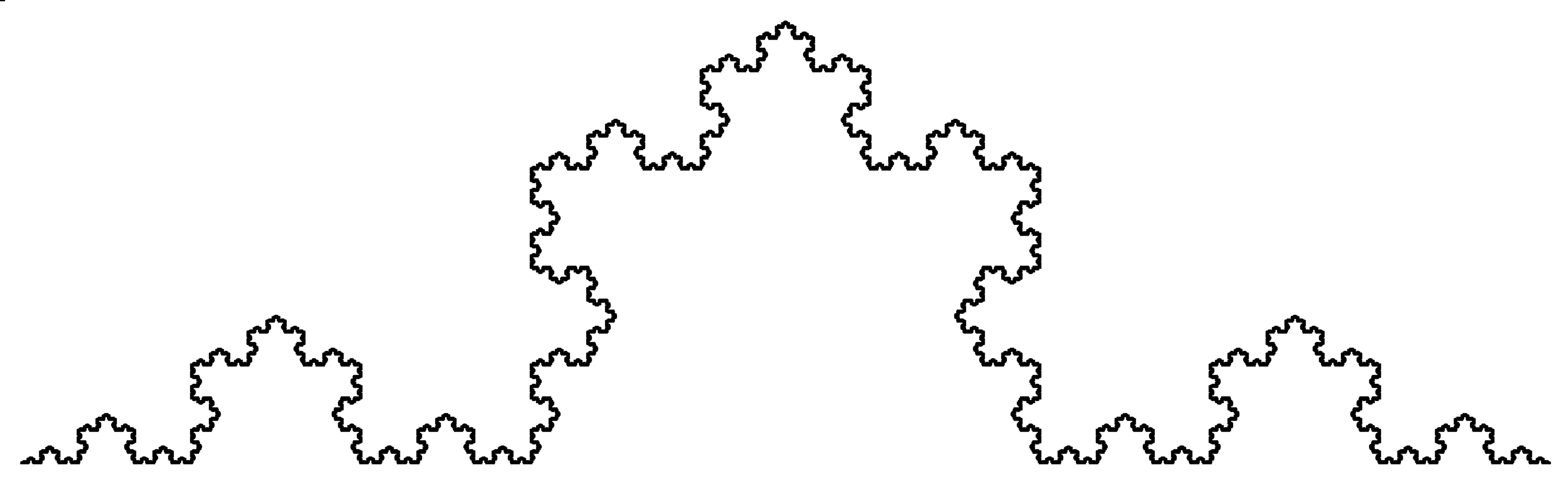}
& 1.262	& 1.270	& 1.354	& 1.268 &(1.237)& N/A	&(1.252)& N/A   \\ \hline
\includegraphics[height=18pt]{tree.png} 
& 1.585 & 1.564	& 1.617 & 1.548	& 1.508	& N/A	& 1.527	& N/A	\\ \hline
\includegraphics[height=18pt]{gasket.png}
& 1.585	& 1.539	& 1.568	& 1.585	& 1.555	& 1.607	& 1.570	& 1.582 \\ \hline
\includegraphics[height=18pt]{gasket_schief.png}
& 1.585	& 1.588	& 1.563	& 1.583	& 1.552	& 1.557	& 1.574	& 1.564 \\ \hline
\includegraphics[height=18pt]{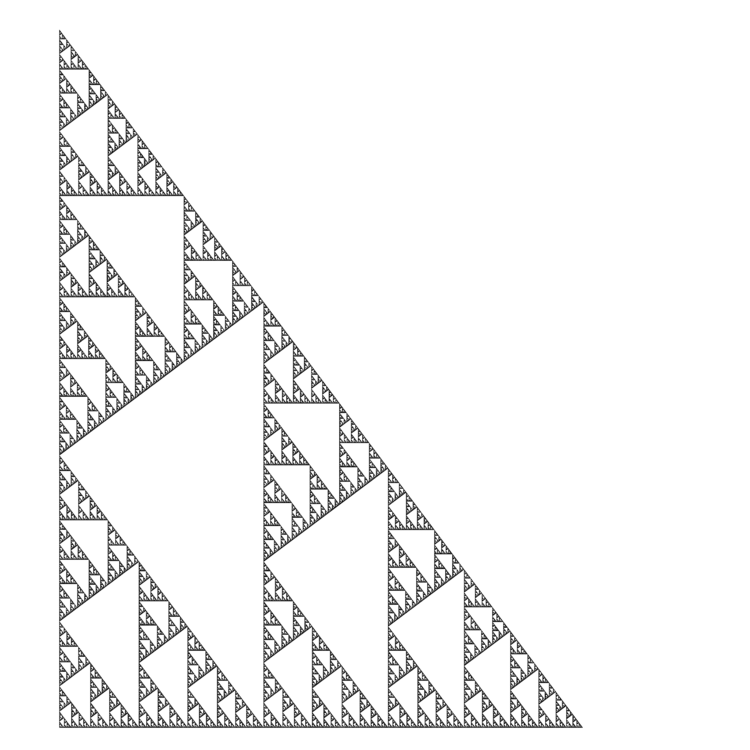} 
& 1.588 & 1.573	& 1.555	& 1.585	& 1.553	& 1.571	& 1.569	& 1.569 \\ \hline
\includegraphics[height=18pt]{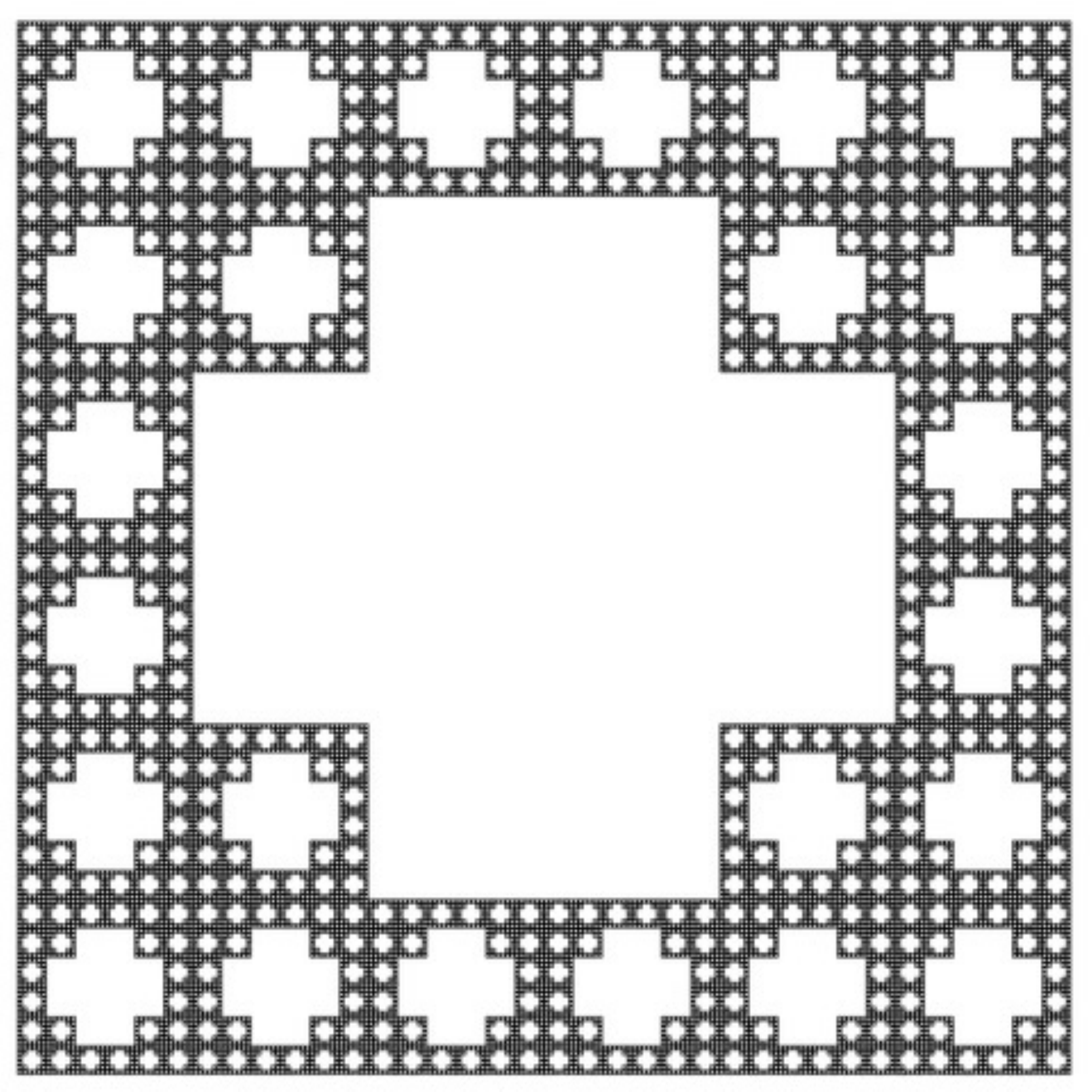}
& 1.774	& 1.706	& 1.760	& 1.732	& 1.656	& 1.635	& 1.694	& 1.715	\\ \hline
\includegraphics[height=18pt]{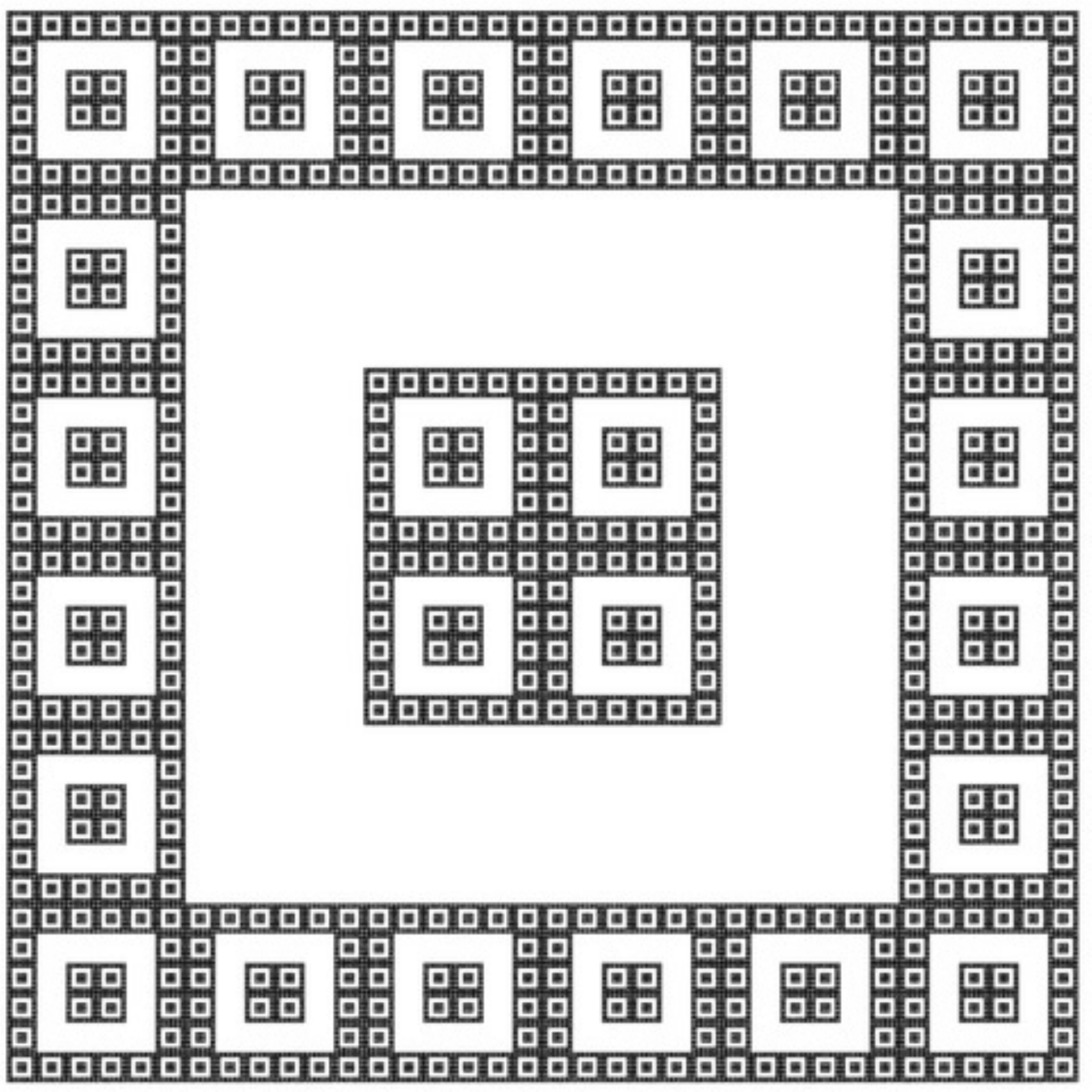} 
& 1.774	& 1.730	& 1.791 & 1.744	& 1.642	& 1.758	& 1.692	& 1.670	\\ \hline
\includegraphics[height=18pt]{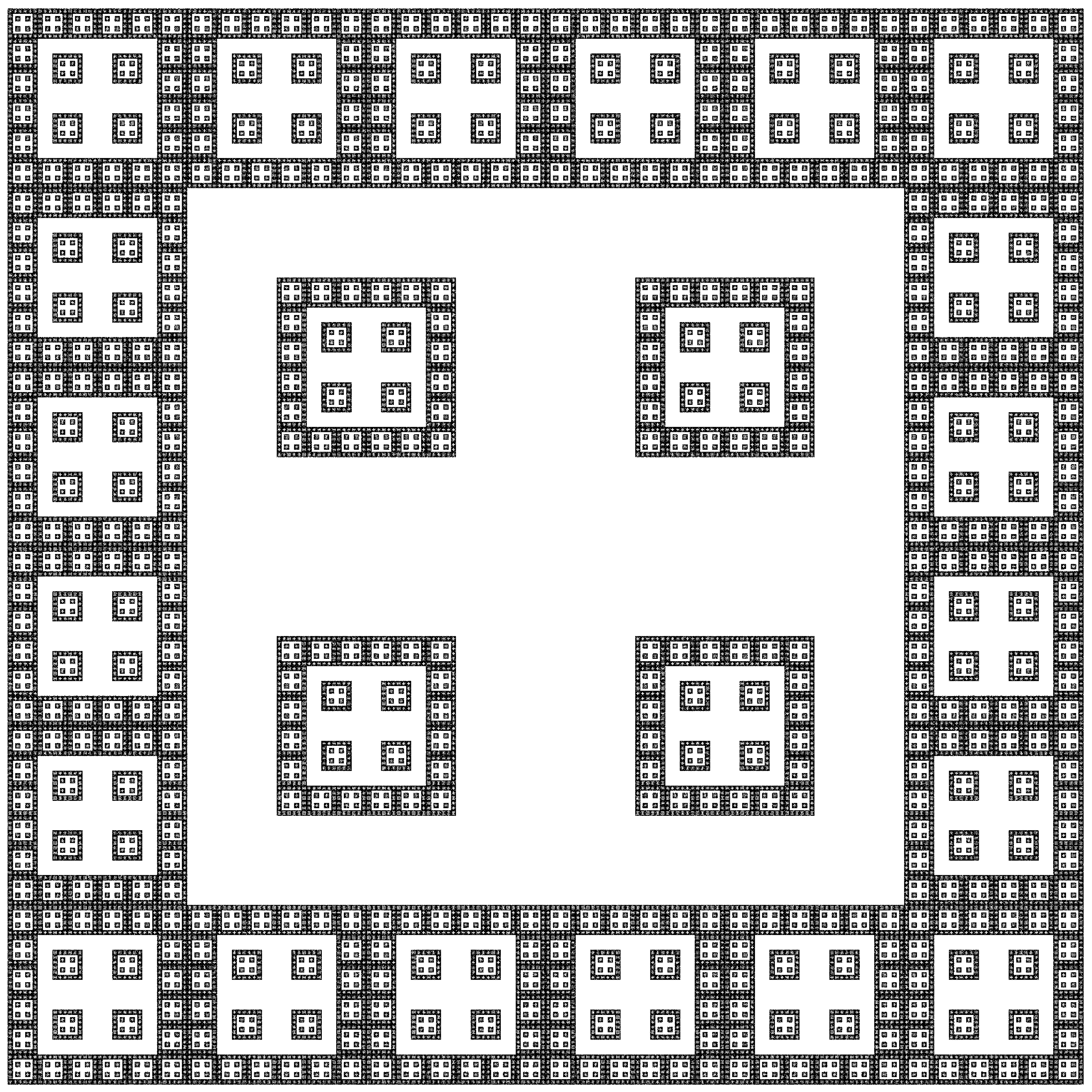} 
& 1.774	& 1.751	& 1.790	& 1.748	& 1.667	& 1.809	& 1.708	& 1.752	\\ \hline
\includegraphics[height=18pt]{quadrate.png} 
& 1.794 & 1.781	& 1.802	& 1.772	& 1.710	& 1.897	& 1.741	& 1.815 \\ \hline
\includegraphics[height=18pt]{carpet.png}
& 1.893	& 1.866	& 1.872	& 1.857	& 1.723	& 1.923	& 1.790	& 1.893	\\ \hline
\includegraphics[height=18pt]{modcarpet.png}
& 1.893	& 1.826	& 1.879	& 1.857	& 1.715	& 1.914 & 1.788	& 1.912	\\ \hline
\includegraphics[height=18pt]{tricarpet.pdf}
& 1.893	& 1.778	& 1.785	& 1.839	& 1.605	& 1.851 & 1.766	& 1.819	\\ \hline
\includegraphics[height=18pt]{fortyAndSevenHigh.png}
& 1.896	& 1.819	& 1.859	& 1.854	& 1.690	& 1.786 & 1.854	& 1.854	\\ \hline
\includegraphics[height=18pt]{fortyAndSevenLow.png}
& 1.896	& 1.839	& 1.875	& 1.863	& 1.716	& 1.766	& 1.789	& 1.945	\\ \hline

\end{tabular}
\caption[Comparison of Dimension Estimates: images of size $3000\times3000$]{Comparison of dimension estimates for large-sized images.}
\end{table}

\begin{table}[!h]
\label{tab:small}
\centering
\begin{tabular}{c||*{7}{p{1.1cm}|}c}
	\multicolumn{8}{c}{images of size $250\times250$} \\ \hline
& theor.& box 	& local	& area	& b'dary& euler	& joint2& joint3 \\ \hline\hline 
\includegraphics[height=18pt]{CD.png} 
& 1.262	& 1.312	& 1.366	& 1.130	&(1.268)&(1.489)&(1.187)& (1.275)	\\ \hline
\includegraphics[height=18pt]{kochcurve_klein.pdf}
& 1.262	& 1.220	& 1.339	& 1.228	&(1.170)& N/A	& 1.204	& N/A   \\ \hline
\includegraphics[height=18pt]{tree.png} 
& 1.585	& 1.571	& 1.673	& 1.570	& 1.545	& N/A	& 1.539	& N/A	\\ \hline
\includegraphics[height=18pt]{gasket.png}
& 1.585	& 1.548	& 1.660 & 1.559	& 1.519	& 1.551	& 1.521	& 1.499 \\ \hline
\includegraphics[height=18pt]{gasket_schief.png}
& 1.585	& 1.516	& 1.657	& 1.553	& 1.513	& 1.512	& 1.515	& 1.505 \\ \hline
\includegraphics[height=18pt]{triangle2.png} 
& 1.588 & 1.535	& 1.627	& 1.545	& 1.495	& 1.517	& 1.498	& 1.537	\\ \hline
\includegraphics[height=18pt]{M1.pdf}
& 1.774	& 1.709	& 1.792	& 1.721	& 1.650	& 1.844	& 1.648	& 1.749	\\ \hline
\includegraphics[height=18pt]{M2.pdf} 
& 1.774	& 1.738	& 1.791 & 1.733	& 1.601	& 1.555	& 1.670	& 1.556	\\ \hline
\includegraphics[height=18pt]{M3.png} 
& 1.774	& 1.799	& 1.798	& 1.727	& 1.616	& 1.791	& 1.632	& 1.722	\\ \hline
\includegraphics[height=18pt]{quadrate.png} 
& 1.794 & 1.726	& 1.803	& 1.748	& 1.704	& 1.652	& 1.680	& 1.745 \\ \hline
\includegraphics[height=18pt]{carpet.png}
& 1.893	& 1.797	& 1.874	& 1.821	& 1.735	& 1.665	& 1.713	& 1.894	\\ \hline
\includegraphics[height=18pt]{modcarpet.png}
& 1.893	& 1.747	& 1.878	& 1.813	& 1.610	& 1.678 & 1.662	& 1.899	\\ \hline
\includegraphics[height=18pt]{tricarpet.pdf}
& 1.893	& 1.823	& 1.880	& 1.824	& 1.607	& 1.780	& 1.675	& 1.814	\\ \hline
\includegraphics[height=18pt]{fortyAndSevenHigh.png}
& 1.896	& 1.805	& 1.891	& 1.819	& 1.739	& 1.279	& 1.694	& 1.876	\\ \hline
\includegraphics[height=18pt]{fortyAndSevenLow.png}
& 1.896	& 1.834	& 1.855	& 1.825	& 1.775	& 1.710	& 1.712	& 1.895	\\ \hline
\end{tabular}
\caption[Comparison of Dimension Estimates: images of size $250\times250$]{Comparison of dimension estimates for small-sized images.} 
\end{table}

It seems that data given by the box-counts $N_\delta$ and by the rescaled areas $y_{2j}$ were the most reliable, as they allowed a good fit and yielded an accurate slope. The rescaled boundary lengths $y_{1j}$ allowed almost as good a fit, however the slope seemed to be systematically to low by around $0.05$, especially for the square-like sets, which tended to have a slope too low by almost $0.15$. The rescaled euler numbers $y_{0j}$ were the least reliable data set. The goodness of its fit strongly depended on the non-arithmeticity of the sets, the best fit being achieved with $\triangle$, the self-similar triangle. As more ``realistic'' fractals will almost certainly be non-arithmetic, there is the hope that this data set will perform better on non-artificial sample images.

\section{Estimates of the Fractal Curvatures.}

Here we describe the first attempt we know of to measure fractal curvatures with an image analyser. Since we expect this task to be very difficult, we used only the highest resolution versions of the sample images ($3000\times3000$ pixels).


Like most image analysers, the algorithms from this thesis all use the pixel scale, i.e. the length $1$ corresponds to the distance between two (horizontally or vertically) neighbouring pixels, meaning that on this scale the base of e.g. the Sierpi\'{n}ski Carpet has length $2980$ and not $1$ as assumed in Winter's calculations. 

Since the sample sets all satisfy $s_k = \underline s_k = s-k$, by proposition~\ref{prop:rescaling} we can switch between the two scales by multiplying by $\lambda^s$, where $\lambda$ is the scale ratio and $s$ is the fractal dimension:
\[C^f_k(\lambda F) = \lambda^{s_k+k} C^f_k(F) = \lambda^s C^f_k(F).\]
Table \ref{tab:curvature-estimates} shows a comparison between the fractal curvature measures which have been calculated this way and the estimated fractal curvatures, obtained by the algorithm from chapter~\ref{ch:main}. When comparing the fractal curvatures, one should keep in mind that the ``unit'' is $pixel^s$, and thus curvatures are only comparable if the corresponding sets have the same dimension. Note also that one should make sure that each set has equal base length, as is the case here (all bases have length 2920 pixels). 

\begin{table}[!h]
\centering
\begin{tabular}{r|l|*{3}{|c}}
\label{tab:curvature-estimates}
&dim 
&$\overline C_0^f$
&$\overline C_1^f$
&$\overline C_2^f$\\ 
\hline 
\hline
\multirow{2}{18pt}{\includegraphics[height=18pt]{tree.png}} 	
 &\multirow{2}{1.1cm}{1.585}	
	&	?	&	?	&	?\\ 
 & 	&	0.0	&$145556$	&$631342$	\\ 
\hline
\multirow{2}{18pt}{\includegraphics[height=18pt]{gasket.png}} 	
 &\multirow{2}{1.1cm}{1.585}	
 	&$-13197$	&$117230$	&$564100$	\\ 
 & 	&$-11109$	&$124064$	&$568985$		\\ 
\hline
\multirow{2}{18pt}{\includegraphics[height=18pt]{gasket_schief.png}} 	
 &\multirow{2}{1.1cm}{1.585}	
  &?			&	?	&	?	\\ 
&  &$-9214 $		&$106250 $	&$498975 $	\\ 
\hline
\multirow{2}{18pt}{\includegraphics[height=18pt]{triangle2.png}} 	
 &\multirow{2}{1.1cm}{1.588}	
 &$-9843$		&$100416$	&$487649$	\\
&&$ -7733$			&$94940 $		&$ 448388 $ 	\\
\hline
\multirow{2}{18pt}{\includegraphics[height=18pt]{carpet.png}} 	
 &\multirow{2}{1.1cm}{1.893}	
 &$-58716$		&$262770$	&$4900200$	\\ 
&  &$-59997$		&$386937$	&$4861736$	\\ 
\hline
\multirow{2}{18pt}{\includegraphics[height=18pt]{modcarpet.png}} 	
 &\multirow{2}{1.1cm}{1.893}	
 &$-50742$		&$260960$	&$4871275$	\\
&  &$-48113$		&$339448$	&$4627544$	\\ 
\hline
\multirow{2}{18pt}{\includegraphics[height=18pt]{fortyAndSevenHigh.png}} 	
 &\multirow{2}{1.1cm}{1.896}	
 &$-54228$		&$242166$	&$4644054$	\\
&  &$-42878$		&$292915$	&$4306752$	\\ 
\hline
\multirow{2}{18pt}{\includegraphics[height=18pt]{fortyAndSevenLow.png}} 	
 &\multirow{2}{1.1cm}{1.896}	
 &$-60811$		&$271565$	&$5207828$	\\
&  &$-65169$		&$371606$	&$5843391$	\\
\hline
\multirow{2}{18pt}{\includegraphics[height=18pt]{M1.pdf}} 	
 &\multirow{2}{1.1cm}{1.774}	
&?&?&?\\
&  &$-30530 $		&$280567 $	&$2056714$ \\
\hline
\multirow{2}{18pt}{\includegraphics[height=18pt]{M2.pdf}} 	
 &\multirow{2}{1.1cm}{1.774}	
&?&?&?\\
&  	&$-6699$	&$278438 $	&$2093975$	\\
\hline
\multirow{2}{18pt}{\includegraphics[height=18pt]{M3.png}} 	
 &\multirow{2}{1.1cm}{1.774}	
&	?	&	?	&	?	\\
&  		& $-9380$ 	& $351984$	& $2754060$	\\
\end{tabular}
\caption[Estimated Curvatures]{Estimated Fractal Curvatures. In each row, the upper number is the theoretically calculated curvature, whereas the lower number is the estimate returned by the algorithm from chapter~\ref{ch:main}. A question mark denotes curvatures that have not been calculated yet, and an ``N/A'' mark means that for the corresponding fractal either $\overline C_0^f(F) = 0$ or the asymptotic of $y_{0j}$ is not good enough.}
\end{table}

A typical relative error lies at around $10\%$, however, for the rectangular-like sets, 
$\overline C_1^f$ seems to be systematically overestimated. This corresponds with the observation from the previous section, where the measured growth of the rescaled boundary lengths $y_{1j}$ was too low, also especially for the rectangular-like sets: Recall the calculation
\[\Gamma_k^{(m)} := \frac{1}{t_m - t_0} \sum_{j=1}^m \exp (y_{kj} - \hat s^{(m)}x_j) (t_{j}-t_{j-1})\]
of the $k$-th fractal curvature estimate on page~\pageref{eq:gamma}. If $\hat s$ is assumed higher than a least-squares fit of $y_{1j}$ vs. $x_j$ would suggest, the convexity of the exponential function will yield a relatively too high value, since the $x_j$ are all negative. The estimator $\exp (\hat D_1)$ shows similar behaviour.

The last four rows in table \ref{tab:curvature-estimates} show the two pairs ``Window and Gate'' and ``M1 and M2''. In both pairs, the two sets have the same dimension, but the first one \textit{seems more lacunar}. This is reflected in a lower Minkowski-content, as has been already calculated for the first pair. The author's confident guess is that it is also the case for the second pair. In any case, lacunarity seems to have the effect of reducing \textit{all} of the curvature measure estimates. As this behaviour is parallel to the dimension estimates, the absolute values of the curvatures might not be suitable to discern different fractals of equal dimension. However, their size relative to each other yields an important geometric invariant (see nxt section).

\section{Characterization of $\Sd$-sets beyond dimension}
\label{sec:beyond}

\paragraph{Specific fractal curvatures.}

For a systematic categorization of fractal sets, \textit{geometrically invariant characteristics} are useful. Here, by a characteristic we simply mean a functional
\[F:\mathbb K^d \rightarrow \mathbb R\] 
on the class $\mathbb K^d$ of compact sets in $\Rd$. $F$ is called 
\begin{itemize}
  \item \textit{motion invariant} if $F(g(K)) = F(K)$ for all $K \in \mathbb K^d$ and for all euclidean motions $g$ on $\Rd$, 
  \item \textit{homogeneous of degree k} if $F(\lambda K) = \lambda^k F(K)$ for all $\lambda \in (0,\infty)$, and 
  \item \textit{scaling invariant} if $F$ is homogeneous of degree 0. 
\end{itemize}
For example, all definitions of fractal dimension satisfy motion and scaling invariance. We now define the following characteristics:

\begin{definition}
Let $A \in \mathbb K^d$ be such that $\overline C_k^f(A)$ exists for all $k\in\{0,\ldots,d\}$, and assume $\overline C_d^f(A) > 0$. Then for each $k \in \{0,\ldots, d-1\}$ 
\[\Xi_k(A) :=  \dfrac{\overline C_k^f(A)}{\overline C_d^f(A)} \]
is called the $k$-th specific fractal curvature of $A$.
\end{definition}

Note that $\overline C_d^f(F)$ is always positive for self-similar sets satisfying the OSC as shown by Gatzouras in \cite{gatzouras}, so $\Xi_k$ is well-defined for all $\Sd$-sets. 

Since the fractal curvatures are motion-invariant and homogeneous for all $k\in \{0,\ldots,d-1\}$, so are the specific fractal curvatures. There are the following two important cases: 
\begin{proposition}
  Assume that $\Xi_k(A)$ is well-defined, and let $\lambda >0$.
\begin{enumerate}
  \item 	If $A \in \mathcal R^d$, then \[\Xi_k(\lambda A) = \lambda^{k-d}\Xi_k(A).\]
  \item		If $s_k(A) = s-k$ for all $k\in\{0,\ldots,d\}$,
		then \[\Xi_k(\lambda A) = \Xi_k(A).\]
\end{enumerate}
\end{proposition}

\textit{Proof:} This is a simple consequence of the $s_k + k$- homogeneity of $C_k^f$, the $s_k$ being 0 for sets in the convex ring. $\square$

\paragraph{Local observations.}

Let ``$\sim$'' be the usual geometrical similarity relation, restricted to $\Sd$. If for two sets $F, G \in \Sd$ their binary representations $\hat F$ and $\hat G$ are given and one has to decide whether or not $F \sim G$, then $\hat F$ and $\hat G$ can be rescaled to have the same diameter, and their (average) Minkowski-contents can be estimated. If they differ significantly, this could be evidence against $F \sim G$. 
\begin{figure}
\centering
  \includegraphics[width=4cm]{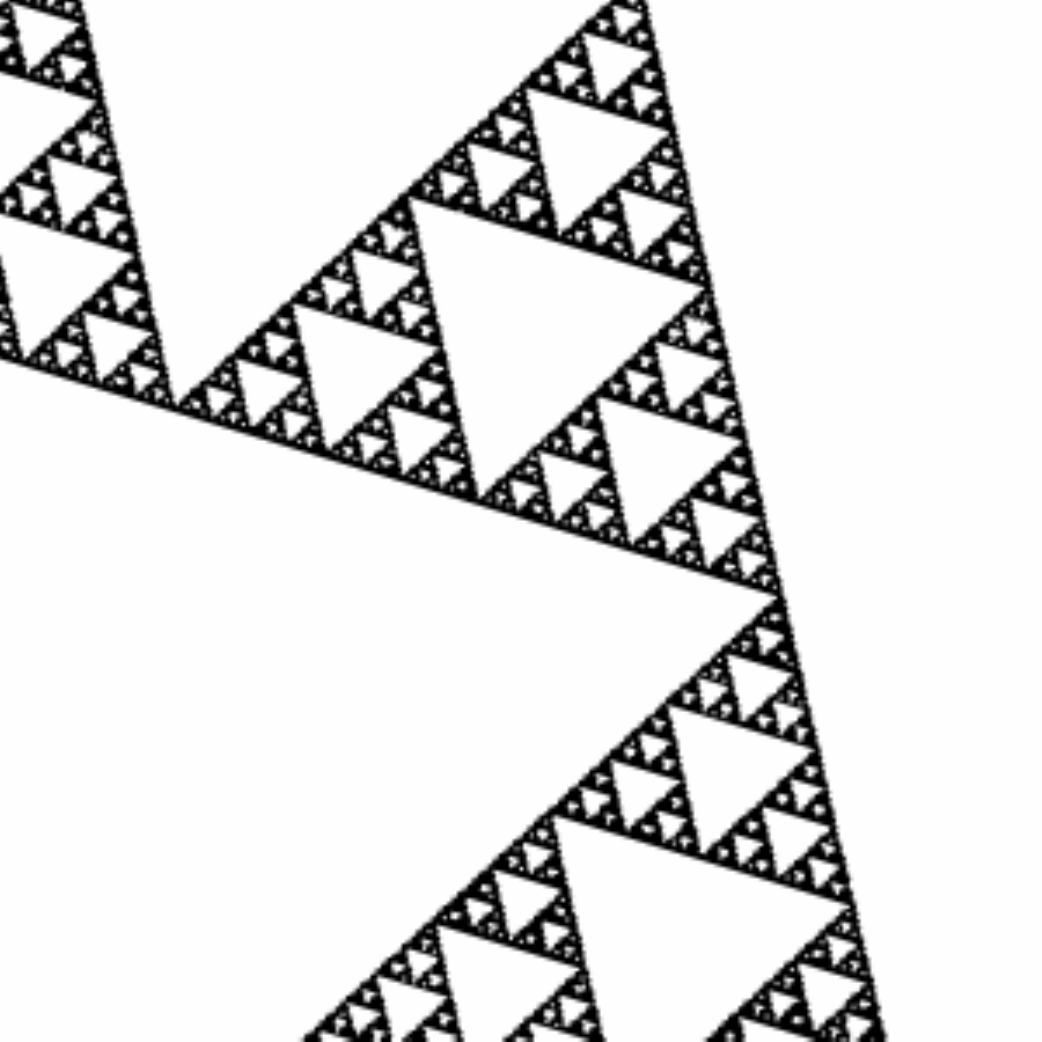}
  \includegraphics[width=4cm]{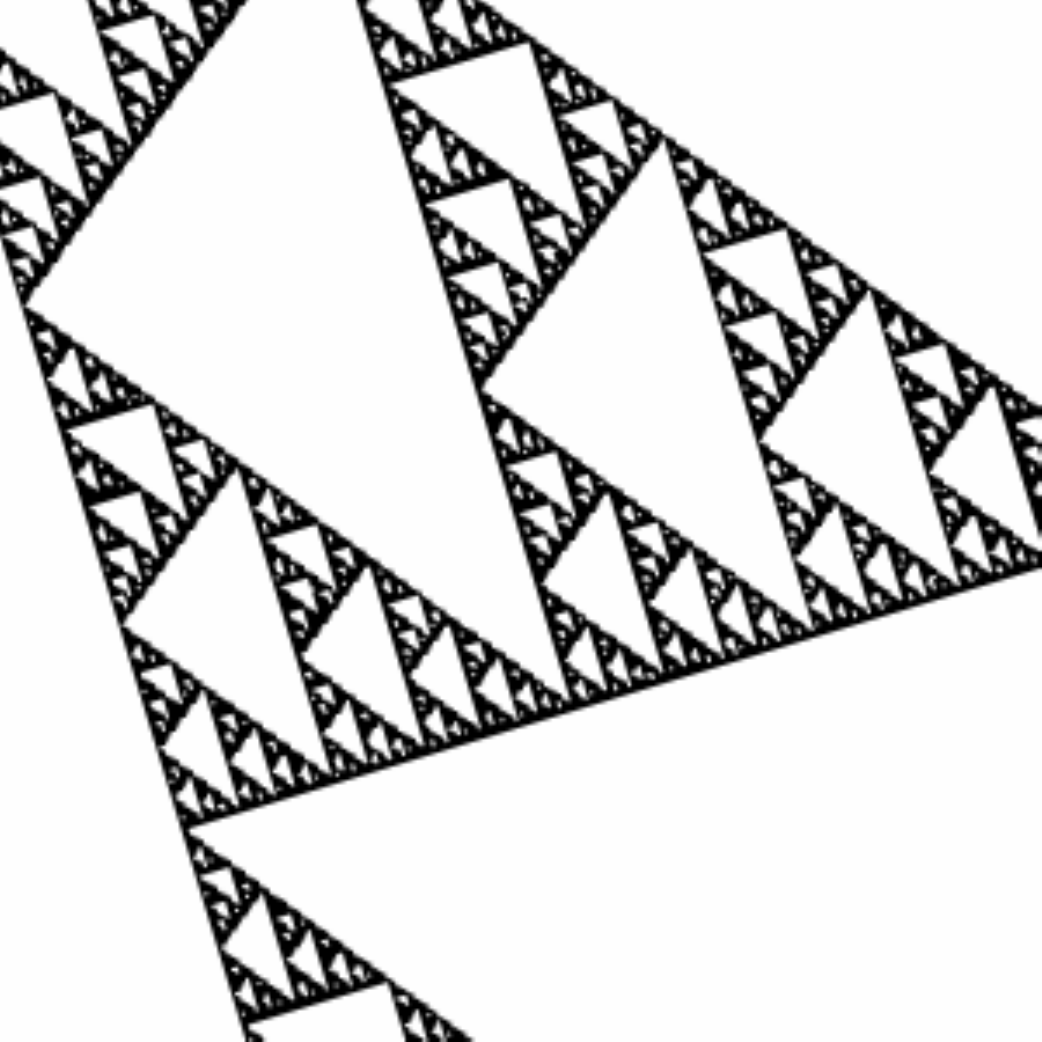}
  \caption{Local observations of the Sierpinski gasket (left) and $\triangle$ (right).}
  \label{fig:cut-out}
\end{figure}

\begin{table}[!ht]
\centering
\begin{tabular}{r|r|r||r|r}
\label{tab:specific-curvatures}
& $\widehat{\dim _B}$ & $\dim $ &$\Xi_0$ & $\Xi_1$\\ 
\hline 
\hline
\multirow{2}{18pt}{\includegraphics[height=20pt]{CD.png}} 	
& \multirow{2}{1.1cm}{1.254} &\multirow{2}{1.1cm}{1.262}   
&N/A	    &N/A\\ &&&		& $(0.412)$\\ 
\hline
\multirow{2}{18pt}{\includegraphics[width=20pt]{kochcurve_klein.pdf}} 
& \multirow{2}{1.1cm}{1.270} &\multirow{2}{1.1cm}{1.262}	
	&N/A	&N/A\\
&&  	&$ $	&$(0.373) $\\ 
\hline
\multirow{2}{18pt}{\includegraphics[height=20pt]{gasket.png}} 
& \multirow{2}{1.1cm}{1.539} &\multirow{2}{1.1cm}{1.585}	
	&$-0.0234$	&$0.208$\\ 
&&  	&$-0.0195$	&$0.218$\\ 
\hline
\multirow{2}{18pt}{\includegraphics[height=20pt]{tree.png}} 
& \multirow{2}{1.1cm}{1.564} &\multirow{2}{1.1cm}{1.585}	
	&0	&?\\
&&  	&0	&$0.230$\\ 
\hline
\multirow{2}{18pt}{\includegraphics[height=20pt]{triangle2.png}} 	
& \multirow{2}{1.1cm}{1.573} &\multirow{2}{1.1cm}{1.588}	
	&$-0.0202$	&$0.206$\\
&& 	&$-0.0172$	&$0.211$\\
\hline
\multirow{2}{18pt}{\includegraphics[height=20pt]{gasket_schief.png}} 	
& \multirow{2}{1.1cm}{1.588} &\multirow{2}{1.1cm}{1.585}	
	&	?	&	?\\ 
&&  	&$-0.0185$	&$0.213$\\ 
\hline
\multirow{2}{18pt}{\includegraphics[height=20pt]{M1.pdf}} 	
& \multirow{2}{1.1cm}{1.706} &\multirow{2}{1.1cm}{1.774}	
	&?	&?\\
&  &	&$-0.0148 $	&$0.136 $\\ 
\hline
\multirow{2}{18pt}{\includegraphics[height=20pt]{M2.pdf}}
& \multirow{2}{1.1cm}{1.730} &\multirow{2}{1.1cm}{1.774}	
	& ?	&?\\
&  &	& $-0.00319$	&$0.1329 $\\ 
\hline
\multirow{2}{18pt}{\includegraphics[height=20pt]{M3.png}}
& \multirow{2}{1.1cm}{1.730} &\multirow{2}{1.1cm}{1.774}	
	& ?	& ?\\
&  &	& $-0.00340$	&$0.1278  $\\ 
\hline
\multirow{2}{18pt}{\includegraphics[height=20pt]{tricarpet.pdf}} 	
& \multirow{2}{1.1cm}{1.778} &\multirow{2}{1.1cm}{1.893}	
	&?		&?		\\
&  &	&$-0.00715$	&$0.0675$	\\ 
\hline
\multirow{2}{18pt}{\includegraphics[height=20pt]{quadrate.png}} 
& \multirow{2}{1.1cm}{1.781} &\multirow{2}{1.1cm}{1.794}	
	&?	&?			\\
&  &	&$-0.00577$	&$0.1114$	\\ 
\hline
\multirow{2}{18pt}{\includegraphics[height=20pt]{fortyAndSevenHigh.png}} 	
& \multirow{2}{1.1cm}{1.819} &\multirow{2}{1.1cm}{1.896}	
	&$-0.0117$	&$0.0521$\\
&  &	&$-0.0100$	&$0.0680$\\ 
\hline
\multirow{2}{18pt}{\includegraphics[height=20pt]{modcarpet.png}} 
& \multirow{2}{1.1cm}{1.826} &\multirow{2}{1.1cm}{1.893}	
	&$-0.01041$	&$0.0536$\\
&  &	&$-0.01039$	&$0.0733$\\ 
\hline
\multirow{2}{18pt}{\includegraphics[height=20pt]{fortyAndSevenLow.png}} 	
& \multirow{2}{1.1cm}{1.839} &\multirow{2}{1.1cm}{1.896}	
	&$-0.0117$	&$0.0521$\\
&  &	&$-0.0112$	&$0.0636$\\
\hline
\multirow{2}{18pt}{\includegraphics[height=20pt]{carpet.png}} 
& \multirow{2}{1.1cm}{1.866} &\multirow{2}{1.1cm}{1.893}	
	&$-0.01183$	&$0.0536$\\ 
&  &	&$-0.01234$	&$0.07958$
\end{tabular}
\caption[Estimated specific fractal curvatures]{Estimated specific fractal curvatures: The first two columns contain the estimated dimension via box-counting and the theoretical value, respectively. In the right two columns, the upper value of each row is theoretical and the lower value is the returned estimate. The Cantor Dust and the Koch Curve are marked ``N/A'' since they are not in $\Sd$.}
\end{table}

However, in many applications $\hat F$ and $\hat G$ will only be local representations, i.e. binary images of $F \cap W_1$ and $G \cap W_2$, respectively, where the $W_i$ are (e.g. rectangular) bounded observation windows. 
If $\hat F$ is to be compared to $\hat G$, one faces the problem that along with the position of the window also the amount of white space will vary, and thus lacunarity analysis by the box-gliding algorithm (chapter~\ref{ch:discreteness}) and any total curvature $\overline C_k^f(\hat F)$ will not be a reliable source of information. (Note that it might happen that $F \cap W_1$ does not have polyconvex parallel sets anymore even though $F$ does, and that thus $\overline C_k^f(F \cap W_1)$ might not be defined any more; but this problem gets lost in the discretization procedure anyways.)
However, the (average) Minkowski-content $\overline C_2^f(F)$ can serve as a normalization factor for the other fractal curvatures, i.e. the specific fractal curvatures can still be calculated. If they are significantly different then $F$ and $G$ might still be distinguished from each other even though the dimension appears to be the same.

\section{Discussion}

\paragraph{Dimension estimates.}

The tests on the sample images suggest that the overall accuracy of dimension estimates via the simultaneous regression on all data sets ``euler number'', ``boundary length'' and ``area'' is comparable to the accuracy of the sausage method and the box-counting method. More data do not lead to more accuracy here:
\begin{itemize}
\item The measured growth of the boundary length of the parallel sets as $\eps \downarrow 0$ is slightly lower than the theory would suggest (see column ``b'dary'' of table~\ref{tab:big}). Measurements based on the algorithms of references \cite{minkdiscr} and \cite{ohser} both qualitatively yielded almost the same results; the latter algorithm showed slightly slower growth. A partial reason for this might be the convexity of the plot  $y_1 = \log\frac{(C_1(F_\eps)}{\eps}$ against $x = -\log\eps$, which does not die off quickly enough as $x$ approaches infinity. Cutting off the data \textit{before} the break condition $C_0^{var}(F_\eps) \leq 2$ applies might help getting rid of some of the negative bias that the slope of the plot has, but the overall stability of the estimate will become worse. 
  \item For the measured\footnote{Recall that $N(F_\eps)$ is the number of connected components of $F_\eps$.} $2N(F_\eps)-C_0(F_\eps)$ to be close to the total variational measure $C_0^{var}(F_\eps)$, the holes of $F_\eps$ need to be convex (also see appendix). But even if they are convex, they might exhibit stairlike behaviour due to arithmeticity, which does not allow for a good regression fit. Of the considered sets, only the self-similar triangle allowed for a good fit of all three data sets.
\end{itemize}

\paragraph{Curvature estimates.}

The accuracy of fractal curvature estimates strongly depends on the accuracy of the dimension estimate. As noted before, the data set of $C_0^{var}(F_\eps)$ can seldomly be measured properly, and thus it should be included in the regression analysis only in special cases. Note however that even if $C_0^{var}(F_\eps)$ is excluded from the regression estimate of the dimension, then it is still possible to calculate the 0-th fractal curvature estimate $\Gamma_0(F)$ (compare equation (\ref{eq:gamma})).

The functional dependence of the curvature estimates on the dimension estimate can be read off the following alternative formulation of equation (\ref{eq:gamma}):
\[\Gamma_k^{(m)} = \frac{1}{t_m - t_0} \sum_{j=1}^m  {\rm sgn}(C_k(F_{e^{-x_j}})) \exp \left(y_{kj}-sx_j\right) (t_{j}-t_{j-1})\]
If the slope $s$ is higher or lower than the plot of $y_{kj}$ against $x_j$ suggests, there will be higher and lower arguments to the exponential function, and since it is convex, positive deviations will come out stronger than negative deviations. Thus a wrong dimension estimate will always result in too high values for $\Gamma_k^{(m)}$, no matter if the dimension estimate itself will be too high or too low. As this behaviour is the same for all $k \in \{0,1,2\}$, one might hope that the negative influence of a bad dimension estimate will cancel out for the specific curvatures $\Xi_0(F)$ and $\Xi_1(F)$, and that thus $\Xi_0(F)$ and $\Xi_1(F)$ are less susceptible to wrong dimension estimates.

\paragraph{Specific curvatures.}

As pointed out before, it is hard to compare two sets by their (average) Minkowski-content, if one cannot be sure that the same scale should be used on both sets; in this case, the specific curvatures seem more appropriate, though the Minkowski-contents should be kept in mind. 

Table~\ref{tab:specific-curvatures} shows the specific curvatures of the sample sets in the right columns, where in each row the upper value is theoretical and the lower value is estimated from the binary image. The general observation is that the specific curvatures increase as the dimension of a self-similar set decreases from $2$ to $1$. Also note that for the observed sets, if two fractals have indistinguishable dimensions then their specific curvatures may be different; however, in most cases like this the specific curvatures will be too close together for an image analyser to depict the difference with high enough confidence.
There are cases, however, where the specific curvature estimates can distinguish sets where the dimension estimates cannot: As an example, compare the Tripet to the self-similar square $\square$ in table \ref{tab:specific-curvatures}. The box-counting estimates are almost the same, but $\Xi_1$ is rather different. 

\chapter*{Concluding Remarks}
\addcontentsline{toc}{chapter}{Concluding Remarks}
\pagestyle{myheadings}
\markboth{CONCLUDING REMARKS}{CONCLUDING REMARKS}

In most of the literature on estimates of fractal dimension, the slope of a logarithmic plot of $C_d(F_r)$ against $r$ earned almost all attention, whereas the intersection point of the fitted line with the vertical axis was mostly ignored. Furthermore, in the context of fractal analysis the author has found almost no occurence of any other generalized volumes than $V_d$ ($ = C_d$) in the literature; The only exception was B. Mandelbrot's exposition \cite{mandelbrot:mfl} on ``gap-lacunarity'' of Fractals on the one-dimensional unit interval, which is essentially the $0$-th fractal curvature $C_0$. This thesis shows that, at least for the special case of $\Sd$-sets, a large amount of information is ignored if only the $d$-dimensional volume of the parallel sets is analysed. 

\begin{figure}[h]
  \centering
  \subfigure[Only two affine maps on $\mathbb R^2$ generate this attractor.] {
  	\includegraphics[width=10cm]{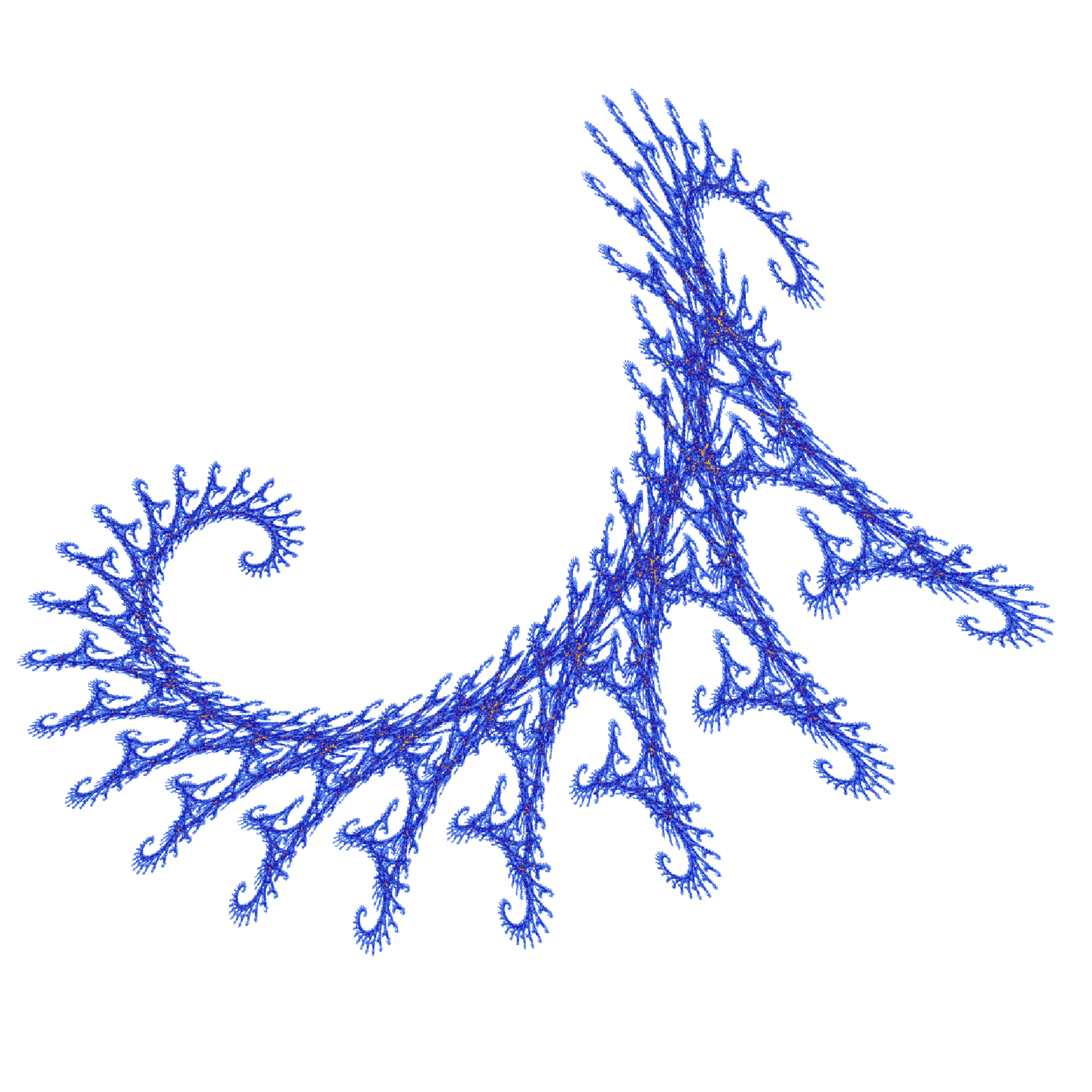}
  }
	\hspace{1cm}
  \subfigure[The usual plot of $y_{kj}$ against $x_j$, $k \in\{0,1,2\}$. The fit seems to be just as good as for an $\Sd$-set.] {
	\includegraphics[width=\textwidth]{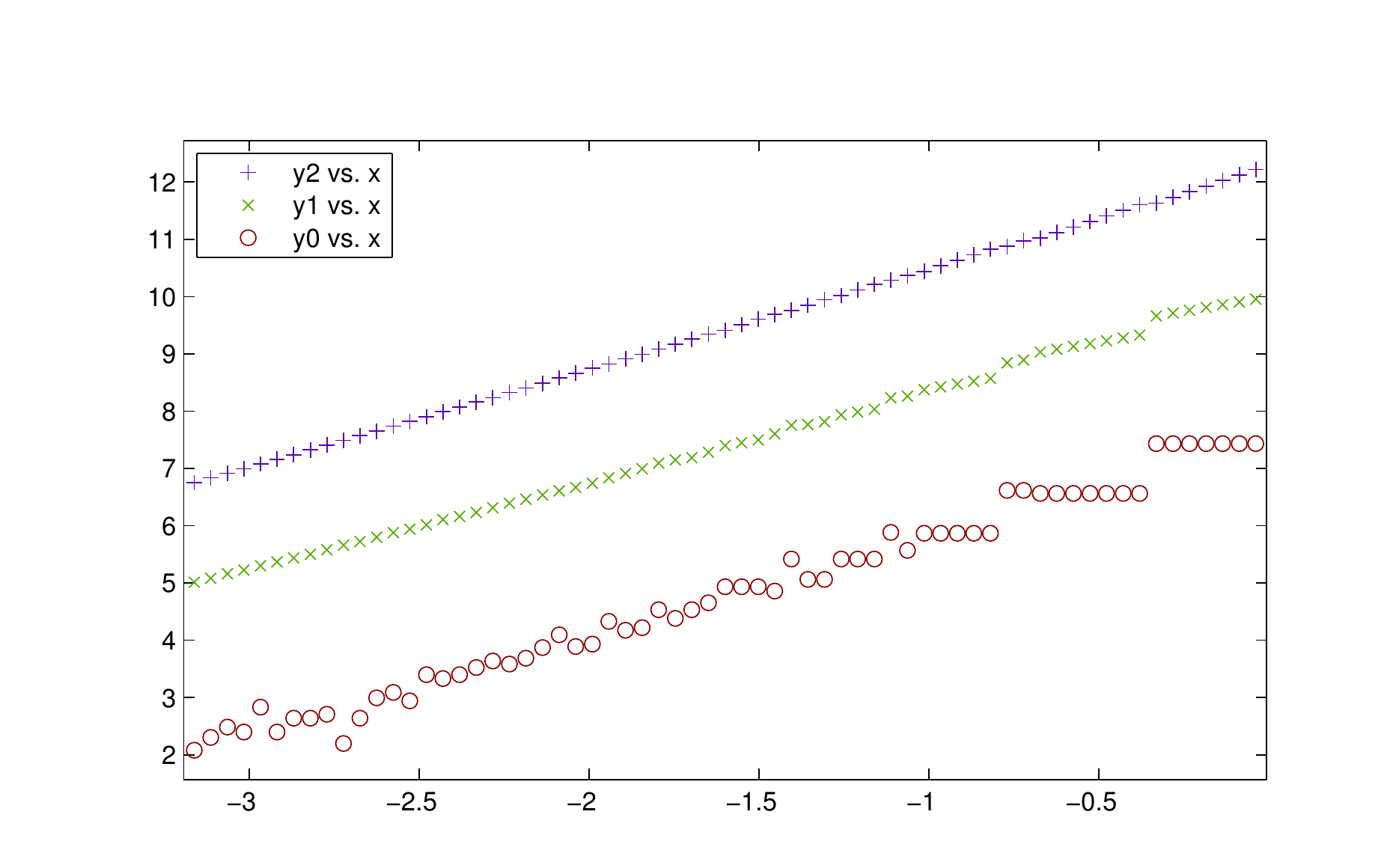} 
  }
  \caption{The Dragon, a self-affine set}
  \label{fig:dragon}
\end{figure}

The theory of fractal curvatures measures is still at its beginning stage, and only for self-similar sets there are numerical expressions for their calculations. Their definition can be somewhat extended, namely to sets whose parallel sets are finite unions of sets of positive reach; However in that case little to nothing is known about the scaling exponents $s_k$, or even if the growth of $C_k(F_r)$ can be compared to $r^s$ for some $s$. But there does seem to be a similar behaviour of self-\textit{affine} sets in terms of the scaling behaviour of the curvatures of parallel sets, as Figure \ref{fig:dragon} shows. 

An example for the case where $C_k(F_r) \sim r^s$ does \textit{not} hold true for any $s \in \mathbb R$ is a brownian path $B_t([0,1])$ on the time interval $[0,1]$ in $\mathbb R^2$: For this statistically self-similar set, 
\begin{equation}
\label{eq:sausage-asymptotics}
  \mathbb E C_2(B_t([0,1])_r) \sim \frac{\pi}{\vert \log r \vert}, ~~ 
	\mathbb E C_1(B_t([0,1])_r) \sim \frac{\pi}{2r\log^2 r} ~~ (r\downarrow 0),
\end{equation}
the second asymptotic having been proved only recently in \cite{rataj2005esa}. Numerical simulations of $\mathbb E C_0(B_t([0,1])_r$ can be found in \cite{meschenmoser}, as well as an empirical confirmation of equations (\ref{eq:sausage-asymptotics}). In this case, regression analysis on the log-log scale does not seem apropriate any more, and other means have to be found to measure Fractal Curvature.


\bibliographystyle{alpha}  
\bibliography{petersbibliography}

\chapter*{Appendix}
\addcontentsline{toc}{chapter}{Appendix}

\pagestyle{myheadings}
\markboth{APPENDIX}{APPENDIX}

\section*{Estimating $C_0^{var}(F_\eps)$}
\label{app:C0var}

In general, with an image analyser only the Euler-number $C_0(F_\eps)$ of the parallel sets can be determined; for the regression analysis of chapter~\ref{ch:main}, however, $C_0^{var}(F_\eps)$ is needed in order to find the scaling exponent $s_0$. Since in 2D, the Euler-number equals the number $N$ of connected components minus the number $Q$ of holes, one could think that the total variational measure $C_0^{var}(F_\eps)$ would equal $N+Q$, but this is false in general. The following lemma gives a sufficient condition for this to be true.

\begin{lemma*}
\label{lem:C0var}
  Let $F\in\mathcal R^2$, and let $N$ be the number of connected components of $F$ and $Q$ the number of holes of $F$ (i.e. the number of connected components of $F^c$ minus 1). Moreover, assume that the closures of the components of $F^c$ are disjoint. Then 
\begin{enumerate}
  \item $\vert C_0(F) \vert \leq N + Q \leq C_0^{var}(F).$
  \item Let $\bd F = \bigcup_{i=1}^M \Gamma_i$ be the decomposition of the boundary of $F$ into its connected components $\Gamma_i$. If 
\begin{equation}
\label{app:equality}
C_0^{var}(\Gamma_i) = 1 \textnormal{ for all } i,
\end{equation}
then 
\[C_0^{var}(F) = N+Q.\]
\end{enumerate}
\end{lemma*}

\textit{Proof:} Since $C_0(F) = N - Q$ the left inequality in 1 is clear. Now we use the result $M = N + Q$, where $M$ is the number of connected components of the boundary of $F$: 

To prove this, start off by labelling the components of $F$ and of $F^c$ in the following way: Assign 0 to the only unbounded component. Assign 1 to all remaining components of $F$ that have a boundary with 0, and proceed with assigning $i$ to the remaining components that have a boundary with the components labelled $i-1$ until all components of both $F$ and $F^c$ have been labelled. By assumption, no component of $F^c$ has a boundary with any other component of $F^c$, and thus components of $F^c$ have even labels, whereas components of $F$ will have odd labels. There will be a maximum label $m$, since $F \in \mathcal R^2$. Now for all odd $k$ define $F_k$ to be the union of all components with odd labels up to $k$. Then the number of components of the boundary of $F_k$ will be the equal to the number of components with  label $\leq (k+1)$ minus 1, as can be checked by induction on $k$. In particular, the boundary of $F = F_{m-1}$ will have as many components as the total number of components of both $F$ and $F^c$ together, minus 1. But this is the number of connected components of $F$ plus the number of holes of $F$. 

Since always $C_0^{var}(\Gamma_i) \geq 1$, the right inequality in 1 follows, and 2 is just a consequence of $M = N + Q$. $\square$

Note that the condition $C_0^{var}(\Gamma_i) = 1$ is equivalent to the interior of the curve $\Gamma_i$ being convex. Thus out of all sets considered in section \ref{sec:images}, equality (\ref{app:equality})  applies in the following way to our sample sets: 
\begin{itemize}
  \item It does apply to Sierpi\'nski Gasket and Carpet, Tripet, $\triangle$, Window and Gate. 
  \item It does not apply to the parallel sets of the Cantor dust and Koch curve, since the parallel sets are not polyconvex.  
  \item For the set M1, in fact $C_0^{var}(\Gamma_i) = 3$, so, for small $\eps$, the value $C_0^{var}(F_\eps)$ will be 3 times as high as $N+Q$. This will reflect in a much too low estimate of the 0-th fractal total variational curvature, but it should not affect the estimate of $C_0^f({\rm M1})$.
  \item For the modified Sierpi\'nski Carpet, assuming equality (\ref{app:equality}) is only a little mistake: Here 
\[0 \leq C_0^{var}(F_\eps) - (N+Q) \leq \frac{3}{4}\] 
for all $\eps > 0$, as only the above opening is causing trouble. 
  \item Unfortunately, for the set $\square$ the value $N+Q$ is not close to $C_0^{var}(F_\eps)$ at all: Depending on $\eps$, the holes will seldom be convex and have the total variational curvature 1 as demanded by the above condition. Instead, often they are L-shaped, or cross-shaped, thus having total variational curvature 1.5 or 3 respectively, and other shapes are possible. Thus, for this set \texttt{useEuler} is best set to \texttt{false}. 
  \item The Sierpi\'nski Tree is another example where the 0-th total variational curvature cannot be measured; all parallel sets have one connected component and no holes, whereas the boundary is obviously very curved.

\end{itemize}

\section*{Calculation of the total Fractal Curvature of $\triangle$}

\label{app:triangle}

 We need the curvature scaling functions 

\begin{equation}R_k(\eps) = C_k(\triangle_\eps) - \sum_{i=1}^N
\label{eq:scaling-function}
\textbf{\textup{1}}_{(0,r_i]}(\eps)C_k((S_i\triangle)_\eps)
\end{equation}
for $k\in\{0,1,2\}$ on the interval $(0,1]$. It turns out that the $R_k(\eps)$ are piecewise polynomials of degree $k$ whose transition points are the discontinuities of the indicator functions at $r_1 = \frac{25}{41}$, $r_2 = \frac{20}{41}$ and $r_3 = \frac{16}{41}$, and the radius of the in-circle $\omega$ (see figure~\ref{fig:triangle-dil}). The inner triangle is in fact congruent to the upper (red) smaller copy, and thus we can compute (\cite{wiki:incircle})
\[\omega = 2 r_2 \frac{\frac{1}{2} ac }{a+b+c} = \frac{4}{41}.\]
Now we determine the curvature scaling functions directly via formula~(\ref{eq:scaling-function}). Note that for $0<\eps<r_3$ the inclusion-exclusion formula says that
\[R_k(\eps)=\sum_{\vert I \vert \geq 2} (-1)^{\vert I \vert -1}C_k
	\left(\bigcap_{i \in I}(S_i(\triangle)_\eps \right), \]
the sum being taken over all subsets $I \subseteq \{1,2,3\}$ with at least two elements, and that for $0<\eps<\omega$ the intersection 
\[(S_1(\triangle))_\eps \cap (S_2(\triangle))_\eps \cap (S_3(\triangle))_\eps \]
is empty, so that in this case 
$R_k(\eps) = -\sum_{i\neq j} C_k \left( E_{ij}\right)$ 
where $E_{ij}(\eps) = S_i(\triangle)_\eps\cap S_j(\triangle)_\eps$. 

We find
\begin{equation*}
\begin{array}{rclrclrcl}
 C_0(E_{12}) &=& 1,
&C_0(E_{13}) &=& 1,
&C_0(E_{23}) &=& 1, \\
C_1(E_{12})&=& (1+\frac{\pi+\gamma}{2})\eps, 
&C_1(E_{13})&=& (2+\frac{\pi+\beta}{2})\eps, 
&C_1(E_{23})&=& (3+\frac{\pi+\alpha}{2})\eps, \\
 C_2(E_{12})&=& (1+\frac{\pi}{2}+\gamma) \eps^2,
&C_2(E_{13})&=& (2+\frac{\pi}{2}+\beta) \eps^2,
&C_2(E_{13})&=& (3+\frac{\pi}{2}+\alpha) \eps^2.
\end{array}
\end{equation*}
Now, for each $k\in \{0,1,2\}$ we use the representation
\begin{equation*}
  R_k(\eps)= \left\lbrace
\begin{array}{lrcl}
	C_k(\triangle_\eps), &  r_1 < &\eps&\leq 1 \\
	C_k(\triangle_\eps)- C_k((S_1\triangle)_\eps), & r_2 < &\eps& \leq r_1\\
	C_k(\triangle_\eps)- C_k((S_1\triangle)_\eps) - C_k((S_2\triangle)_\eps), 
		& r_3 < &\eps& \leq r_2\\
	C_k(\triangle_\eps)- \sum_{k=1}^3 C_k((S_k\triangle)_\eps),
		& \omega<&\eps&\leq r_3\\
	-\sum_{i\neq j} C_k(E_{ij}),
		& 0<&\eps&\leq \omega
\end{array}
\right.
\end{equation*}
and putting 
\[A := \frac{ac}{2} = 0.24, ~~~~~ S := a+b+c = 2.4\]
for the area and the circumference of the big triangle we arrive at 
\begin{equation*}
 R_0(\eps)= \begin{cases} 
	1, 	&   r_1 < \eps \leq 1 \\ 
  	0,  	& r_2 < \eps \leq r_1\\
	-1	,& r_3 < \eps \leq r_2\\
	-2	,& \omega<\eps\leq r_3\\
	-3	,& 0 < \eps < \omega
\end{cases}
\end{equation*}
\begin{equation*}
 R_1(\eps)= \left\lbrace
\begin{array}{lrcl}
 \frac{S}{2} + \pi\eps,			& r_1 < \eps \leq 1 ~ \\
(1-r_1)\frac{S}{2},			& r_2 < \eps \leq r_1\\
(1-r_1-r_2)\frac{S}{2} - \pi\eps, 	& r_3 < \eps \leq r_2\\
(1-r_1-r_2-r_3)\frac{S}{2} - 2\pi\eps,  	& \omega<\eps\leq r_3\\
-(6 + 2\pi)\eps,		& 0   < \eps \leq \omega
\end{array}
\right.
\end{equation*}
\begin{equation*}
 R_2(\eps)= \left\lbrace
\begin{array}{lrcl}
 A + S\eps + \pi\eps^2,		& r_1 < \eps \leq 1 ~\\
(1-r_1^2)A + (1-r_1)S\eps,	& r_2 < \eps \leq r_1\\
(1-r_1^2-r_2^2)A+(1-r_1-r_2)S\eps- \pi\eps^2, 	
				& r_3 < \eps \leq r_2\\
(1-r_1^2-r_2^2-r_3^2)A+(1-r_1-r_2-r_3)S-2\pi\eps^2,
				& \omega<\eps\leq r_3\\
-(6 + 2\pi)\eps^2,		& 0   < \eps \leq \omega
\end{array}
\right.
\end{equation*}

Now the fractal curvatures $C^f_k(\triangle)$ can be computed as in theorem~\ref{th:curvature-formula}, but the author used a shortcut-formula from \cite{diss_winter}. 

\newpage

\section*{Signed measures}

\begin{definition*}[signed measure]
\label{def:signed-measure}
  Let $(X,\Sigma)$ be a measurable space. If the mapping
\[\mu: \Sigma \rightarrow \R\]
satisfies 
\[\mu\left(\bigcup_{i=1}^{\infty}A_i\right) = \sum_{i=1}^{\infty} \mu(A_i)\] for a sequence of disjoint sets $A_i \in \Sigma$, then it is called a signed measure. 
\end{definition*}

\begin{definition*}[positive, negative and total variational measure]
  \label{def:pos-neg-var}
Let $\mu$ be a signed measure on the measure space $(X,\Sigma)$. Then the positive variational measure, negative variational measure and the total variational measure are defined respectively by
\begin{eqnarray*}
  \mu^+(A) 	&:=&  \sup\{\mu(B): B \in \Sigma, B\subseteq A  \}  \\
\mu^-(A) 	&:=& -\inf\{\mu(B): B \in \Sigma, B\subseteq A  \} \\
\mu^{var}(A) 	&:=&    \mu^+(A) + \mu^-(A).
\end{eqnarray*}
which are non-negative measures on $(X,\Sigma)$.
\end{definition*}

By the Hahn decomposition theorem, $\mu^+ - \mu^- = \mu$ for every signed measure $\mu$.

\clearpage




\selectlanguage{german}

\thispagestyle{empty} \vspace*{3cm} \hspace{1cm}
\begin{minipage}[t]{12cm}
\vspace*{2cm} {\LARGE {\bf Erkl\"{a}rung}}

\bigskip
\bigskip

Hiermit versichere ich, dass ich die vorliegende Arbeit
\begin{center} \vspace*{0.5cm} \tisimple \vspace*{0.5cm} \end{center} selbst"andig verfasst habe.
Es wurden nur die angegebenen ver"offentlichten und nicht ver"offentlichten Quellen verwendet.

\bigskip
Ferner erkl"are ich, dass die vorliegende Arbeit in keinem anderen Studiengang und an keiner anderen Stelle als
Pr"ufungsleistung verwendet wurde.

\bigskip
\bigskip

Sydney, 30. Oktober 2007

\bigskip
\bigskip
\bigskip
\bigskip
\bigskip

Peter Straka \vspace*{2cm}

\end{minipage}

\clearpage

\end{document}